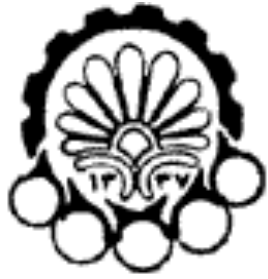

# AMIRKABIR UNIVERSITY OF TECHNOLOGY

(Tehran Polytechnic)

SCHOOL OF INDUSTRIAL ENGINEERING

MASTER THYSES

INDUSTRIAL-INDUSTRIAL ENGINEERING

## Modeling and Developing Appropriate Algorithm to Solve Generalized Probabilistic Vehicle Routing Problem

Student: Faraz Dadgostari

Assistant Professor: Mir Mehdi Seyd Esfehani

2009

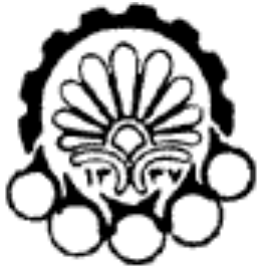

**دانشگاه صنعتی امیرکبیر**

دانشکده مهندسی صنایع

پایان‌نامه کارشناسی ارشد

مهندسی صنایع- صنایع

# ارائه الگوریتم مناسب برای حل مسأله مسیریابی توسعه‌یافته احتمالی

استاد راهنما: دکتر میرمهدی سید اصفهانی

دانشجو: فراز دادگستری

۱۳۸۸

ب

تقدیم به

خانواده‌ی عزیزم

به نام خدا

چکیده

در این رساله مسأله‌ی مسیریابی توسعه یافته در حالت احتمالی معرفی شده، مدل‌سازی گردیده و الگوریتم‌های دقیق و ابتکاری برای حل آن در ابعاد مختلف، توسعه یافته است.

در ابتدا مسأله‌ی مسیریابی کلاسیک معرفی شده، انواع معمول آن در حالت قطعی مورد بررسی قرار گرفته، مدل ریاضی آنها ارائه شده و الگوریتم‌های دقیق و ابتکاری موجود برای حل آنها شرح داده شده است. پس از شناخت مسأله کلاسیک و در حالت قطعی مسأله مسیریابی توسعه یافته احتمالی مدنظر قرار گرفته است.

نظر به اینکه این مسأله برای اولین بار در این رساله معرفی شده و نیاز به طرح مبانی تئوریک در حوضه بهینه‌سازی احتمالی مسائل عدد صحیح و در نتیجه‌ی آن مبانی جبر خطی فضاهای گسسته خواهد بود. لذا پیش از مدل سازی و ارائه الگوریتم حل دقیق، مبانی مورد نیاز برای درک مسأله و الگوریتم پیشنهادی طرح شده‌اند. مفاهیم ریاضی به کار برده شده نیز به پیوست آمده است. در مرحله بعد و با توجه به طبیعت NP-Hard مسائل مسیریابی برای حل مسائل در ابعاد بزرگ، الگوریتم‌های حل ابتکاری نیز پیشنهاد شده‌اند که در واقع توسعه‌ی الگوریتم‌های موجود در مسائل کلاسیک برای حالت توسعه یافته و احتمالی هستند. در بخش آخر این رساله نتایج محاسباتی بدست آمده از حل مسائل در ابعاد مختلف مورد بررسی و تحلیل قرار گرفته است.

**کلمات کلیدی:** بهینه‌سازی احتمالی عدد صحیح – مسیریابی توسعه یافته احتمالی – برش بهینه‌یاب – روش L-Shape

**Abstract**

This dissertation introduces stochastic generalized routing problem model and proposes exact and heuristic algorithms to solve it efficiently, in wide range of problem size.

At first, classic routing problem with its common variations in deterministic form was reviewed. Its mathematical models were demonstrated and exact and heuristic algorithms were described. After it, stochastic generalized routing was discussed.

Since this problem is introducing for the first time in this dissertation, it was necessary to dispute over theoretical principles in stochastic integer optimization and in consequence linear algebra concepts in discrete spaces. Thus before modeling the problem and developing exact and heuristic algorithms, required bases to understand the problem's proposed model and algorithms to solve it discussed. Dissertation appendix contains mathematical concepts which used in this part. In next stage with regard to NP-Hard nature of routing problem, heuristic algorithms were proposed to solve it in large scale size. Computational results from solving it in different sizes were analyzed.

**Résumé**

Cette thèse introduit stochastique a généralisé le modèle de problème routing et propose les algorithmes exacts et heuristiques pour le résoudre efficacement, dans la large gamme de grandeur de problème.

Au début, le problème routing classique avec ses variations communes dans la forme déterministe a été reconsidéré. Ce sont des modèles mathématiques ont été démontrés et les algorithmes exacts et heuristiques ont été décrits. Après cela, stochastique a généralisé routing an été discuté.

Comme ce problème introduit pour la première fois dans cette thèse, il était nécessaire de débattre sur les principes théoriques dans l'optimisation de nombre entier stochastique et dans la conséquence les concepts d'algèbre linéaires dans les espaces distincts. Ainsi avant de modeler le problème et le développement des algorithmes exacts et heuristiques, les bases exigées pour comprendre le modèle proposé du problème et les algorithmes pour le résoudre discuté. L'appendice de thèse contient des concepts mathématiques qui ont utilisé dans cette partie. Dans le stade suivant en ce qui concerne la nature NP-Hard de problème routing, on a proposé des algorithmes heuristiques pour le résoudre dans la grande grandeur d'échelle. Les résultats quantificatifs de la solution de cela dans de différentes grandeurs ont été analysés.

Mots clé: l'optimisation de nombre entier stochastique, Routing Généralisé Stochastique, la Coupe Optimale, la Méthode de L-forme

# فهرست مندرجات

## فصل دوم: بهینه‌سازی تصادفی

## فصل سوم: مسأله مسیریابی توسعه یافته احتمالی

## فصل چهارم: طراحی الگوریتم‌های مناسب (دقیق و ابتکاری) برای حل مسأله مسیریابی توسعه‌یافته احتمالی

# فهرست جداول

# فهرست اشکال

# فصل اول

# معرفی مسأله مسیریابی کلاسیک

**مقدمه**

مسأله مسیریابی وسیله نقلیه از شناخته شده‌ترین مسائل در حیطه برنامه‌ریزی حمل و نقل و توزیع است که در رابطه با مباحث مکان‌یابی یا اصولاً هر مسأله‌ای که با مسیریابی و حمل و نقل به طور سیستماتیک مرتبط بوده و نیاز به حل ریاضی دارد مورد استفاده می‌گیرد. آنچه در ادامه خواهد آمد حداقل‌های لازم برای آشنایی با مسائل VRP کلاسیک، الگوریتم‌های دقیق و تقریبی حل آنها است.

## ۱-۱- تعریف مسأله

### ۱-۱-۱- تعریف مسأله فروشنده دوره گرد (TSP)[1]

مسأله TSP یک مسأله عام در تحقیق در عملیات می‌باشد که در آن فروشنده‌ای باید به چند شهر مراجعه نموده و به آنها خدمتی ارائه نماید. هدف از مسأله تعیین بهترین مسیر و تقدم و تأخر مراجعه به شهرها می‌باشد به نحوی که دربرگیرنده حداقل مسافت (هزینه) باشد. لازم به ذکر است که فروشنده از یک شهر شروع کرده و بعد از مراجعه به سایر شهرها، به همان شهر برمی‌گردد.

---

1- Traveling salesman problem

۲

این مسأله کاربردهای زیادی دارد و تحقیقات گسترده‌ای درباره آن انجام شده است. برای درک بهتر به شکلهای (۱-۱) و (۱-۲) توجه کنید:

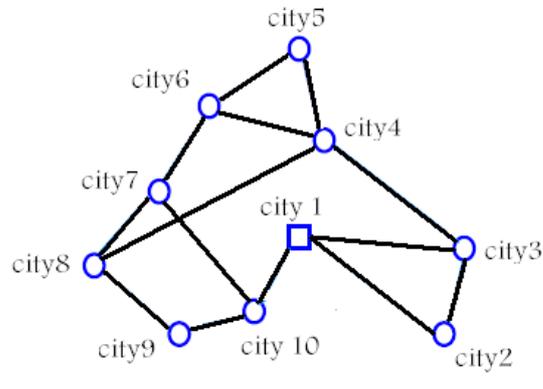 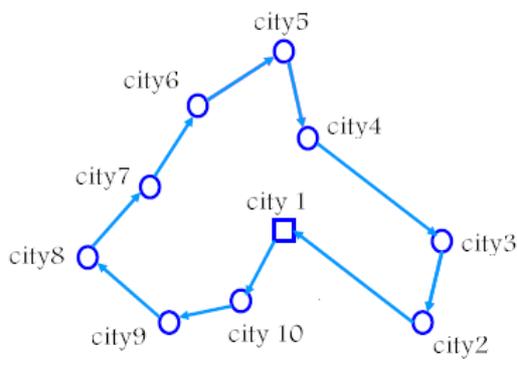

**شکل۱-۱- شبکه راههای موجود**          **شکل ۱-۲- مسیر بهینه**

**داده‌ها:**         **سوال:**

۱- مجموعه شهرها و شبکه راههای موجود

۲- هزینه سفر از شهر i به شهر j ($C_{ij}$)

۳- تعیین نقطه شروع و پایان

ارزان‌ترین (یا کوتاه‌ترین) مسیر برای طی همه شهرها فقط برای یک بار و بازگشت به شهر مبدا. به زبان تئوری گراف کوتاه‌ترین دور همیلتونی در یک گراف کامل وزنی کدام است.

### ۱-۱-۲- تعریف مسأله مسیریابی وسائل نقلیه کلاسیک (VRP کلاسیک)

این مسأله اولین بار توسط Dantzig and Ramser و در سال ۱۹۵۹ معرفی شد و یک مدل ریاضی و الگوریتم حل دقیق برای آن پیشنهاد گردید [۶۰] [۲۰۰۲ .Vigo و Toth]. هدف از این مسأله تعیین مسیر خودروهایی است که هر یک قرار است به دسته‌ای از مشتریان در یک گستره جغرافیایی مراجعه نموده وخدمتی را ارائه نمایند. لذا در VRP هم تعیین مجموعه مشتریان اختصاصی به هر خودرو مدنظر بوده و هم تعیین مسیر رفت و برگشت و تقدم و تأخر مراجعه به آنها مطرح می‌باشد. در واقع علاوه بر دسته‌بندی مشتریان بایستی مسیریابی خودروها هم انجام شود.

همچنین مجموعه‌ای از فرضها برای تعریف VRP کلاسیک وجود دارد که به این شرح است:

- **فرض ۱**: وسایل نقلیه ما از یک نوع هستند.
- **فرض ۲**: ظرفیت وسیله نقلیه نا محدود است.
- **فرض ۳**: در گره‌ها توقف نداریم.

با فرض داشتن شبکه راه‌های موجود (شکل ۱-۱) نتیجه حل یک مسأله VRP میتواند مطابق شکل ۱-۳ باشد.

| داده‌ها: | سوال: |
|---|---|
| ۱- مجموعه گره‌ها[1] و مسیرها[2] | چه مجموعه‌ای از مسیر ها به وسایل نقلیه (هر وسیله یک مسیر) اختصاص دهیم تا هزینه ارایه خدمت به مجموعه گره‌ها (مشتریان) کمینه شود. |
| ۲- هزینه سفر از شهر i به شهر j ($C_{ij}$) | |
| ۳- تعیین نقطه شروع و پایان (دپو) | |
| ۴- تعداد وسایل نقلیه در دسترس در دپو[3] | |

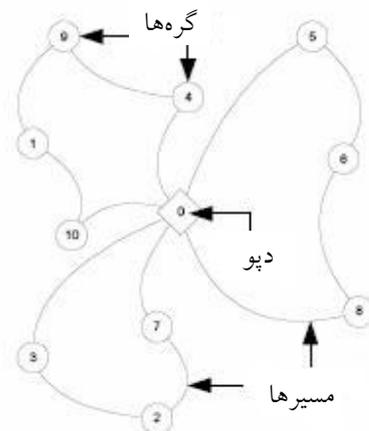

**شکل ۱-۳- با مقایسه VRP و TSP می‌توان دریافت که شالوده‌ی VRP بر TSP استوار بوده و از آن پیچیده‌تر می‌باشد و در واقع TSP حالت خاصی از VRP است. در ادامه به صورت مفصل به هر دو مسأله مذکور پرداخته می‌شود.**

---

1- Nodes
2- Routes
3- Depot

### ۳-۱-۱- معرفی مسأله VRP در حالت کلی

در منابع مطالعاتی تحقیق در عملیات، حمل و نقل و مدیریت توزیع، مسأله VRP به عنوان یک حوزه مطالعاتی مهم همواره مورد توجه پژوهشگران و مدیران بوده است. VRP بر مجموعه‌ای از مسائل که در آنها تعدادی خودرو متمرکز در یک یا چند قرارگاه (دپو) بایستی به مجموعه‌ای از مشتریان مراجعه نموده و خدمتی را ارائه دهند، مشتمل می‌گردد. مشتریان در یک گستره جغرافیایی پراکنده هستند و هر یک دارای تقاضای معینی می‌باشند. خدمت ارائه شونده ممکن است شامل دریافت و یا تحویل کالا و یا هر نوع دیگر باشد. ظرفیت خودروها هم مشخص بوده و همگی مسیرهای مربوطه را از مبدأ (قرارگاه مرکزی) شروع کرده و به آن ختم می‌کنند.

مبدأ ممکن است انبار مرکزی یا دفتر مرکزی اداره خدماتی و ... باشد [۶۰] [۲۰۰۲ .Vigo و Toth]. تابع هدف مسأله می‌تواند برآورده شدن تقاضای مشتریان و ارائه خدمات با حداقل تعداد خودرو و حداقل مسافت طی شونده، تعریف گردد.

این مسأله به مرور زمان پیچیده‌تر شده است و جنبه‌های دیگری به لحاظ تنوع خودروها، احتمالی بودن تقاضاها، تعدد مبادی شروع مسیرها و ... را دربرگرفته است.

#### ۳-۱-۱-۱- تعریف ریاضی مسأله VRP

فرض کنید $G=(V,E)$ یک نگاره (گراف)[1] باشد که در آن $V=\{v_1, v_2, ..., v_m\}$ مجموعه گره‌ها و $E=\{(v_i, v_j) \mid v_i, v_j \in V, i<j\}$ مجموعه سویه‌های (کمان‌های) نگاره می‌باشند. گره $v_1$ نشان دهنده انبار یا دفتر مرکزی (قرارگاه مرکزی) است. هر یک از گره‌های V (به غیر از $v_1$) دارای تقاضای نامنفی $q_i$ ($i=۲,۳,...,m$) می‌باشند و زمان ارائه خدمات به آنها برابر $S_i$ است.

ماتریس $C=\{C_{ij} \mid \forall i<j\}$ نشان دهنده هزینه یا زمان سفر و یا فاصله بین دو گره i و j می‌باشد که روی مجموعه E تعریف شده است. در صورتی که ماتریس C متقارن باشد به جای مجموعه سویه‌های E، می‌توان مجموعه خطوط بدون جهت E را جایگزین کرد. تعداد خودروها ممکن است ثابت باشد و یا به عنوان متغیر تصمیم‌گیری تلقی گردد. ظرفیت هر خودرو هم برابر $Q_i$ ($i=1,2,...,n$) می‌باشد.

---

1- Graph

۵

مسأله VRP عبارت از تعیین مجموعه‌ای از مسیرها برای خودروها با حداقل هزینه یا زمان می‌باشد. در خلال تعیین مسیرها بایستی محدودیت‌های مسأله ارضاء گردند تا جواب حاصل موجه[1] باشد. در زیر به تعدادی از محدودیت‌ها که ممکن است در مسائل با جنبه‌های مختلف مطرح شوند، اشاره می‌گردد:

۱) هر گره از مجموعه V (به غیر از $v_1$) دقیقاً یکبار و به وسیله یک خودرو ملاقات شده و به آن خدمتی ارائه می‌گردد.

۲) همه مسیرها از گره $v_1$ شروع و به آن ختم می‌شوند.

۳) مجموعه تقاضاهای گره‌های اختصاصی به هر خودرو نبایستی از ظرفیت آن خودرو بیشتر گردد.[2]

۴) طول هر مسیر نباید از یک مقدار از پیش تعیین شده بیشتر گردد.[3] این طول با توجه به سرعت (متوسط) خودروها قابل تبدیل به دیمانسیون زمان می‌باشد.

۵) به گره i بایستی در یک فاصله زمانی $(t_i, t_j)$ خدمتی ارائه گردد.[4]

۶) بعضی از گره‌ها (به عنوان مثال گروه i) بایستی قبل از گره‌های دیگری (به عنوان مثال گروه j) خدمات را دریافت نمایند.

۷) گره‌ای ممکن است از یک خودرو خاص و از پیش تعیین شده، خدمات دریافت کند.

محدودیت‌های ۱ و ۲ در همه مسائل VRP مطرح می‌باشند ولی بعضی یا همه از محدودیت‌های ۳ تا ۷ و یا هر محدودیت قابل تعریف دیگر، بسته به نوع مسأله به مدل اضافه می‌شوند.

### ۱-۱-۳-۲- طبقه‌بندی مسائل VRP

همه‌ی مسائل VRP شامل محدودیت‌های ۱ و ۲ که در تعریف ریاضی مسأله مورد اشاره قرار گرفته است می‌شوند. لیکن با توجه به طبیعت مسائل مختلف محدودیت‌های دیگری نیز اعمال می‌گردد و مسائل VRP را به چند دسته تقسیم‌بندی می‌کند. از آنجایی که این تقسیم‌بندی براساس محدودیت‌های مسأله انجام می‌گیرد بر شیوه‌ی مدل‌سازی و همچنین الگوریتم‌های حل هر یک تأثیر گذارده و لذا اهمیت می‌یابد. از این جهت در ادامه به دسته‌بندی کلی مسائل VRP اشاره کرده و در رابطه با اصطلاحات مرتبط توضیح مختصری ارائه می‌نماییم.

---

1- Feasible Solution

۲- مسائلی شامل این محدودیت به CVRP (ظرفیت محدود ، Capacity (C)) معروف هستند .

۳- مسائلی شامل این محدودیت به DVRP (حداکثر فاصله طی شونده محدود ، Distance (D)) معروف هستند .

4- Vehicle Routing Problems with Time Windows (VRPTW)

۶

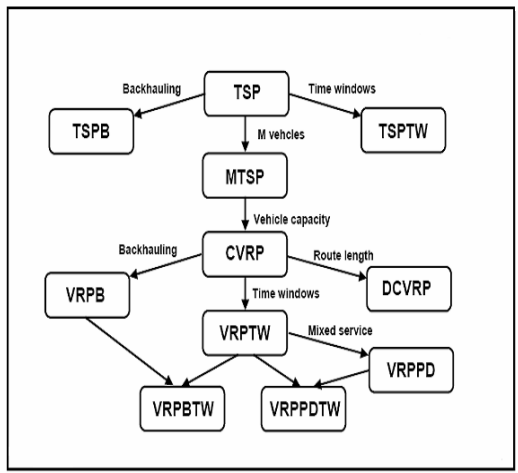

- Traveling Salesman Problem (TSP)
- Traveling Salesman Problem with Backhauls (TSPB)
- Traveling Salesman Problem with Time Windows (TSPTW)
- Multiple Traveling Salesman Problem (MTSP)
- Capacitated Vehicle Routing Problem (CVRP)
- Distance Constrained Vehicle Routing Problem (DCVRP)
- Vehicle Routing Problem with Backhauls (VRPB)
- Vehicle Routing Problem with Time Windows (VRPTW)
- Vehicle Routing Problem with Pickup and Delivery (VRPPD)

شکل ۱-۴- دسته‌بندی کلی مسائل VRP

در رابطه با این تقسیم‌بندی اصطلاحاتی که در ادامه تعریف شده‌اند برای درک مفاهیم از اهمیت بالایی برخوردارند و با اضافه شدن به اصطلاح VRP آن را خاص کرده و به دسته‌ی مشخصی از مسائل VRP مرتبط می‌نماید. این تقسیم‌بندی به طور شماتیک توسط شکل (۱-۴) نمایش داده شده است.

C (Vehicle Capacity) : ممکن است وسایل نقلیه‌ی مسئله مورد نظر محدودیت ظرفیت داشته باشند. این محدودیت می‌تواند برای وسایل مختلف متفاوت باشد.

TW (Time windows) : همانطور که پیش از این گفته شد ممکن است تمام یا بعضی از گره‌ها برای دریافت خدمات در یک (یا چند) فاصله زمانی آمادگی دارند و لذا بایستی این محدودیت در مدل‌سازی مسأله لحاظ گردد.

DC (Route length) : ممکن است در مسئله مورد نظر محدودیتی برای حداکثر طول مسیر (یا زمان طی مسیر) برای خودرو وجود داشته باشد.

B (With Backhaul) : انواع مختلفی از محدودیت‌ها را با این عنوان معرفی می‌کنند. عام‌ترین محدودیت با این عنوان به این صورت است که وسایل نقلیه ابتدا مجبور باشند تا محموله‌های خود را تحویل داده و پس تخلیه کامل شروع به تحویل گرفتن محموله نمایند.

PD (Pickup and Delivery) : در این حالت هر وسیله در هر گره هم می‌تواند کالا تحویل داده و هم بگیرد که با توجه به عرضه و تقاضای گره و ظرفیت وسیله باید محدودیت در مدل اعمال شود.

۷

### ۱-۳-۱-۳- مسیریابی گره‌ای[1] و مسیریابی سویه‌ای[2]

در یک تقسیم‌بندی کلی می‌توان مسائل مسیریابی را به دو دسته مسیریابی گره‌ای و مسیریابی سویه‌ای تقسیم نمود. تعریفی که پیشتر از مسائل مسیریابی ارائه شد مبتنی بر مسیریابی گره‌ای می‌باشد که در آنها تمام گره‌ها به وسیله خودروها خدمت داده می‌شوند. به طوری که تابع هدف در راستای محدودیت‌ها برآورده شود. مسأله فروشنده دوره گرد (TSP) نمونه‌ای از این مسیریابی است. در مسیریابی سویه‌ای، هدف پوشش‌دهی به همه سویه‌های شبکه (نگاره) با رعایت تمامی محدودیت‌ها می‌باشد. مسأله پستچی چینی (CPP)[3] از این نوع مسیریابی است.

در مسائل واقعی ممکن است که فرض نمودن هر مشتری به عنوان گره‌ای از شبکه که بایستی خدمتی را دریافت نماید، مشکل به نظر برسد. زیرا این فرض ممکن است حجم مسأله را زیاد و حل آن را دشوار نماید. در این موارد، مسیریابی سویه‌ای را می‌توان به جای مسیریابی گره‌ای به کار برد و چند مشتری را در قالب یک سویه‌ی (کمان) واحد تعریف نمود و به این ترتیب حجم مسأله را کاهش داد. در واقع مسیریابی سویه‌ای حالت پیوسته مسیریابی گره‌ای می‌باشد. در ادامه شکلهای (۱-۵) و (۱-۶) به طور شماتیک مسیریابی گره‌ای و مسیریابی سویه‌ای نشان می‌دهند.

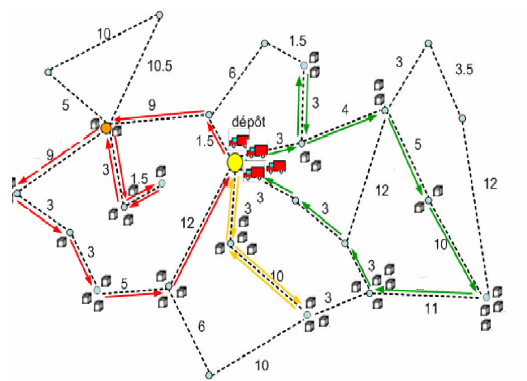 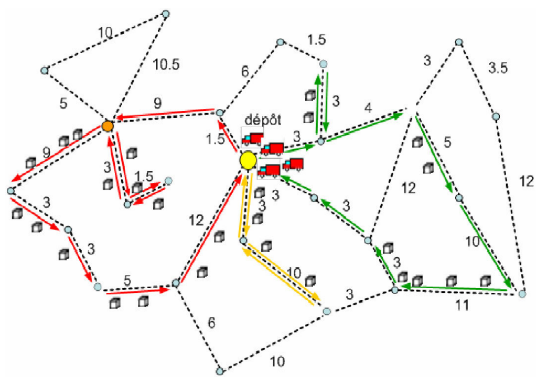

شکل ۱-۶- مسیریابی سویه‌ای         شکل ۱-۵- مسیریابی گره‌ای

---

1- Node routing problem
2- Arc routing problem
3- Chinese Postman Problem

## ۱-۲- مدل سازی ریاضی مسأله

### ۱-۲-۱- مدل ریاضی مسأله فروشنده دوره گرد

شالوده روشهای حل مسأله VRP بر مسأله TSP استوار است که مدل سازی دقیق آن به صورت برنامه‌ریزی صفر - یک[1] می‌باشد. اگر $G = [N, V, C]$ نشان‌دهنده نگاره‌ای با مجموعه گره‌های $N = \{1, 2, ..., n\}$ و مجموعه سویه‌های $V = \{(i-j) \mid i, j \in N\}$ باشد و ماتریس $C = \{C_{ij} \mid i, j \in N\}$ نشان‌دهنده هزینه سفر یا فاصله بین هر دو گره $i$ و $j$ باشد، آنگاه هدف از حل مسأله پیدا کردن مسیری دوری با حداقل هزینه یا طول می‌باشد به شرطی که از هر گره دقیقاً یکبار عبور نموده و ختم مسیر به گره‌ای باشد که از آن شروع گردیده است. این مسیر دوری را گردش (تور)[2] می‌نامند. n هم نشان‌دهنده تعداد گره‌ها است. گردشی که دربرگیرنده تمام n گره نگاره نباشد، اصطلاحاً زیر گردش[3] نامیده می‌شود. روابط زیر صورت‌بندی ریاضی (بهینه‌یاب) مسأله را بیان می‌کند.

$$\min \sum_{(i,j) \in A'} c_{ij} x_{ij}$$

S.T.

$$\sum_{i \in V \setminus \{j\}} x_{ij} = 1 \qquad \forall j \in V \setminus \{\circ\} \qquad (1-1) \text{ (محدودیت ورود به هر گره)}$$

$$\sum_{j \in V \setminus \{i\}} x_{ij} = 1 \qquad \forall i \in V \setminus \{\circ\} \qquad (1-2) \text{ (محدودیت خروج از هر گره)}$$

$$\sum_{i \in S} \sum_{j \in S} x_{ij} \geq r(S) \qquad \forall S \subseteq V \setminus \{\circ\}, S \neq \phi \qquad (1-3) \text{ (محدودیت حذف کننده زیر دورها)}$$

$$x_{ij} \in B \qquad \forall (i, j) \in A'$$

اگر سویه (i-j) در مسیر وجود داشته باشد آنگاه $x_{ij} = 1$ است و در غیر این صورت $x_{ij} = 0$ می‌باشد. روابط (۱-۱) و (۱-۲) باعث می‌شوند که مسیر ایجاد شده از هر نقطه $i \in N$ دقیقاً یکبار عبور نماید.

---

1- Binary Programming
2- Tour
3- Subtour

آنچه گفته شد در ادامه به صورت شماتیک آمده است.

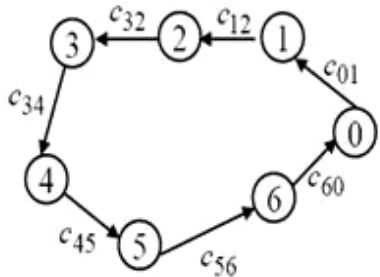

به ازای $\forall (i,j) \in A$

$$x_{ij} = \begin{cases} 1 & : \text{اگر گره } j \text{ دقیقاً بعد از گره } i \text{ طی شود} \\ 0 & : \text{در غیر اینصورت} \end{cases}$$

$c_{ij}$ : هزینه رفتن از گره $i$ به گره $j$

رابطه (۳-۱) به محدودیت زیر گردش‌ها[1] معروف می‌باشد [۶۰] [Toth و Vigo. ۲۰۰۲]. اگر مسأله بدون این محدودیت حل شود و فاقد زیر گردش باشد آنگاه جواب حاصل موجه است. در غیر این صورت بایستی محدودیت فوق الذکر لحاظ گردد. محدودیت‌های زیرگردش‌ها به صورت‌های مختلفی قابل تعریف می‌باشند به عنوان مثال می‌توان آنها را به صورت

$$\sum_{i,j \in s} x_{ij} \leq |S| - V(S) \qquad (|S| \geq 2) \qquad (۴-۱)$$

نیز تعریف نمود که در آن $S$ هر زیرمجموعه قابل قبول (با تعداد اعضا $|S|$) از $N$ می‌باشد و $V(S)$ یک حد پایین برای جواب مسأله است که به ساختار آن بستگی دارد.

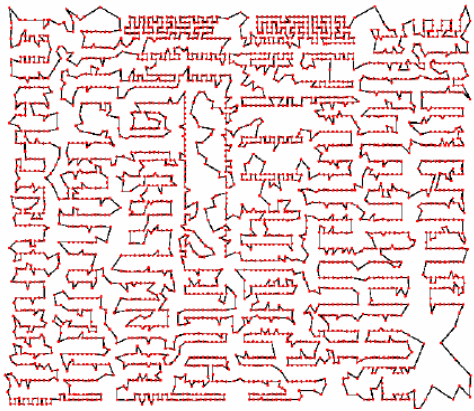

**از کاربردهای TSP می‌توان به موارد زیر اشاره کرد:**

۱- مسیریابی اتوبوس دانش آموزان

۲- مسیریابی ماشین جمع‌آوری سکه‌های تلفن عمومی

۳- جمع‌آوری بسته‌های پستی

۴- استفاده در طراحی مدارهای الکتریکی (شکل ۱-۷)

**شکل ۱-۷- طراحی مدار های الکتریکی**

---

[1] Subtours restrictions

همانطور که پیش از این نیز اشاره شد مسائل TSP و VRP جزو مسائل NP-hard بوده و حل آنها بطور دقیق در تعداد گره‌های بالا بسیار سخت است. به عنوان مثال بزرگترین TSP حل شده به صورت بهینه تاکنون در سال ۲۰۰۲ و برای یافتن جواب بهینه یک تور در کشور سوئد با ۲۴۹۷۸ گره (شهر) بوده است. برای حل این مسئله از ۹۶ عدد Intel ; dual processor Xeon 2.8 GHz در Georgia Tech's School of Industrial Engineering و به روش B&C استفاده شده است. در صورتی که فقط از یک CUP 2.8GHz استفاده میشد ۸٫۸۴ سال حل این مساله طول می‌کشید. شکل (۱-۸) حل نهایی این مسأله را نشان می‌دهد.

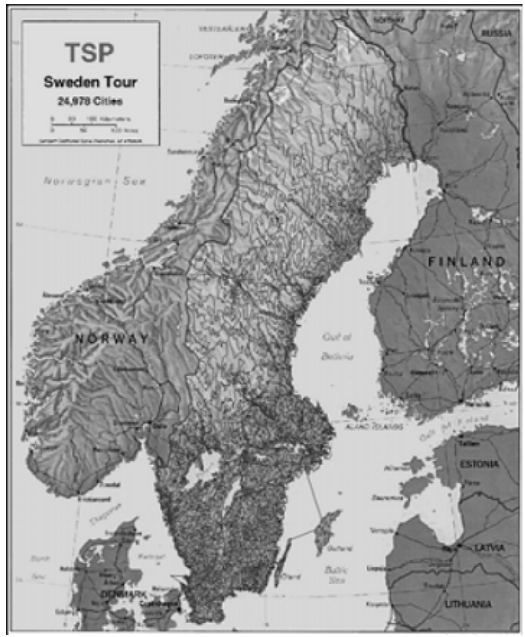

**شکل ۱-۸- بزرگترین TSP حل شده به صورت بهینه تاکنون برای یافتن جواب بهینه یک تور در کشور سوئد**

## ۱-۲-۲- مدل‌های ریاضی انواع مسائل VRP[1]

در این بخش ابتدا به مدلسازی مسئله VRP ساده می‌پردازیم و سپس با اضافه کردن محدودیت‌های جدید انواع دیگر مسائل VRP مدلسازی می‌شوند.

پارامترهای لازم برای همه‌ی مدل شامل موارد زیر است:

n : تعداد نقاط (گره‌های) تقاضا (قرارگاه مرکزی در گره ۱=i قرار دارد)

V : تعداد خودروهای در دسترس

$C_v$ : ظرفیت خودروی v ام

T : حداکثر زمانی که خودرو v می‌تواند طی مسیر کند .

$d_i$ : تقاضای گره i ( $d_1=0$ )

$t_i^v$ : زمان مورد نیاز برای ارائه خدمت به گره i به وسیله خودرو $v (t_1^v = 0)$

$t_{ij}^v$ : زمان مورد نیاز جهت طی سویه (i-j) به وسیله خودرو $v (t_{ij}^v = \infty)$

---

1- Generic model of VRP

۱۱

$C_{ij}$ : هزینه سفر (طول فاصله) بین گره‌های i و j

X : ماتریسی با درایه‌های $X_{ij} = \sum_{v=1}^{V} x_{ij}^{v}$

هدف از این مدل‌سازی حداقل کردن هزینه‌های سفر (مسافتهای طی‌شونده) برای تمامی گره‌ها و خودروها می‌باشد که در زیر آمده است:

### مدل VRP ساده:

$$\min \sum_{(i,j)\in A'} c_{ij} x_{ij}$$

S.T.

(1-5) (محدودیت ورود به هر گره) $\sum_{i\in V\setminus\{j\}} x_{ij} = 1 \qquad \forall j \in V\setminus\{\circ\}$

(1-6) (محدودیت خروج به هر گره) $\sum_{j\in V\setminus\{j\}} x_{ij} = 1 \qquad \forall i \in V\setminus\{\circ\}$

(1-7) (محدودیت تعداد وسایل وارد شده به دپو) $\sum_{i\in V\setminus\{\circ\}} x_{i\circ} = K$

(1-8) (محدودیت تعداد وسایل خارج شده از دپو) $\sum_{j\in V\setminus\{\circ\}} x_{\circ j} = K$

(1-9) (محدودیت حذف کننده زیردورها) $\sum_{i\notin S}\sum_{j\in S} x_{ij} \geq r(S) \qquad \forall S \subseteq V\setminus\{\circ\}, S \neq \phi$

$x_{ij} \in B \qquad \forall (i,j) \in A'$

روابط (1-5) و (1-6) باعث می‌شود که هر گره تقاضا فقط از یک خودرو خدمت دریافت کند. رابطه (1-8) و (1-9) متضمن این است که اگر خودرویی به گره‌ای وارد شود بایستی از آن خارج گردد و به این ترتیب پیوستگی مسیرها برقرار می‌باشد. همچنین رابطه (1-7) مربوط به حذف زیرگردش‌ها است. در ادامه با اضافه کردن محدودیت

$$\sum_{i=1}^{n} d_i (\sum_{j=1}^{n} x_{ij}^{v}) \leq c^{v} \qquad (1-10)$$

که مربوط به حداکثر ظرفیت خودروها است می‌توانیم مسأله VRP را به مسأله CVRP تبدیل کنیم. در ادامه با اضافه کردن محدودیت

$$\sum_{i=1}^{n} t_i^{v} \sum_{j=1}^{n} x_{ij}^{v} + \sum_{i}^{n}\sum_{j}^{n} t_{ij}^{v} x_{ij}^{v} \leq T \qquad (1-11)$$

که مربوط به حداکثر زمان طی مسیرها توسط هر وسیله نقلیه است به DCVRP میرسیم.

## مدل کامل DCVRP:

$$\min \sum_{(i,j) \in A'} c_{ij} x_{ij}$$

S.T.

(۱۲-۱) (محدودیت ورود به هر گره) $\quad \sum_{i \in V \setminus \{j\}} x_{ij} = 1 \qquad \forall j \in V \setminus \{\circ\}$

(۱۳-۱) (محدودیت خروج به هر گره) $\quad \sum_{j \in V \setminus \{j\}} x_{ij} = 1 \qquad \forall i \in V \setminus \{\circ\}$

(۱۴-۱) (محدودیت تعداد وسایل وارد شده به دپو) $\quad \sum_{i \in V \setminus \{\circ\}} x_{i\circ} = K$

(۱۵-۱) (محدودیت تعداد وسایل خارج شده از دپو) $\quad \sum_{j \in V \setminus \{\circ\}} x_{\circ j} = K$

(۱۶-۱) (محدودیت ظرفیت هر وسیله) $\quad \sum_{i=1}^{n} d_i (\sum_{j=1}^{n} x_{ij}^v) \leq c^v$

(۱۷-۱) (محدودیت حداکثر زمان طی مسیرها) $\quad \sum_{i=1}^{n} t_i^v \sum_{j=1}^{n} x_{ij}^v + \sum_{i}^{n} \sum_{j}^{n} t_{ij}^v x_{ij}^v \leq T$

(۱۸-۱) (محدودیت حذف کننده زیردورها) $\quad \sum_{i \notin S} \sum_{j \in S} x_{ij} \geq r(S) \qquad \forall S \subseteq V \setminus \{\circ\}, S \neq \phi$

$$x_{ij} \in B \qquad \forall (i,j) \in A'$$

## مدل VRPTW:

در ادامه با تعریف متغیرهای جدید و در نهایت اضافه کردن ۲ محدودیت به مدل DCVRP به مدل VRPTW می‌رسیم.

تعریف متغیرها و محدودیت‌های جدید برای ساختن مدل ریاضی VRPTW:

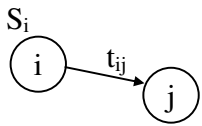

۱) $S_i$: مدت زمان لازم برای سرویس دهی در گره $i$

۲) پنجره زمانی که گره $i$ میتواند سرویس دریافت کند: $\quad a_i \leq \; \leq b_i$

۳) $w_j^v$: زمان شروع سرویس رهی در گره $i$ برای خودرو $v$ $\quad a_i \sum_{j \in v} x_{ij}^v \leq w_i^k \leq b_i \sum_{j \in v} x_{ij}^v$

۴) $E$: زودترین زمان شروع سرویس‌دهی در دپو $\quad E \leq w_i^v \leq L$

$L$: دیرترین زمان شروع سرویس‌دهی در دپو

## 1-3- الگوریتم‌های حل CVRP

در این بخش به ارائه الگوریتم‌های حل مسأله VRP و در حالت خاص CVRP می‌پردازیم. در برخی الگوریتم‌ها، کاربرد آنها را در انواع دیگر مسأله VRP بیان می‌کنیم. اما به هر حال چون CVRP عمومی‌ترین نوع این مسایل می‌باشد و اغلب الگوریتم‌ها ابتدا برای حل این مسئله به کار گرفته شده، مبنای بحثمان را بر آن قرار می‌دهیم.

اصولاً سه نوع الگوریتم برای حل مسائل بهینه‌سازی در نظر گرفته می‌شود [۶۰] [Toth و Vigo. ۲۰۰۲]:

(۱) الگوریتم‌های حل دقیق
(۲) الگوریتم‌های هیوریستیک کلاسیک[1]
(۳) الگوریتم‌های متاهیوریستیک[2]

در بخشهای آتی انواع مختلف هر یک از این الگوریتم‌ها معرفی می‌شود.

## 1-3-1- الگوریتم‌های حل دقیق

از انواع الگوریتم‌های حل دقیق مسأله CVRP می‌توان به الگوریتم‌های شاخه و کران[3]، شاخه و برش[4] و الگوریتم‌های بر مبنای پوشش کامل[5] اشاره کرد. به هر کدام از این الگوریتم در ذیل اشاره‌ای می‌کنیم.

### 1-3-1-1- الگوریتم‌های شاخه و کران (B&B) برای CVRP

متد B&B در دهه‌های اخیر برای حل مسئله CVRP و انواع اصلی دیگر آن به طور گسترده‌ای به کار گرفته شده است. Laporte و Nobert در یک مطالعه‌ی جامع، یک تحلیل کامل از الگوریتم‌های شاخه و کران ارائه شده برای CVRP ارائه کرده‌اند [۳۴] [Teixeira; Pais Antunes; Pinho de Sousa. 2004]. برای یادآوری، می‌دانیم که مسئله CVRP، گسترش یافته‌ی مسئله TSP می‌باشد که مسأله TSP به دنبال یک دور همیلتونی با کمینه کردن هزینه ملاقات‌های یک سری محل می‌باشد. بنابراین بسیاری از رویکردهای دقیق برای CVRP از راه‌حل‌های موفق TSP و گسترش آنها بدست آمده است. اخیراً هم الگوریتم‌های

---

[1] Classical Heuristics
[2] Metaheuristics
[3] Branch & Bound
[4] Branch & Cut
[5] Set-covering-based

۱٤

B&B پیچیده‌تری بر اساس ساده‌سازی لاگرانژ[1] و رویکرد فزاینده[2] ارائه شده‌اند که هر چند ابعاد مسئله را افزایش می‌دهند ولی با این کار مسایل بوسیله‌ی B&B به طور دقیق حل می‌شود. به هر حال حل دقیق مسایل CVRP در ابعاد بزرگ بوسیله‌ی الگوریتم‌های B&B ممکن نیست و حل مسایل با این روش‌ها اصولاً به شدت زمان‌بر است. اما به هر حال این مزیت که مسئله را به طور دقیق حل می‌کند، را داراست. از نتایج خوب بدست آمده از این روش، حل مسئله CVRP با ۱۰۰ گره در سال ۱۹۹۴ توسط Fisher بوده است.

### ۱-۳-۱-۲- الگوریتم شاخه و برش (B&C) برای CVRP

این روش در یافتن جواب‌های بهینه مسائل نزدیک به CVRP مانند STSP[3] بسیار موفق عمل کرده است. اما مطالعات کمی نسبت به STSP سعی در حل CVRP با این روش کرده‌اند، بنابراین نباید انتظار جواب‌های خیره‌کننده روش B&C در CVRP، مانند STSP را داشته باشیم. تعداد تحقیقات انجام شده بر روی الگوریتم B&C برای CVRP هنوز بسیار محدود است و بیشتر آنها به چاپ نرسیده (تا سال ۲۰۰۲). اکنون بسیاری از گروه‌های تحقیقاتی روی این موضوع در حال کارند و البته نتایجشان هنوز در دست نیست.

استفاده از روش B&C برای حل CVRP در اول راه توسعه‌اش می‌باشد. درک بهتر و تلاش بیشتر در روش‌های تفکیک مؤثر می‌تواند نتایج محاسباتی بهتری تولید کند.

### ۱-۳-۱-۳- الگوریتم‌ها بر مبنای پوشش کامل برای CVRP

لیست مجموعه‌ای از n مشتری که با $V = \{1, 2, ..., n\}$ نشان داده می‌شود را در نظر بگیرید. « ۰ » را برای نشان دادن دپو استفاده می‌کنیم و نقاط گراف کامل مربوطه مان به صورت $V^{(0)} = \{0\} \cup V$ می‌شود. به مشتری $i \in V$ تقاضای $d_i > 0$ را نسبت می‌دهیم که نشان می‌دهد که در محل i بار باید برداشته شود. C ظرفیت وسیله نقلیه می‌باشد و فرض می‌کنیم که برای عمل تحویل k وسیله در اختیار داریم. به وضوح، برای شدنی بودن لازم است که برای هر $i \in V$، $d_i \leq C$ باشد. $t_{ij}$ را طول یال $(i, j)$، $i, j \in V$ قرار می‌دهیم و فرض می‌کنیم که $t_{ij}$ها نامساوی مثلث را ارضا می‌کنند. بنابراین می‌توان یک

---
1- Lagrangian Relaxation
2- additive approach
3- Symmetric Traveling Salesman Problem

عدد یا هزینه ثابت را به همه یال‌ها یا گره‌ها اضافه کنیم. یک روش اولین بار توسط Quandt و Balinski برای حل مسئله CVRP بر اساس فرمولاسیون مسئله مانند یک مسئله پوشش کامل پیشنهاد شد. ایده به این صورت بود که:

تمام مسیرهای شدنی شمرده شوند. یک مسیر شدنی عبارتست از مسیری که از دپو آغاز شود و به آن ختم شود و مجموع بار برداشته شده از ظرفیت C بیشتر نشود. مجموعه مسیرهای شدنی را $\Re = \{1, 2, ..., R\}$ قرار می‌دهیم. $C_r$ هزینه (یا طول) مسیر r و $S_r$ که زیر مجموعه‌ای از V می‌باشد، نشان دهنده مشتریانی است که در مسیر r ظاهر شده‌اند. تعریف می‌کنیم:

$$\alpha_{ir} = \begin{cases} 1 & \text{اگر مشتری i در مسیر r سرویس بگیرد} \\ 0 & \text{در غیر این صورت} \end{cases}$$

$$y_r = \begin{cases} 1 & \text{در مسیر بهینه باشد} \\ 0 & \text{در غیر این صورت} \end{cases}$$

در فرمولاسیون بر مبنای پوشش کامل مسئله CVRP هدف انتخاب یک مجموعه از مسیرهای شدنی با کمینه کردن هزینه (یا مسافت) می‌باشد که هر مشتری باید در یک مسیر شامل شده باشد. بیان این قضیه به صورت ریاضی به صورت زیر است:

$$\text{Min} \sum_{r \in R} C_r y_r$$

S.T. 
$$\sum_{r \in R} a_{ir} y_r \geq 1 \quad (1-19)$$

$$\sum_{r \in R} y_r \leq k \quad (1-20)$$

$$y_r \in \{0, 1\} \quad \forall r \in R$$

$$a_{ir} \in \{0, 1\} \quad \forall r \in R, \forall i \in V$$

محدودیت اول می‌گوید که هر مشتری باید حداقل در یک مسیر شدنی ظاهر شود و محدودیت دوم می‌گوید حداکثر k مسیر شدنی می‌تواند مورد استفاده قرار گیرد (چون کلاً k خودرو در اختیار داریم و هر خودرو فقط یک مسیر را می‌تواند طی کند). از لحاظ منطقی به نظر می‌رسد که محدودیت اول باید به صورت تساوی برقرار باشد. در واقع فرمولاسیون با محدودیت‌های مساوی هم معادل همین فرمولاسیون

۱٦

بالاست چرا که فرض کرده‌ایم که ماتریس {$t_{ij}$} نامساوی مثلث را ارضا میکند، پس بنابراین در جواب بهینه هر مشتری فقط یک بار ملاقات می‌شود (ملاقات‌های بیشتر هزینه بیشتر به همراه دارد). بدین علت اینجا از محدودیت‌های نامساوی استفاده کردیم چون کار کردن با آنها راحت‌تر است. بر مبنای این فرمولاسیون برنامه‌ریزی ریاضی الگوریتم‌های حل دقیق و هیوریستیک متعددی ایجاد شده‌اند. این فرمولاسیون به طور موفق توسط Ratliff، Jarvis و Culleu برای طراحی روشهای هیوریستیک برای VRP استفاده شده است. الگوریتم‌های دقیقی هم بر این مبنا توسط Agarwal، Mathur و Salkin و اخیراً توسط Bixby و Hadji Constantinon, Christofides & Mingozzi ایجاد و توسعه یافته‌اند. [۶۰] این روش بسیار کلی است و می‌تواند در بسیاری از مسائل VRP مانند VRPTW (البته بدون محدودیت ظرفیت) و یا CVRPTW مورد استفاده قرار گیرد. الگوریتم عمومی برای روش‌های بر مبنای پوشش کامل مراحلی دارد که در مرحله اول یک مسأله آزاد شده[1] از این فرمولاسیون را در نظر می‌گیرند و از حل آن یک حد پایین برای جواب بهینه بدست می‌آورند. در جداول ۱-۱ و ۱-۲، کد مسئله، تعداد مشتریان (n)، مقدار حد پایین ($Z^{LB}$) و مقدار حد بالا ($Z^{UB}$) با متدهای مختلف بدست آمده‌اند. در بیشتر این الگوریتم‌ها شاید به خود جواب بهینه نرسیم ولی حد بالا و پایین بدست آمده تا حد زیای به هم نزدیک می‌شوند و بنابراین به مقدار بهینه بسیار نزدیک می‌شویم. مؤثر بودن[2] حد پایین به صورت زیر تعریف می‌شود: ($Z^{LB}/Z^{UB}$)۱۰۰ یا ($Z^{LB}/Z^{*}$)۱۰۰ بسته به اینکه جواب بهینه مسئله را داشته باشیم یا نه تعریف می‌شود. همان‌طور که ملاحظه می‌شود تقریباً در همه موارد یک چیز مشترک است و آن این است که حد پایین حاصل از آزادسازی فرمولاسیون پوشش کامل مبنا بسیار قوی می‌باشد (یعنی بسیار به حد بالا و در نتیجه جواب بهینه نزدیک است).

---

1- relaxed
2- effectiveness

۱۷

جدول ۱-۱- نتایج الگوریتم پوشش کامل مبنای Agarwal، Mathur و Salkin

| مسئله | n (تعداد مشتریان) | $Z^{LB}$ | $Z^{UB}$ | مؤثربودن حد پایین |
|---|---|---|---|---|
| E۰۱۶-۰۳M | ۱۵ | ۲۶۸ | ۲۷۶ | ۹۷/۱٪ |
| E۰۱۶-۰۵M | ۱۵ | ۳۲۶ | ۳۳۲ | ۹۸/۲٪ |
| E۰۲۱-۰۴M | ۲۰ | ۳۵۱ | ۳۵۸ | ۹۸/۰٪ |
| E۰۲۱-۰۶M | ۲۰ | ۴۳۰ | ۴۳۰ | ۱۰۰/۰٪ |
| E۰۲۲-۰۴M | ۲۱ | ۳۷۴ | ۳۷۵ | ۹۹/۷٪ |
| E۰۲۲-۰۶M | ۲۱ | ۴۷۹ | ۴۹۴ | ۹۷/۰٪ |
| E۰۲۶-۰۸M | ۲۵ | ۶۰۶ | ۶۰۷ | ۹۹/۸٪ |

جدول ۱-۲- نتایج الگوریتم پوشش کامل مبنای Bixby

| مسئله | n (تعداد مشتریان) | $Z^{LB}$ | $Z^{UB}$ | مؤثربودن حد پایین |
|---|---|---|---|---|
| S۰۰۷-۰۲a | ۶ | ۱۴۴ | ۱۱۴ | ۱۰۰/۰٪ |
| S۰۱۳-۰۴d | ۱۲ | ۲۷۹ | ۲۹۰ | ۹۶/۲٪ |
| E۰۲۱-۰۶m | ۲۰ | ۴۳۰ | ۴۳۰ | ۱۰۰/۰٪ |
| E۰۲۲-۰۴g | ۲۱ | ۳۷۵ | ۳۷۵ | ۱۰۰/۰٪ |
| E۰۲۳-۰۳g | ۲۲ | ۵۶۶ | ۵۶۹ | ۹۹/۵٪ |
| E۰۳۰-۰۴s | ۲۹ | ۵۰۳ | ۵۰۳ | ۱۰۰/۰٪ |
| E۰۵۱-۰۵e | ۵۰ | ۵۱۸ | ۵۲۱ | ۹۹/۴٪ |

## ۱-۳-۲- هیوریستیک‌های کلاسیک برای CVRP

خانواده‌های متعددی از روش‌های ابتکاری یا خلاق برای VRP پیشنهاد شده است که عمدتاً می‌توان به دو کلاس اصلی دسته‌بندی‌شان کرد:

(۱) هیوریستیک‌های کلاسیک که بیشتر در سالهای ۱۹۶۰ تا ۱۹۹۰ به وجود آمده‌اند و

(۲) متاهیوریستیک‌ها که رشدشان در دهه‌های اخیر بوده است. روش‌های هیوریستیک، عموماً شناسایی محدود و جستجوی متناسبی از فضای حل را انجام می‌دهند و نوعاً جواب‌های خوبی در زمان‌های محاسباتی نسبتاً کم تولید می‌کند. بعلاوه، اکثر آنها را می‌توان به سادگی برای محاسبه‌ی محدودیت‌های مختلف مسایل دنیای واقعی، گسترش داد. بنابراین این دسته هنوز هم به طور گسترده‌ای در بسته‌های تجاری استفاده می‌شود. در متاهیوریستیک‌ها، تأکید روی انجام جستجوی عمیق‌تر در مناطقی از فضای حل می‌باشد که احتمال یافتن جواب بهینه در آنجا بیشتر است. این متدها نوعاً قوانین پیچیده‌ای برای ترکیب جواب‌ها با هم و بدست آوردن یک جواب بهتر دارند. بعلاوه این متدها بسته به مسأله مورد نظر احتیاج به طراحی و تنظیم پارامترها دارند و همین موضوع گسترش آنها را در کاربرد برای موارد دیگر مشکل می‌سازد. در بخش‌های آتی در مورد متاهیوریستیک‌ها صحبت خواهیم کرد.

هیوریستیک‌های کلاسیک برای VRP عمدتاً به سه دسته تقسیم می‌شوند:

۱- **هیوریستیک‌های سازنده**[1]: یک حل شدنی را می‌سازند و در حین ساختن نیم نگاهی هم به هزینه حل دارد. اما فی‌نفسه فاز بهبود را شامل نمی‌شود. روش‌های سازنده دو تکنیک کلی در ساختن جواب‌های VRP دارند: (۱) ادغام کردن مسیرهای موجود با استفاده از معیار صرفه‌جویی (مانند الگوریتم صرفه‌جویی Clarke[2] و Wright) و (۲) تخصیص تدریجی گره‌ها به مسیرهای وسایل نقلیه با استفاده از هزینه اضافه کردن[3]

۲- **هیوریستیک‌های دو مرحله‌ای**[4] : مسئله را به دو جزء طبیعی خود منفک می‌کند: (۱) دسته کردن گره‌ها به مسیرهای شدنی (یعنی مجموع تقاضایشان بیشتر از C نشود) و (۲) ساختن مسیر واقعی و در عین حال امکان بازخورد بین دو مرحله را فراهم می‌کند.

---

1- constructive
2- saving
3- insertion cost
4- two-phase

۱۹

هیوریستیک‌های دو فازی خودشان به دو دسته تقسیم می‌شوند:

**الف- اول دسته کردن- دوم مسیریابی[1]**: اول گره‌ها را به دسته‌های شدنی ($\sum d_i \leq C$) تقسیم می‌کند و بعد یک مسیر برای هر کدام از دسته‌ها ساخته می‌شود (مانند الگوریتم جارو[2])

**ب- اول مسیریابی- دوم دسته کردن[3]**: ابتدا یک تور (مسیری که از همه‌ی گره‌ها فقط یک بار عبور کند و به گره اولش باز گردد) روی همه گره‌ها ساخته می‌شود، سپس به مسیرهای شدنی در راستای تور اولیه شکسته می‌شود. مسیر اولیه که از همه گره‌ها می‌گذرد لزوماً شدنی نیست و فقط برای شروع کار از آن استفاده می‌کنیم.

**۳- هیوریستیک‌های بهبوددهنده**: سعی می‌کنند بوسیله تغییر در ترتیب یال‌ها یا گره‌ها در یک مسیر یا مسیرهای مختلف با هم، حل شدنی‌مان را بهبود دهند.

حال به بررسی دو الگوریتم هیوریستیک معروف و پراستفاده در VRP می‌پردازیم.

### ۱-۳-۲-۱- الگوریتم صرفه‌جویی Clark و Wright

این الگوریتم برای مسیریابی گره‌ای توسط Clarke و Wright ارائه گردید. که بر استراتژی صرفه‌جویی/ اضافه کردن، استوار است. در این الگوریتم ابتدا فرض می‌شود که هر گره تقاضا به همراه قرارگاه مرکزی در یک مسیر قرار می‌گیرند. سپس بر مبنای صرفه‌جویی‌های ممکنه و رعایت محدودیت‌های مسأله (حداکثر ظرفیت خودرو یا حداکثر طول هر مسیر) مسیرها با هم ادغام می‌شوند. اگر شبکه دارای n گره بوده و قرارگاه مرکزی در گره ۱ باشد، مراحل این الگوریتم به شرح زیر است:

قدم ۱: برای گره‌های $i = 2,...,n$، $n-1$ مسیر جداگانه که هر کدام شامل قرارگاه مرکزی می‌باشند، به صورت (۱-i-۱) ایجاد نمائید.

قدم ۲: برای هر زوج گره j و i ، مقدار صرفه جویی را به صورت زیر محاسبه نمائید:

$$S_{ij} = C_{i1} + C_{j1} - C_{ij}, (i,j = 2,...,n, i \neq j) \qquad (2-21)$$

که در آن $C_{ij}$ مسافت و یا هزینه مرتبط با گره‌های i و j می‌باشد و رابطه $C_{ij} = C_{ji}$ برقرار است.

قدم ۳: $S_{ij}$ را به صورت نزولی مرتب کنید.

---

1- Cluster first - Route second
2- Sweep Algorithm
3- Route first - Cluster second

قدم ۴: از بیشترین مقدار $S_{ij}$ شروع کرده و برای هر i و j که گره‌های ابتدایی یا انتهایی مسیرهای جداگانه می‌باشند. چنانچه $S_{ij} > 0$ است آنگاه مسیرهای مرتبط با آنها را موقتاً با هم ادغام کرده و به عنوان مسیر واحد معرفی نمائید. چنانچه مسیر اخیر محدودیت‌های مسأله را ارضاء کرده و زیر گردشی ایجاد نمی‌گردد، آن را بپذیرد. در غیر اینصورت ادغامی صورت نگیرد.

این مرحله را تا جایی که تمامی بهبودهای ممکنه انجام گردد، ادامه دهید.

قدم ۵: توقف

قدم ۴ الگوریتم به فرض نامحدود بودن تعداد خودروها عنوان شده است. اما اگر تعداد خودروها محدود باشد می‌توان مرحله‌ی یاد شده را تا رسیدن به حداکثر تعداد خودرو ادامه داد.

به طور مشابه ساختار این الگوریتم می‌تواند برای مسیریابی سویه‌ای بکار رود.

### ۱-۳-۲-۲- الگوریتم جارو

این الگوریتم در سال ۱۹۷۲ مطرح و در سال ۱۹۷۴ توسط Gillett و Miller بسط داده شد. [۶۰] در این الگوریتم ابتدا مجموعه گره‌های اختصاصی به هر خودرو مشخص شده و بعد از آن مسیر مناسب و تقدم و تأخر ملاقات گره‌ها به وسیله خودروها، مشخص می‌گردد.

اطلاعات لازم برای این الگوریتم شامل مختصات و تقاضای گره‌ها و حداکثر خودروهای در دسترس می‌باشد. مختصات قطبی گره‌ها از روی مختصات دکارتی قابل محاسبه است. برای این کار یکی از گره‌ها که معمولاً قرارگاه مرکزی است را در مبدأ مختصات دکارتی فرض کرده و زوایای قطبی گره‌ها نسبت به آن سنجیده می‌شود قدم‌های الگوریتم به شرح زیر می‌باشد:

قدم ۱: مختصات قطبی گره‌ها را محاسبه و آنها را برحسب زاویه قطبی و به صورت صعودی مرتب نمائید.

قدم ۲: خودرویی که هنوز استفاده نشده است، انتخاب نمائید.

قدم ۳: از بین نقاطی که هنوز به هیچ خودرویی اختصاص نیافته‌اند، نقطه‌ای که دارای کمترین زاویه قطبی است را به خودرو انتخابی در قدم ۲ اختصاص دهید. به همین ترتیب سایر نقاط اختصاص نیافته را تا جایی که محدودیت‌های مسأله اجازه میدهد به مجموعه نقاط اختصاص یافته به این خودرو اضافه نمایید. (محدودیت‌ها شامل حداکثر ظرفیت هر خودرو، حداکثر طول هر مسیر و ... می‌باشد.)

قدم ۴: چنانچه نقاط تخصیص نیافته‌ای وجود دارند به قدم ۲ بروید در غیر این صورت به قدم ۵ بروید.

قدم ۵: مبدأ مختصات (قرارگاه مرکزی) را جزء مجموعه نقاط اختصاصی به هر خودرو لحاظ نمائید و برای تمام مجموعه‌ها بهترین مسیر را با بکارگیری روش‌های مسألهٔ TSP بدست آورید. در صورتی که جابجایی بعضی از نقاط در مسیرها باعث کاهش مسافت کل می‌شود، تغییرات را انجام داده و مسیرهای خودروها را مجدداً تعیین نمائید. سپس محورهای مختصات برخلاف جهت عقربه‌های ساعت چرخش داده می‌شوند. و در هر مرحله چرخش، قدم‌های ۱ تا ۵ مجدداً انجام می‌گیرند تا جواب جدیدی برای مسأله حاصل شود. فرآیند چرخش محورها و حصول جواب‌های مختلف تا جایی ادامه پیدا می‌کند که هر گره شبکه به علت چرخش محورها، دارای کمترین زاویه قطبی می‌گردد.

این مرحله به مرحله‌ی حرکت رو به جلو[1] معروف است که بعد از پایان آن مرحله‌ی حرکت رو به عقب انجام می‌شود[2]. در مرحله اخیر و در قدم (۱) الگوریتم، زوایای قطبی گره‌ها به صورت نزولی مرتب می‌شوند. در این مرحله هم فرآیند چرخش محورهای مختصات (در جهت حرکت عقربه‌های ساعت) و ایجاد جواب‌های متعدد صورت می‌گیرد تا نهایتاً جواب نهایی مسأله از بین تمام جواب‌های تولیدی در مرحله‌ی رو به جلو و رو به عقب حاصل آید. طرح الگوریتم صرفه‌جویی و جارو در این قسمت به خاطر مقایسه‌ای است که در فصل بعدی بین آن‌ها و الگوریتم‌های پیشنهادی به عمل می‌آید. برای آشنایی با الگوریتم‌های ابتکاری دیگر می‌توان به مراجع مربوطه رجوع کرد.

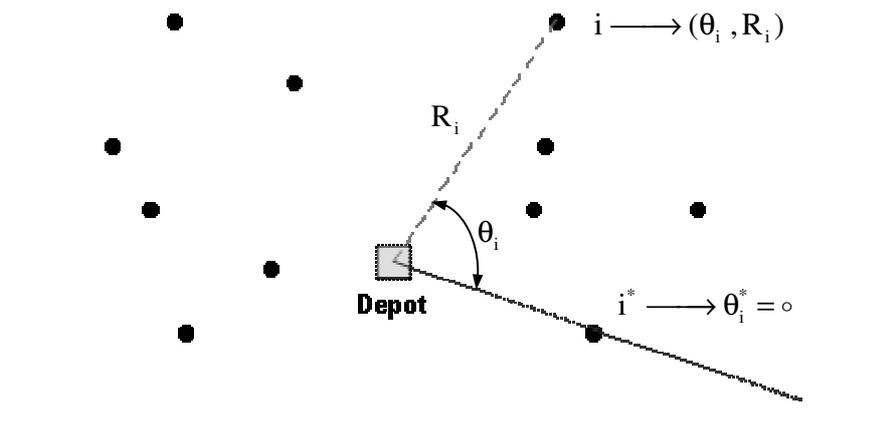

شکل ۱-۱۰- نمایش تصویری الگوریتم جارو

---

1- Forward
2- Backward

۲۲

### ۱-۳-۲-۳- نتیجه‌گیری

بیش از ۴۰ سال از چاپ روش هیوریستیک برای VRP می‌گذرد و در این مدت انواع مختلفی از تکنیک‌های حل پیشنهاد شده‌اند. مقایسه بین این روش‌ها همیشه آسان نیست، خصوصاً چون خصوصیات اجرایی مختلف می‌تواند روی عملکرد الگوریتم اثر بگذارد. در ضمن، تعداد و اندازه مسائل تست برای مقایسه محدود بودند و محققان هم از قراردادهای یکسانی در روند کار تحقیقاتشان پیروی نکرده‌اند. در سال‌های اخیر این مشکل برطرف شده است. امروزه روشن است که از نظر کیفیت جواب، هیوریستیک‌های کلاسیک که بر مبنای ساخت‌های ساده و تکنیک‌های بهبود موضعی هستند، قابل مقایسه با الگوریتم‌های جستجوی ممنوع[1] در VRP نیستند. البته این روش‌ها به سادگی قابل انطباق بر سایر مدل‌های VRP هستند و اجرای آنها ساده است و در ضمن در زمان کوتاه‌تری به یک جواب نسبتاً خوب می‌رسند. این موضوع استفاده گسترده از این روش‌ها را در نرم‌افزارهای تجاری توضیح می‌دهد. برای مثال الگوریتم جاروی Clarke و Wright در عمل هنوز پرکاربردترین روش در عامه می‌باشد.

اما به هر حال، به علت اینکه متاهیوریستیک‌ها در CVRP از نظر کیفیت جواب از روش‌های هیوریستیک کلاسیک (و در بعضی موارد حتی در زمان محاسبات) پیشی گرفته‌اند، به نظر می‌رسد که اهمیت کمی برای کار بر روی بهبود روش‌های هیوریستیک کلاسیک احساس می‌شود.

---

[1] - Tabu search

### 1-3-3- متاهیوریستیک‌ها برای CVRP

در سال‌های اخیر متاهیوریستیک‌های بسیاری برای VRP پیشنهاد شده است. متاهیوریستیک‌ها رویکردهای حل کلی هستند که فضای حل را برای تشخیص جواب‌های خوب می‌گردند و اغلب ساختن مسیرشدنی و بهبودهای هیوریستیک را در درون خود دارند. اختلاف متاهیوریستیک‌ها و هیوریستیک‌های کلاسیک در این است که متاهیوریستیک‌ها، جواب‌های بدتر و حتی جواب‌های میانجی نشدنی را در مسیر فرآیند جستجو مجاز می‌دانند. مشهورترین متاهیوریستیک‌ها که برای VRP جواب بهتری نسبت به آخرین هیوریستیک‌ها می‌دهند، اما اصولاً متاهیوریستیک‌ها زمان برتر هستند.

شش نوع اصلی از متاهیوریستیک‌ها برای VRP مورد استفاده قرار گرفته است:

(1) Simulated Annealing (SA)

(2) Deterministic Annealing (DA)

(3) Tabu Search (TS)

(4) Genetic Algorithm (GA)

(5) Ant Systems (AS)

(6) Neural Network (NN)

سه نوع اول از یک جواب اولیه $x_1$ شروع می‌کنند و در هر تکرار $t$ از $x_t$ به یک جواب $x_{t+1}$ در همسایگی $x_t$ حرکت می‌کند که در $x_{t+1}$ یک شرط توقف ارضا می‌شود. اگر $f(x)$ هزینه جواب در $x$ باشد لزوماً $f(x_{t+1})$ از $f(x_t)$ کمتر نمی‌باشد. در ضمن باید الگوریتم طوری طراحی شود که از افتادن در حلقه‌ی تکرار اجتناب کند. GA در هر مرحله، جامعه‌ای[1] از جواب‌ها را تست می‌کند. هر جامعه از ترکیب بهترین المان‌ها و دور انداختن بدترین المان‌های جامعه قبلی‌اش حاصل می‌شود. SA یک رویکرد سازنده است که در هر تکرار، با استفاده از برخی اطلاعات تکرار قبلی جواب‌های متعددی ساخته می‌شود. همان طور که Taillard می‌گوید، TS، GA و AS روش‌هایی هستند که اطلاعات جواب‌های بدست آمده را بایگانی می‌کنند و از آن برای بدست آوردن جواب بهتر استفاده می‌کنند. [60] NN یک مکانیزم یادگیرنده[2] دارد که به تدریج یک مجموعه از وزن‌ها را تعدیل می‌کند تا به یک جواب قابل قبول برسد. قوانین حاکم بر جستجو در هر یک از این موارد متفاوت است و هم چنین هر متاهیوریستیک باید به تناسب مسأله در دست طراحی شود. این کار احتیاج به تجربه و خلاقیت دارد.

---
1- population
2- learning

با توجه به نتایج بدست آمده، ملاحظه می‌شود که متاهیوریستیک SA زمان انجام محاسباتش نسبت به هیوریستیک‌ها خیلی بالاتر نیست و البته جواب‌هایش به نسبت بهتر هستند (جدول ۱-۳). متاهیوریستیک AS هم زمان کمتری نسبت به بقیه متاهیوریستیک‌ها دارد اما جواب‌های خوبی بدست نمی‌دهد (جدول ۱- ۴). در مورد TS، زمان‌ها به نظر طولانی‌تر هستند اما جواب‌های بسیار خوبی بدست آمده است. در واقع جواب‌های TS از همه متاهیوریستیک‌های دیگر بهتر است (جدول ۱-۵).

جدول ۱-۳- نتایج محاسباتی الگوریتم SA، Osman

| زمان | بهترین جواب شناخته شده | $f^*$ | کد مسئله |
|---|---|---|---|
| ۱۶۷/۴ | ۵۲۴/۶۱ | ۵۲۸ | E۰۵۱-۰۵e |
| ۶۴۳۴/۳ | ۸۳۵/۲۶ | ۸۳۸/۶۲ | E۰۷۶-۱۰e |
| ۹۳۳۴/۰ | ۸۲۶/۱۴ | ۸۲۹/۱۸ | E۱۰۱-۰۸e |
| ۶۳۲/۰ | ۸۱۹/۵۶ | ۸۲۶ | E۱۰۱-۱۰c |
| ۳۱۵/۸ | ۱۰۴۲/۱۱ | ۱۱۷۶ | E۱۲۱-۰۷c |
| ۵۰۱۲/۳ | ۱۰۲۸/۴۲ | ۱۰۵۸ | E۱۵۱-۱۲c |
| ۲۳۱۸/۱ | ۱۲۹۱/۴۵ | ۱۳۷۸ | E۲۰۰-۱۷c |
| ۳۴۱۰/۲ | ۵۵۵/۴۳ | ۵۵/۴۳ | E۰۵۱-۰۶c |
| ۶۲۶/۵ | ۹۰۹/۶۸ | ۹۰۹/۶۸ | E۰۷۶-۱۱c |
| ۹۵۷/۲ | ۸۶۵/۹۴ | ۸۶۶/۷۵ | E۱۰۱-۰۹c |
| ۳۰۵/۲ | ۸۶۶/۳۷ | ۸۹۰ | E۱۰۱-۱۱c |
| ۷۶۲۲/۵ | ۱۵۴۱/۱۴ | ۱۵۴۵/۹۸ | E۱۲۱-۱۱c |
| ۸۴۳۰۱/۲ | ۱۱۶۲/۵۵ | ۱۱۶۴/۱۲ | E۱۵۱-۱۴c |
| ۵۷۰۸/۰ | ۱۳۹۵/۸۵ | ۱۴۱۷/۸۵ | E۲۰۰-۱۸c |

### ۱-۳-۳-۱- نتیجه‌گیری

بررسی متاهیوریستیک‌ها برای VRP نشان می‌دهد که بهترین این روش‌ها می‌تواند جواب‌های عالی و گاهی اوقات بهینه را برای مثال‌هایی با چند صد مشتری پیدا کند، البته با صرف یک زمان قابل توجه. TS امروزه به عنوان مؤثرترین روش در حل VRP ظهور پیدا کرده است. رویکردهای براساس متاهیوریستیک‌های GA و یا NN به وضوح از نظر عملکرد عقب‌تر هستند و آنهایی که بر اساس

متاهیوریستیک‌های SA یا DA یا AS هستند واقعاً نمی‌توانند با TS رقابت کنند. البته با توجه به بهبودهای عملکرد بدست آمده از AS و GA و اینکه این روش‌ها هنوز کاملاً بهره‌برداری نشده‌اند، این قابلیت را دارند که در آینده جواب‌هایی به خوبی TS بدست دهند. موضوعی که وجود دارد این است که مجموعه‌های داده‌ها[1] در حال حاضر برای مقایسه استفاده می‌شوند، ابعاد کوچکی دارند و به یک مدل اجازه نمی‌دهند تا تفاوتش را به وضوح نسبت به بقیه متاهیوریستیک‌ها نشان دهد. بنابراین به مجموعه داده‌های مربوط به مثال‌ها در ابعاد بزرگتر نیاز داریم.

جدول 1-4- نتایج محاسباتی الگوریتم‌های AS

| کد مسئله | Hybrid ant systems | | AS بهبودیافته | | بهترین جواب شناخته شده |
|---|---|---|---|---|---|
| | f* | زمان | f* | زمان | |
| E051-05e | 524/61 | 0/6 | 524/61 | 0/1 | 524/61 |
| E076-10e | 870/8 | 2/4 | 844/31 | 1/3 | 835/26 |
| E101-08e | 879/43 | 11/3 | 832/32 | 3/8 | 826/14 |
| E101-10c | 819/96 | 10/1 | 819/56 | 5/0 | 819/56 |
| E121-07c | 1072/45 | 16/2 | 1065/21 | 9/2 | 1042/11 |
| E151-12c | 1147/41 | 28/5 | 1061/55 | 18/4 | 1028/42 |
| E200-17c | 1473/40 | 82/2 | 1343/49 | 87/6 | 1291/45 |
| E051-06c | 562/93 | 0/2 | 560/24 | 0/1 | 555/43 |
| E076-11c | 948/16 | 3/5 | 916/21 | 1/7 | 909/68 |
| E101-09c | 886/17 | 7/3 | 866/74 | 4/8 | 865/94 |
| E101-11c | 869/86 | 3/1 | 867/74 | 5/8 | 866/37 |
| E121-11c | 1590/52 | 4/3 | 1559/92 | 11/0 | 1541/14 |
| E151-14c | 1202/01 | 26/6 | 1195/99 | 27/5 | 1162/55 |
| E200-18c | 1504/79 | 57/3 | 1451/65 | 81/8 | 1395/85 |

---

[1] - data sets

جدول ۱-۵- نتایج محاسباتی الگوریتم‌های جستجوی ممنوع

| Problem | Osman (BA) | | Taillard[8] | Taburoute[9] | | Rochat and Taillard[10] | Xu and Kelly[4,5] | | Rego and Roucairol[11] | Toth and Vigu[12] | |
|---|---|---|---|---|---|---|---|---|---|---|---|
| | f* | Time¹ | f* | f* | Time² | f* | f* | Time³ | f* | f* | Time⁴ |
| E۰۵۱-۰۵e | ۵۲۴/۶۱ | ۱/۱۲ | ۵۲۴/۶۱ | ۵۲۴/۶۱ | ۶/۰ | ۵۲۴/۶۱ | ۵۲۴/۶۱^{۴/۵} | ۲۹/۲۲^{۴/۵} | ۵۲۴/۶۱ | ۵۲۴/۶۱ | ۰/۸۱ |
| E۰۷۶-۱۰e | ۸۴۴ | ۱/۱۸ | ۸۳۵/۲۶ | ۸۳۵/۷۷ | ۵۳/۸ | ۸۳۵/۳۲ | ۸۳۵/۲۶^{۴/۵} | ۴۸/۸۰^{۴/۵} | ۸۳۵/۳۲ | ۸۳۸/۶۰ | ۲/۲۱ |
| E۱۰۱-۰۸e | ۸۳۵ | ۱۱/۲۵ | ۸۲۶/۱۴ | ۸۲۹/۴۵ | ۱۸/۴ | ۸۲۶/۱۴ | ۸۲۶/۱۴^{۴/۵} | ۷۱/۹۳^{۴/۵} | ۸۲۷/۵۳ | ۸۲۸/۵۶ | ۲/۳۹ |
| E۱۰۱-۱۰c | ۸۱۹/۵۹ | ۶/۷۹ | ۸۱۹/۵۶ | ۸۱۹/۵۶ | ۱۶/۰ | ۸۱۹/۵۶ | ۸۱۹/۵۶^{۴/۵} | ۵۶/۶۱^{۴/۵} | ۸۱۹/۵۶ | ۸۱۹/۵۶ | ۱/۱۰ |
| E۱۲۱-۰۷c | ۱۰۴۲/۱۱ | ۲۳/۳۱ | ۱۰۴۲/۱۱ | ۱۰۷۳/۴۷ | ۲۲/۲ | ۱۰۴۲/۱۱ | ۱۰۴۲/۱۱^{۴/۵} | ۹۱/۲۳^{۴/۵} | ۱۰۴۲/۱۱ | ۱۰۴۲/۸۷ | ۳/۱۸ |
| E۱۵۱-۱۲c | ۱۰۵۲ | ۵۱/۲۵ | ۱۰۲۸/۴۲ | ۱۰۳۶/۱۶ | ۵۸/۸ | ۱۰۳۱/۰۷ | ۱۰۲۹/۵۶^{۴/۵} | ۱۴۹/۹۰^{۴/۵} | ۱۰۴۴/۳۵ | ۱۰۳۳/۲۱ | ۴/۵۱ |
| E۲۰۰-۱۷c | ۱۳۵۴ | ۳۲/۸۸ | ۱۲۹۸/۷۹ | ۱۳۲۲/۶۵ | ۹۰/۹ | ۱۳۱۱/۳۵ | ۱۲۹۸/۵۸^{۴/۵} | ۲۷۲/۵۲^{۴/۵} | ۱۳۳۴/۵۵ | ۱۳۱۸/۲۵ | ۷/۵۰ |
| E۰۵۱-۰۶c | ۵۵۵/۴۴ | ۲/۳۴ | ۵۵۵/۴۳ | ۵۵۵/۴۳ | ۱۳/۵ | ۵۵۵/۴۳ | ۵۵۵/۴۳^{۴/۵} | ۳۰/۶۷^{۴/۵} | ۵۵۵/۴۳ | ۵۵۵/۴۳ | ۰/۸۶ |
| E۰۷۶-۱۱c | ۹۱۳ | ۳/۳۸ | ۹۰۹/۶۸ | ۹۱۳/۲۳ | ۵۴/۶ | ۹۰۹/۶۸ | ۹۶۵/۶۲^{۴/۵} | ۱۰۲/۱۳^{۴/۵} | ۹۰۹/۶۸ | ۹۲۰/۷۲ | ۲/۷۵ |
| E۱۰۱-۰۹c | ۸۶۶/۷۵ | ۲۰/۰۰ | ۸۶۵/۹۴ | ۸۶۵/۹۴ | ۲۵/۶ | ۸۶۵/۹۴ | ۸۸۱/۳۸^{۴/۵} | ۹۸/۱۵^{۴/۵} | ۸۶۶/۷۵ | ۸۶۹/۴۸ | ۲/۹۰ |
| E۱۰۱-۱۱c | ۸۶۶/۳۷ | ۹۲/۹۸ | ۸۶۶/۳۷ | ۸۶۶/۳۷ | ۶۵/۷ | ۸۶۶/۳۷ | ۹۱۵/۲۴^{۴/۵} | ۱۵۲/۹۸^{۴/۵} | ۸۶۶/۳۷ | ۸۶۶/۳۷ | ۱/۴۱ |
| E۱۲۱-۱۱c | ۱۵۴۷ | ۲۲/۳۸ | ۱۵۴۱/۱۴ | ۱۵۷۳/۸۱ | ۵۹/۲ | ۱۵۴۵/۹۳ | ۱۶۱۸/۵۵^{۴/۵} | ۲۰۱/۷۵^{۴/۵} | ۱۵۵۰/۱۷ | ۱۵۴۵/۵۱ | ۹/۳۴ |
| E۱۵۱-۱۴c | ۱۱۸۸ | ۴۰/۷۳ | ۱۱۶۲/۵۵ | ۱۱۷۷/۷۶ | ۷۱/۰ | ۱۱۶۲/۸۹ | ۱۴۳۹/۲۹^{۴/۵} | ۱۶۸/۰۸^{۴/۵} | ۱۱۶۴/۱۲ | ۱۱۷۳/۱۲ | ۵/۶۷ |
| E۲۰۰-۱۸c | ۱۴۲۲ | ۵۵/۱۷ | ۱۳۹۷/۹۴ | ۱۴۱۸/۵۱ | ۹۹/۸ | ۱۴۰۴/۷۵ | | ۳۶۸/۳۷^{۴/۵} | ۱۴۲۰/۸۴ | ۱۴۳۵/۷۴ | ۹/۱۱ |

حال ببینیم در عمل این متاهیوریستیک‌ها در مثال‌هایی در ابعاد خیلی بزرگ که اغلب در کاربردهای عملی با آنها روبرو هستیم، چگونه عمل می‌کند. شاید در مثال‌های خیلی بزرگ، با توجه به پیچیدگی محاسباتی این متاهیوریستیک‌ها، نتوانند در زمان معقولی به جواب برسند. بنابراین در کاربردهایی که سرعت تصمیم‌گیری اهمیت بالایی دارد، این متدها کارا نیستند. اصولاً متاهیوریستیک‌ها نسبت به هیوریستیک‌های کلاسیک زمان برتر می‌باشند، اما درعوض جواب‌های بسیار بهتری را فراهم می‌کنند. نوعاً، متدهای کلاسیک جواب‌های تولیدی‌شان بین ۲ تا ۱۰٪ بدتر از مقدار بهینه (یا بهترین جواب شناخته شده) می‌باشد،

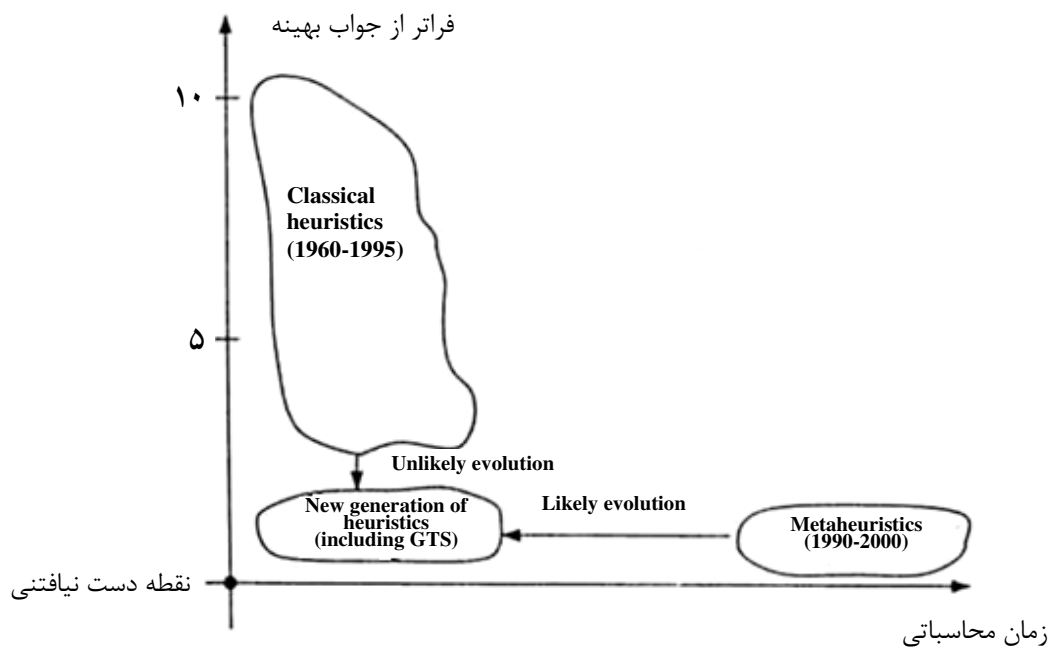

**شکل ۱-۹- نمودار تکامل تدریجی روش‌های ابتکاری برای حل مسائل مسیریابی**

در حالیکه همین رقم برای بهترین متاهیوریستیک‌ها اغلب کمتر از ۰٫۵٪ است (شکل ۱-۹). حال زمان آن است که روش‌های ساده‌تر که توانایی فراهم کردن جواب‌ها با کیفیت خوب در زمان کوتاه را دارند، ایجاد گردند. احتمال انجام این کار با سرعت بخشیدن به بهترین متاهیوریستیک‌های حاضر محتمل‌تر است تا اینکه بخواهیم روی بهبود جواب‌های هیوریستیک‌های کلاسیک کار کنیم. الگوریتم GTS که توسط Toth و Vigo پیشنهاد شده یک گام مهم در این راستا می‌باشد[۶۰]. GTS از تجارب بسیار زیاد در زمینه‌ی متاهیوریستیک‌ها و استخراج بهترین مفاهیم آنها بدست آمده. GTS با شناسایی دقیق ساختار مسئله موفق شده است که از بسیاری از محاسبات غیرضروری که در الگوریتم TS قبلی صورت می‌گرفت، اجتناب کند.

# فصل دوم

# بهینه‌سازی تصادفی

## مقدمه

با در نظر گرفتن یک مسأله‌ی برنامه‌ریزی ریاضی، پارامترهای مختلفی از مسأله ممکن است خاصیت تصادفی داشته باشند و به صورت متغیرهای تصادفی توصیف گردند. به عنوان مثال هزینه‌ی حمل و نقل به طور طبیعی وابسته به قیمت سوخت است که خود متغیری تصادفی است. تقاضا در آینده نیز خاصیتی احتمالی دارد و به شرایط غیرقطعی بازار وابسته است. هر دو این پارامترها در مسأله‌ی مسیریابی می‌تواند به صورت تصادفی در نظر گرفته شود و در واقع مسأله ما را با دنیای واقعی و کاربردهای واقعی نزدیک‌تر خواهد کرد.

در این رساله عدم قطعیت با استفاده از پیشامدهای تصادفی‌ای که توسط w نمایش داده می‌شوند، مشخص می‌گردند. مجموعه همه‌ی این پیشامدها، توسط $\Omega$ نمایش داده می‌شود. همچنین مجموعه‌ای از پیشامدها را با $\mathcal{A}$ نمایش می‌دهیم. بر این اساس فضای تصادفی را می‌توان با سه تایی مرتب $(\Omega, A, P)$ نمایش داد [۳۶] (Chung-۱۹۷۴) که در آن P تابعی است به هر $A \in \mathcal{A}$ مقدار $P(A)$ را تخصیص می‌دهد، به‌طوری که:

$P(A_1 \cup A_2) = P(A_1) + P(A_2)$ اگر $(A_1 \cap A_2 = \phi)$

$0 \leq P(A) \leq 1$

$P(\phi) = 0$
$P(\Omega) = 1$

در ارتباط با بهینه‌سازی تصادفی، فضای تصادفی بیانگر شرایطی است که مدل براساس آن طراحی شده است و متناظر با هر رخداد عضو فضای تصادفی، شرایطی قطعی را مدل کرده است.

## ۲-۱- مرحله‌ها و تصمیم‌ها

مسائل بهینه‌سازی تصادفی گونه‌های مختلفی دارند. مهمترین گونه این نوع از بهینه‌سازی که در این رساله برای مدل‌سازی و بهینه‌سازی مسأله‌ی مسیریابی مورد استفاده قرار گرفته است، بهینه‌سازی تصادفی دو مرحله‌ای است. در واقع روش‌های چند مرحله توصیف کننده‌ی شرایط تصادفی خاصی هستند که در آن پارامترهای تصادفی به مرور زمان قطعیت می‌یابند. مرحله اوّل تصمیم‌گیری بدون هیچ اطلاعی از پارامترهای تصادفی اتخاذ می‌گردد. در مرحله دوّم بخشی از پارامترها قطعیت یافته و در واقع یکی از رخدادهای محتمل مرحله‌ی اول عینیت می‌یابد. در این مرحله و با توجه به داده‌های جدید مسأله مجدداً بهینه‌سازی شده و تصمیم‌گیری صورت می‌گیرد. پس از آن و در مرحله‌ی سوّم مجدداً بخشی از پارامترهای تصادفی قطعیت یافته و یکی از رخدادهای محتمل مرحله‌ی سوّم عینیت می‌یابد. این کار تا جایی که همه‌ی پارامترهای مسأله قطعیت یابند ادامه می‌یابد. در حالت ۲ مرحله‌ای تصمیم‌هایی که پیش از رخداد پارامترهای تصادفی گرفته می‌شود، تصمیم‌های مرحله‌ی اوّل[1] و تصمیم‌های پس از رخداد، تصمیم‌های مرحله‌ی دوّم[2] خوانده می‌شوند.

## ۲-۲- بهینه‌سازی احتمالی دو مرحله‌ای با Recourse قطعی

برنامه‌ریزی خطی احتمالی دو مرحله‌ای کلاسیک که اولین بار توسط دانتزیک و Beale [۱۹۵۵] پیشنهاد شد، به صورت زیر مدل‌سازی شده بود. [۲۹]

$$\min Z = c^T x + E_\xi \left[ \min q(w)^T y(w) \right]$$

S.T.

$$Ax = b$$

$$T(w)x + Wy(w) = h(w)$$  (۲-۱)

$$x \geq 0, y(w) \geq 0$$

در این مدل متغیرهای مرحله‌ی اول و مرحله‌ی دوّم از هم متمایز شده‌اند. متغیرهای مرحله‌ی اوّل توسط بردار x با اندازه‌ی $n_1 \times 1$ نمایش داده می‌شوند، همچنین متناظر با متغیرهای مرحله‌ی اوّل بردارهای c، b و A پارامترهای قطعی مسأله را مشخص می‌کنند. در مرحله‌ی دوّم ممکن است تعدادی از

---
1- First stage decision
2- Second stage decision

۳۱

پیشامدهای تصادفی ($w \in \Omega$) رخ دهند. متناظر با هر رخداد $w$، پارامترهای مرحله‌ی دوّم $q(w)$، $h(w)$ و $T(w)$ مشخص می‌شوند.

در این مدل، c بردار هزینه‌ی با مقدار معلوم (قطعی) در فضای $\mathbb{R}^{n_1}$، b بردار با مقدار معلوم (قطعی) در فضای $\mathbb{R}^{m_1}$، A و W ماتریس‌های معلوم (قطعی) با اندازه‌های (به ترتیب) $m_1 \times n_1$ و $m_2 \times n_2$ هستند و W ماتریس recourse خواهد بود.

w پارامتر متغیر تصادفی گسسته است و برای هر پیشامد آن $T(w)$ از اندازه‌ی $m_2 \times n_1$، $q(w) \in \mathbb{R}^{n_2}$ و $h(w) \in \mathbb{R}^{m_2}$ خواهند بود. همچنین اگر $\xi^T(w) = (q(w)^T, h(w)^T, T_1(w), ..., T_{m_2}(w))$ را به عنوان متغیر تصادفی گسسته فرض کنیم، $E_\xi$ نشان دهنده‌ی، امید ریاضی با توجه به متغیر تصادفی $\xi$ خواهد بود.

همچنین اگر $N = n_2 + m_2 + (m_2 \times n_1)$ باشد، $\Xi \subseteq \mathbb{R}^N$ همدامنه‌ی متغیر تصادفی $\xi$ خواهد بود، به شکلی که $P\{\xi \in \Xi\} = 1$

(در واقع $\Xi$ کوچکترین زیرمجموعه بسته است در $\mathbb{R}^N$ که $P(\Xi) = 1$ باشد.)

همچنان که گفته شد زمانی که رخداد تصادفی w عینیت یافت، داده‌های مرحله‌ی دوّم، یعنی q، h و T مشخص می‌شوند. آنگاه متغیر تصمیم مرحله‌ی دوّم، یعنی $y(w)$ و در واقع همان $y(w, x)$ بایستی مشخص شود. طبیعت بستگی y به w کاملاً با طبیعت وابستگی q، h و T به w متفاوت است. در واقع q، h و T خود متغیرهای تصادفی وابسته به w هستند، امّا y یک متغیر تصمیم است که پس از مشخص شدن مقدار w و از طریق بهینه‌سازی مسأله مرحله‌ی دوّم مشخص می‌شود.

برای شرح مدل فوق بایستی توجه کرد که تابع هدف از ۲ بخش تشکیل شده است. بخش اوّل $c^T x$ است که در واقع بخش قطعی تابع هدف را تشکیل می‌دهد و بخش دوّم، امید ریاضی تابع هدف $q(w)^T y(w)$ که براساس همه‌ی پیشامدهای متغیر تصادفی w تعیین می‌شود. در واقع سختی مسائل بهینه‌سازی تصادفی در همین بخش دوّم است، چرا که برای هر w، مقدار $y(w)$ بایستی با حل یک مسأله برنامه‌ریزی ریاضی محاسبه شود. برای درک این مسأله مدل بالا را به صورت زیر بازنویسی می‌کنیم:

$$\min Z = C^T x + Q(x)$$
  S.T.
  $$Ax = b \qquad (2-2)$$
  $$x \geq 0$$

$$Q(x, \xi(w)) = \min_y \{ q(w)T_y | W_y = h(w) - T(w)x, y \geq \circ \}$$

در جایی که

$$Q(x) = E_\xi Q(x, \xi(w))$$

در اینجا $Q(x,\xi(w))$، تابع مقدار[1] مرحله‌ی دوّم به حساب می‌آید و $Q(x)$ امید ریاضی آن.
مدل بالا در واقع مسأله قطعی معادل[2] (غیراحتمالی) مسأله احتمالی به حساب می‌آید و عملاً ساده‌ترین نوع نمایش مسائل از این دست است. برای مدل‌سازی مسائل عدد صحیح تنها کافی است متغیرهای تصمیم‌گیری را به شکل عدد صحیح معرفی کنیم. برای تشریح الگوریتم‌های حل مسائل احتمالی دو مرحله‌ای نیاز به تعریف چند مفهوم و بیان چند قضیه است که در ادامه می‌آید.

لازم به ذکر است، آنچه در ادامه می‌آید حداقل‌های ممکن برای تشریح الگوریتم‌هاست و برای اطلاع از مبانی تئوریک دقیق مطالعه منابع مشروح چون [۲۹] [۱۹۹۷] Birge، [۵۹] [۲۰۰۴] Kall و [۶۷] [۲۰۰۲] Shapiro توصیه می‌گردد.

## ۲-۳- مجموعه‌های شدنی

با فرض اینکه $\xi$ یک متغیر گسسته‌ی متناهی باشد، $k_1$ را به صورت زیر تعریف می‌کنیم.

$$k_1 = \{ x | Ax = b, x \geq \circ \} \tag{۲-۳}$$

این مجموعه، فضای شدنی برای مرحله اوّل مسأله را مشخص می‌سازد و عملاً محدوده‌ی شدنی متغیرهای تصمیم قطعی با توجه به محدودیت‌های قطعی در آن مشخص می‌گردند. این مقادیر به اینکه کدام یک از پیشامدها رخ دهد بستگی ندارند. به همین ترتیب مجموعه $k_2$ را به صورت زیر تعریف می‌کنیم.

$$k_2 = \{ x | Q(x) < \infty \} \tag{۲-۴}$$

$k_2$ فضای شدنی مرحله‌ی دوّم را برای متغیرهای تصمیم مرحله‌ی اوّل مشخص می‌کند. در واقع این مجموعه مقادیر شدنی را برای متغیرهای تصمیم قطعی با توجه به پیشامدهای ممکن مشخص می‌کند. xهایی که ما به ازای آنها، فضای شدنی برای محدودیت‌های مرحله‌ی دوّم (محدودیت‌های احتمالی) با توجه به همه‌ی پیشامدهای ممکن، تهی نباشد.

---

1- Value Function
2- deterministic equivalent program (DEP)

حال با توجه به تعاریف بالا مسأله‌ی بهینه‌سازی احتمالی دو مرحله‌ای را به صورت زیر نیز می‌توان نوشت:

$$\min Z(x) = C^T x + Q(x)$$
S.T. (۲-۵)
$$x \in k_1 \cap k_2$$

با توجه به اینکه تعریف فضای $k_2$ عملاً کار ساده‌ای نیست. تعریف جایگزین زیر را ارائه می‌دهیم:

$$k_2(\xi) = \{x \mid Q(x,\xi) < +\infty\} \qquad (۲-۶)$$

بدین ترتیب می‌توان مجموعه‌ی $k_2^p$ را تعریف کرد. این مجموعه، فضای شدنی هر پیشامد را مشخص می‌کند. یعنی در حالیکه برای هر $\xi \in \Xi$، $y \geq 0$ وجود داشته باشد می‌توان $k_2^p$ را بصورت زیر تعریف کرد:

$$k_2^p = \{Wy = h - Tx\} = \bigcap_{\xi \in \Xi} k_2(\xi)$$

بدین ترتیب $k_2^p$ معادل احتمالاتی $k_2$ خواهد بود. در واقع x متعلق به مجموعه‌ی $k_2^p$ خواهد بود اگر برای همه‌ی پیشامدهای ممکن (مقدار بردار تصادفی $\xi$) مقدار شدنی برای y موجود باشد.

قضیه ۱- اگر $\Xi$ متناهی باشد، $k_2^p$ یک چندوجهی است و یک مجموعه‌ی بسته و محدب می‌باشد [۲۹] Birge-[۱۹۹۷].

قضیه ۲- $Q(x,\xi)$ با recourse قطعی، یک تابع محدب خطی قطعه قطعه است برای هر x در $k = k_1 \cap k_2$.

اثبات:

تنها لازم است اثبات کنیم که $f(b) = \min\{q^T y \mid wy = b\}$ یک تابع محدب در b است. با فرض ۲ بردار متفاوت $b_1$ و $b_2$ و ترکیب محدب $b_\lambda = \lambda b_1 + (1-\lambda) b_2$ برای $\lambda \in (0,1)$ فرض کنید $y_1^*$ و $y_2^*$ جواب بهینه مسأله به ازای $b = b_1$ و $b = b_2$ باشند. همچنین $\lambda y_1^* + (1-\lambda) y_2^*$ یک جواب شدنی برای $\min\{q^T y \mid wy = b_\lambda\}$ است. حال فرض کنید $y_\lambda^*$ جواب بهینه‌ی مسأله آخر باشد.

$$f(b_\lambda) = q^T y_\lambda^* \leq q^T(\lambda y_1^* + (1-\lambda) y_2^*) = \lambda q^T y_1^* + (1-\lambda) q^T y_2^* = \lambda f(b_1) + (1-\lambda) f(b_2)$$

با توجه به تعریف بهینه‌سازی تصادفی در حالت دو مرحله‌ای و قضایایی که طرح شد، تکنیک‌های حل مسائل بهینه‌سازی از این دست که در این رساله مورد استفاده قرار گرفته است شرح داده می‌شود. با توجه به پیچیدگی زیاد این مسأله الگوریتم در حالت خطی معرفی شده و سپس به حالت عدد صحیح

توسعه می‌یابد. با توجه به استفاده‌ی گسترده این رساله از مفهوم دوگان در فضاهای گسسته و عدد صحیح و عدم وجود منابع کافی، مفهوم دوگان در فضای گسسته به پیوست شرح داده شده است. همچنین از آنجا که مسائل بهینه‌سازی تصادفی از جمله مسائلی هستند که حتی مثال‌های کوچک آنها، جزو مسائل برنامه‌ریزی ریاضی با ابعاد بزرگ[1] به حساب می‌آیند، لذا تکنیک‌های حل مسأله‌این گونه‌ی بهینه‌سازی غالباً بر گرفته از تکنیک‌های بهینه‌سازی در ابعاد بزرگ[2] است. لذا اصل تجزیه که مبنای اصلی روش‌های بهینه‌سازی در ابعاد بزرگ به حساب می‌آید، از فصل 7 کتاب بازارا، ترجمه و خلاصه و به پیوست آمده است. تکنیک L-Shape که در واقع استفاده از همین مفهوم در بهینه‌سازی تصادفی است در ادامه شرح داده شده است. این تکنیک اصلی‌ترین ابزار برای توسعه الگوریتم حل دقیق مسأله مسیریابی توسعه یافته احتمالی در این رساله است.

## 2-4- روش L-Shape

با در نظر گرفتن مدل بهینه‌سازی احتمالی در مرحله‌ای که در بخش قبل ارائه شد و فرض اینکه $\delta$ هم‌دامنه‌ی متناهی داشته باشد، برای مشخص کردن رخدادها از ایندکس $k=1,...,K$ که هر کدام نمایش دهنده‌ی یک رخداد خواهند بود، استفاده می‌کنیم.

با این فرض می‌توان معادل قطعی مسأله را به صورت گسترده[3] نمایش داد. این کار با تخصیص $y_k$ به هر رخداد به عنوان متغیر تصمیم مرحله‌ی دوّم با پارامترهای $q_k$، $h_k$ و $T_k$ امکان پذیر است. در این صورت می‌توان مسأله Large Scale معادل را در حالت خطی به صورت زیر بازنویسی کرد:

$$\min C^T x + \sum_{k=1}^{k} p_k q_k^T y_k$$

S.T.

$\quad Ax = b$ \hfill (2-7)

$\quad T_k x + W y_k = h_k \qquad k=1,...,k$ \hfill (2-8)

$\quad x \geq 0, y_k \geq 0 \qquad k=1,...,k$

با دقت در فرم مدل بدست آمده و مقایسه آن با ساختار بلوکی شرح داده شده در پیوست، می‌توان از روش تجزیه برای حل مسأله استفاده کرد.

---

1- Large Scale Programming
2- Large Scale Optimization
3- extensive form

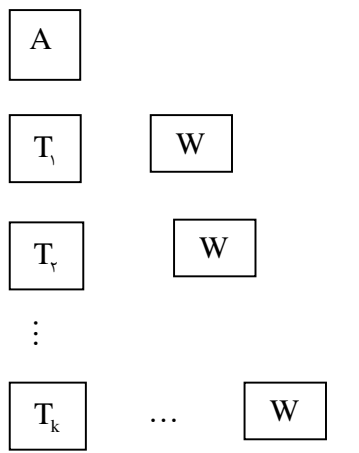

**شکل ۲-۱- ساختار بلوکی فرم گسترده مسأله بهینه‌سازی احتمالی دو مرحله‌ای**

این فرم الهام بخش نام L-Shape برای تکنیکی است که از این ساختار برای حل مسأله استفاده می‌کند. در واقع دوگان مدل فوق همان ساختار دیاگرام بلوک[1] اصل تجزیه خواهد بود. لذا استفاده از تکنیک تجزیه Wolfe-Dantzig [1960] [۱۷] برای دوگان مسأله یا تجزیه Benders [1962] [۱] برای مسأله‌ی اولیه طبیعی است و منجر به طراحی تکنیک L-Shape خواهد شد. در واقع تجزیه Benders[2] توسط Wets [1969] و Van Slyke برای حل مسائل بهینه‌سازی احتمالی توسعه یافته و به روش L-Shape معروف شده است [۲۹].

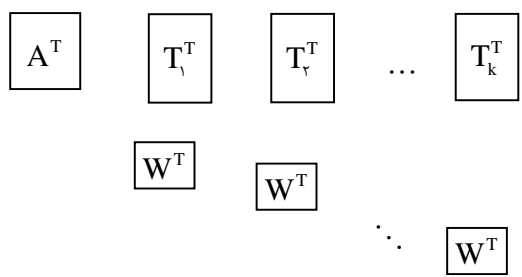

**شکل ۲-۲- ساختار بلوکی دوگان مسأله‌ی بهینه‌سازی احتمالی دومرحله‌ای**

---

1- Block diagram structure
2- Benders decomposition

برای تشریح روش L-Shape برای حل مسائل بهینه‌سازی احتمالی دومرحله‌ای در حالت خطی، ابتدا گام‌های روش ارائه شده و سپس به تشریح هر گام می‌پردازیم.

**گام ٠-** r و s و v را برابر با صفر قرار دهید.

**گام ١-** v را برابر با v+۱ قرار دهید و مسأله‌ی خطی زیر را حل کنید:

$$\min C^T x + \theta \qquad (٢-٩)$$

S.T.

$$Ax = b$$

$$D_\ell x \geq d_\ell \qquad \ell = ١, ..., v \qquad (٢-١٠)$$

$$E_\ell x + \theta \geq e_\ell \qquad \ell = ١, ..., s \qquad (٢-١١)$$

$$x \geq ٠ \qquad \theta \in \mathbb{R}$$

حل بهینه را برابر با $(x^v, \theta^v)$ قرار دهید. اگر محدودیت نوع (٢-١١) وجود نداشته باشد، $\theta^v$ برابر با $-\infty$ خواهد بود.

**گام ٢-** برای $k = ١, ..., K$ مسأله‌ی خطی زیر را حل کنید: ($e^T = (١, ..., ١)$)

$$\min w' = e^T v^+ + e^T v^- \qquad (٢-١٢)$$

S.T.

$$Wy + Iv^+ - Iv^- = h_k - T_k x^v \qquad (٢-١٣)$$

$$y \geq ٠, v^+ \geq ٠, v^- \geq ٠$$

این کار را تا جایی ادامه دهید که به ازای یک k، جواب بهینه‌ی $w'$ مقداری مثبت شود. در صورتی که k ای با مقدار بهینه‌ی مثبت $w'$ وجود داشت، ضرایب دوگان مسأله‌ی معادل آن را برابر با $\sigma^v$ قرار داده و ضرایب زیر را تعریف کنید:

$$D_{r+١} = (\sigma^v)^T T_k \qquad (٢-١۴)$$

$$d_{r+١} = (\sigma^v)^T h_k \qquad (٢-١۵)$$

با استفاده از این ضرایب محدودیت‌های (٢-١٠) ساخته خواهند شد. r را برابر با r+1 قرار دهید و با اضافه کردن محدودیت‌های تولید شده به باقی محدودیت‌های (٢-١٠)، به مرحله‌ی ١ بازگردید. این کار را تا جایی ادامه دهید که به ازای هیچ یک از موارد k، $w'$ بهینه‌ی مثبت وجود نداشته باشد.

**گام ۳**- برای $k = 1, ..., K$ مسأله خطی زیر را حل کنید.

$$\min w = q_k^T y$$
$$\text{S.T.}$$
$$Wy = h_k - T_k x^v \qquad (2-16)$$
$$y \geq 0$$

فرض کنید $\pi_k^v$ ضرایب دوگان مرتبط با جواب بهینه برای $k$ اُمین مسأله باشد. با استفاده از آنها عبارت‌های زیر را تعریف کنید:

$$E_{s+1} = \sum_{k=1}^{k} p_k \cdot (\pi_k^v)^T \cdot T_k \qquad (2-17)$$

$$e_{s+1} = \sum_{k=1}^{k} p_k \cdot (\pi_k^v)^T h_k \qquad (2-18)$$

$w^v$ را به صورت مقابل محاسبه کنید: $w^v = e_{s+1} - E_{s+1} x^v$

اگر $\theta^v \geq w^v$ توقف کنید. $x^v$ جواب بهینه مسأله است. در غیر این صورت $s = s+1$ قرار داده و به گام ۱ باز گردید.

در این ۳ گام ما با دو مفهوم برش‌های فضای شدنی[1] و برش‌های بهینه‌یاب[2] روبرو هستیم. در واقع محدودیت‌های (۲-۱۰) برش‌های فضای شدنی و محدودیت‌های (۲-۱۱) برش‌های بهینه‌یاب هستند. در تشریح کلی روش L-Shape بایستی به مدل صورت گسترده بازگردیم. در این مدل تابع هدف ۲ بخش دارد که بخش اوّل یعنی $c^t x$ تابع هدف مرحله‌ی اوّل و باقی تابع هدف، تابع مقدار مرحله‌ی دوّم است. محدودیت‌های (۲-۷) محدودیت‌های مرحله اوّل و محدودیت‌های (۲-۸) محدودیت‌های مرحله‌ی دوّم هستند. متغیرهای تصمیم‌گیری x، متغیرهای مرحله‌ی اوّل و متغیرهای تصمیم‌گیری y، متغیرهای مرحله‌ی دوّم می‌باشند.

با حل مدل ارائه شده در گام ۱، در اولین مرحله جواب بهینه‌ی حاصل، نقطه‌ای از فضای شدنی متغیر تصمیم مرحله‌ی اوّل ($k_1$) را به ما می‌دهد. این نقطه در گام دوّم با استفاده از مدل گام ۲، برای عضویت در $k_2$ مورد آزمایش قرار می‌گیرد. اگر برای هیچ کدام از پیشامدها $w'$ بهینه مثبت نشود، نقطه‌ی مورد نظر $k_2$ عضو خواهد بود. در غیر این صورت یک برش فضای شدنی با استفاده از ضرایب دوگان حاصل از مدل پیشامدی که $w'$ بهینه‌ی آن مثبت بوده است ساخته می‌شود و به مدل گام ۱ اضافه می‌گردد. این کار تا زمانی صورت می‌پذیرد که متغیر تصمیم مرحله‌ی اوّل متعلق به $k_1 \cap k_2$ باشد.

---

1- Feasible cut
2- Optimal cut

۳۸

در گام سوّم، برای نقطه‌ی مورد نظر با استفاده از ضرایب دوگان حاصل از بهینه‌سازی محدودیت‌های مرحله‌ی دوّم با توجه به تابع هدف مرحله‌ی دوّم برش بهینه‌یاب ایجاد می‌گردد. اگر $\theta$ بهینه‌ی حاصل از حل مدل گام ۱ در این برش صدق کند، نقطه‌ی مزبور حل بهینه‌ی مدل است. در غیراینصورت این برش به مدل اضافه شده و بهینه‌سازی از گام اوّل مجدداً شروع می‌شود. برای بیان ریاضی آنچه شرح داده شد مدل بهینه‌سازی احتمالی دو مرحله‌ای را در حالت خطی به صورت زیر بازنویسی می‌کنیم:

$$\min C^T x + Q(x) \qquad (2\text{-}19)$$
$$\text{S.T.}$$
$$x \in k_1 \cap k_2$$

این مدل معادل مدل زیر خواهد بود:

$$\min C^T x + \theta$$
$$Q(x) \leq \theta$$
$$\text{S.T.}$$
$$x \in k_1 \cap k_2$$

در حالی که

$$Q(x, \xi(w)) = \min_y \{q(w)^T y \mid Wy = h(w) - T(w)x, y \geq 0\}$$

و

$$Q(x) = E_\omega Q(x, \xi(w))$$

در گام سوّم، مدل (۲-۱۶) به طور مکرر برای $k = 1, \ldots, K$ حل می‌شود و برای هر جواب بهینه (برای هر $k$) ضرایب دوگان $\pi_k^\nu$ را خواهیم داشت. برای هر $k$ داریم:

$$Q(x^\nu, \xi_k) = (\pi_k^\nu)^T (h_k - T_k x^\nu)$$

با توجه به محدب بودن $Q(x, \xi_k)$ [رجوع شود به قضیه ۲] و نامساوی زیر گرادیان[1] خواهیم داشت:

$$Q(x, \xi_k) \geq (\pi_k^\nu)^T h_k - (\pi_k^\nu)^T T_k x$$

اثبات:

نامساوی زیر گرادیان: (برای توابع محدب)

$$f(x) \geq f(x^\nu) + \nabla^T f(x^\nu)(x - x^\nu)$$
$$f(x^\nu) = Q(x^\nu, \xi_k) = (\pi_k^\nu)^T (h_k - T_k x^\nu)$$
$$\nabla^T f(x^\nu) = (-(\pi_k^\nu)^T . T_k)$$

---

1- Subgradient Inequality

۳۹

$$f(x) = Q(x, \xi_k)$$

$$\Rightarrow Q(x, \xi_k) \geq (\pi_k^v)^T (h_k - T_k x^v) + (-(\pi_k^v)^T . T_k)(x - x^v)$$

$$= (\pi_k^v)^T h_k - (\pi_k^v)^T T_k x^v - (\pi_k^v)^T T_k x + (\pi_k^v)^T T_k x^v$$

$$= (\pi_k^v)^T h_k - (\pi_k^v)^T T_k x$$

با محاسبه‌ی امید ریاضی ۲ طرف این نامساوی خواهیم داشت:

$$Q(x^v) = E(\pi^v)^T (h - Tx^v) = \sum_{k=1}^{k} p_k (\pi_k^v)^T (h_k - T_k x^v)$$

و در نتیجه

$$Q(x) \geq E(\pi^v)^T (h - Tx) = \sum_{k=1}^{k} p_k (\pi_k^v)^T h_k - \left( \sum_{k=1}^{k} p_k (\pi_k^v)^T T_k \right) x$$

با توجه به اینکه $\theta \geq Q(x)$، $(x, \theta)$ تنها برای (۲-۱۶) زمانی‌شدنی خواهند بود که $\theta \geq E(\pi^v)^T (h - Tx)$ باشد. این نامساوی در صورتی که از تعاریف (۲-۱۷) و (۲-۱۸) استفاده کنیم، معادل نامساوی‌های (۲-۱۰) و (۲-۱۱) است.

همچنین $(x^v, \theta^v)$ جواب بهینه‌ی (۲-۱۹) خواهد بود اگر $Q(x^v) = \theta^v$. دلیل این مسأله این است که $\theta$ در مدل مزبور یک متغیر نامحدود بجز در مورد $\theta \geq Q(x)$ است. بدین ترتیب شرط بهینگی در گام سوّم توجیه می‌شود. همچنین این مسأله بدین معنی است که اگر $Q^v < Q(x^v)$، برش‌های بهینه‌یاب نتوانسته‌اند مسأله را به جواب بهینه برسانند و لذا بایستی از ضرایب دوگان $\pi_k^v$ جدید برای ایجاد برش‌های بهینه‌یاب جدید استفاده کرد. دلیل همگرایی این الگوریتم متناهی بودن تعداد ترکیبات مختلف ضرایب $\pi_k$ برای k است. چرا که برای هر رخداد تعداد متناهی از جواب‌های پایه وجود خواهد داشت.

در گام دوم، زیر مسأله‌ی (۲-۱۲) این که آیا $h_k - T_k x^v$ به pos w به ازای $k = 1, ..., K$ تعلق دارد را آزمون می‌کند. اگر نه این بدین معنی است به ازای یکی از kها، $h_k - T_k x^v \notin \text{pos w}$. در نتیجه حتماً ابرصفحه‌ای بین این pos w و $h_k - T_k x^v$ وجود دارد. با در نظر گرفتن $\sigma^T$ به عنوان ضرایب این ابرصفحه، به ازای $t \in \text{pos w}$، نامساوی $\sigma^T t \leq 0$ همواره درست خواهد بود. همچنین است نامساوی $\sigma^T (h_k - T_k x^v) > 0$. (در مرحله‌ی دوّم این ابرصفحه با استفاده از ضرایب دوگان $(\sigma^v)$ زیر مسأله‌ی (۲-۱۲) در نقطه‌ی بهینه بدست می‌آید.)

دلیل درست بودن این نامساوی این است که براساس تئوری دوگان، در صورتی که $w'$ بطور اکید مثبت باشد، به طور معادل $(\sigma^v)^T (h_k - T_k x^v)$ نیز بطور اکید مثبت خواهد بود. همچنین $\sigma^T w \leq 0$ به

۴۰

این دلیل درست است که $\sigma^v$ ضرایب دوگان نقطه‌ی بهینه هستند و در نقطه‌ی بهینه، قیمت‌های سایه‌ای مرتبط با y بایستی غیرمنفی باشند.

براساس آنچه گفته شد، شرایط لازم برای اینکه x به $k_v$ تعلق داشته باشد این است که $(\sigma^v)^T(h_k - T_k x) \leq 0$. همچنین تعداد محدودیت‌های متناهی از نوع (۲-۱۰) (چرا که تعداد پایه‌های مسأله‌ی (۲-۱۲) متناهی است) باعث خواهد شد تعداد مراحل طی شده برای رسیدن به نقطه‌ای شدنی نیز متناهی باشد.

لذا با توجه به اینکه تعداد مراحل در هر گام ۲ و ۳ متناهی هستند و مجموع مراحل طی شده نیز چندجمله‌ای است که از حاصلضرب این مراحل بدست می‌آید. این تکنیک ما را در مراحل متناهی به جواب بهینه می‌رساند.

برای توسعه روش L-Shape برای فضاهای گسسته و عدد صحیح، بایستی از تئوری دوگان در فضای مونوئیدهای محدب استفاده کرد. از این رو تئوری دوگان (Duality (X,D)) برای فضاهای مجرد در پیوست معرفی شده است. با استفاده از این تئوری تمام مراحل روش L-Shape قابل توسعه برای حالت عدد صحیح (و به همین ترتیب حالت صفر و یک) که در توسعه الگوریتم دقیق برای مسأله مسیریابی توسعه یافته احتمال مورد استفاده قرار گرفته است، می‌باشد.

# فصل سوم

# مسأله مسیریابی توسعه یافته احتمالی

**مقدمه**

توانایی هر سیستم در طراحی فرایندهای خود بستگی مستقیم به طراحی شبکه‌هایی دارد که غالباً برای جابجایی یا حمل و نقل مورد استفاده قرار می‌گیرند. عمومی‌ترین مثال برای توصیف چنین شبکه‌هایی مسأله‌ی فروشنده دوره گرد است که ساده‌ترین نمونه از چنین شبکه‌هایی است. هر چند این مسأله در حوزه‌ی گسترده‌ای از کاربردهای عملی از تولید و توزیع تا طراحی شبکه‌های الکترونیکی مورد استفاده قرار می‌گیرد، لیکن حتی کاربرد محدود آن در مسائل حمل و نقل و توزیع اهمیت آن را توجیه می‌کند. از این رو سهم ۱۰ تا ۲۰ درصدی از قیمت تمام شده کالا [۶۰] [۲۰۰۲ .Vigo و Toth] برای هزینه‌های حمل و نقل در تمام سطوح تولید و توزیع این مسأله را تبدیل به مسأله‌ای مهم برای صنایع ساخته است. لیکن با در نظر گرفتن مسائل واقعی و مدل‌های ریاضی در نظر نگرفتن طبیعت تصادفی بسیاری از پارامترهای مسائل در مدل‌سازی و حل است. علت البته در این واقعیت نهفته است که با در نظر گرفتن جنبه احتمالی مسائل مورد بحث، به شدت به ابعاد مسأله افزوده شده و عملاً امکان حل مسائل واقعی در ابعاد واقعی ممکن نخواهد بود.

همین مسأله در رابطه با مسائل کلی‌تر مسیریابی، یعنی مدل‌های مسیریابی، مکان‌یابی[1] نیز صادق است. این مسأله نیز بطور جدی در مقایسه با مسائل عادی مسیریابی، پیچیده‌تر بود، و غالباً در حل آن یا از الگوریتم‌های ابتکاری و فراابتکاری استفاده می‌شود و یا دو مسأله (مسیریابی و مکان‌یابی) بطور جداگانه حل شده و با هم تلفیق می‌گردند.

در این بخش سعی شده است تا در مواجهه با ترکیب این دو، یعنی در نظر گرفتن جنبه احتمالی تقاضا در مسأله مسیریابی و توسعه‌ی مسأله مسیریابی به مدل بهینه‌سازی همزمان مسأله مسیریابی- مکان‌یابی یا همان مسیریابی توسعه یافته، مدلی واقعی‌تر از مسائل کاربردی ارائه کرده و در ضمن با استفاده از روش‌های بهینه‌سازی در ابعاد بزرگ تکنیک مؤثری برای حل بهینه در ابعاد واقعی ارائه کنیم. در انتها الگوریتم ابتکاری نیز مبتنی بر روش ارزش قرار دادن برای حل این مسأله ارائه شده است.

## ۳-۱- ادبیات موضوع

مسأله مسیریابی توسعه یافته احتمالی تاکنون در ادبیات مسأله مسیریابی مورد مطالعه قرار نگرفته است. از این رو مسائل مرتبط با این مسأله مورد بررسی قرار می‌گیرد. مسیریابی توسعه یافته GVRP که در آن پارامترها حالت قطعی دارند و در واقع نوعی مسأله مکان‌یابی- مسیریابی است [۶۰] [۱۹۸۸، Laporte] دو تصمیم همزمان را با هم ترکیب می‌کند. انتخاب گره‌ها و توالی گره‌ها. مسائل کاربردی چندی از این نوع حل شده‌اند. از جمله مسأله پستچی توسط [۶۰] [۱۹۸۹] Laporte et al که در آن اگر وسائل نقلیه بیش از یکی باشند نامتقارن خواهد بود.

مسأله‌ی GTSP که حالت تک وسیله‌ی نقلیه‌ای GVRP است، بیشتر مورد مطالعه قرار گرفته است [۴۳] [۱۹۸۳] Nobert و Laporte و [۴۲] [۱۹۸۷] Laporte et al اولین کسانی بودند که برای این مسأله (GTSP) مدل برنامه‌ریزی عدد صحیح ارائه کردند. آنها روش شاخه و کران را برای حل مسائل با ابعاد کوچک و متوسط مورد استفاده قرار دارند و مسأله را در حال متقارن و نامتقارن بطور بهینه حل کردند. [۵۸] [۱۹۸۸] Noon یک روش ساده‌سازی لاگرانژ[2] و چند روش ابتکاری برای حل مسأله GTSP نامتقارن پیشنهاد کرد. الگوریتم شاخه و برش با نتایج محاسباتی گسترده در ۲ مقاله‌ی متوالی [۱۵ و ۱۶] [۱۹۹۷ و ۱۹۹۵] Fischetti et al ارائه شد و مسائلی را با شبکه متقارن و تا حداکثر ۴۴۲ گره بطور بهینه حل کرد.

---
1- LRP, location-routing problem
2- lagrangian relaxation

مسأله‌ی دیگری که ارتباط نزدیک با مسأله مسیریابی توسعه یافته احتمالی دارد، مسأله‌ی مسیریابی احتمالی (PVRP) است. این مسأله اولین بار توسط [۶۰] [۱۹۶۹] Tillman مورد توجه قرار گرفت. Tillman مسأله‌ی مسیریابی با چند دپو را با تقاضاهای احتمالی (که از توزیع پواسن پی روی می‌کردند) مطرح کرد. راه‌حل در واقع توسعه روش صرفه‌جویی بود که توسط [۶۰] [۱۹۶۴] Wright و Clarke مطرح شد. اولین الگوریتم مؤثر برای حل دقیق مسأله مسیریابی احتمالی توسط [۴۱] [۱۹۹۳] Laporte و Louveaux مبتنی بر روش L–shape پیشنهاد شد. روش آنها مبتنی بر اضافه کردن برش‌های شدنی بر مدل جریان ساده‌سازی شده CVRPSD تا زمانی که یک جواب شدنی عدد صحیح یافت شود، بود. اگر جواب عدد صحیح شدنی یافت نشد، شاخه زنی[1] صورت می‌گرفت. در نهایت به ازای هر نقطه‌ی شدنی، برش بهینه‌یاب اعمال شده و تا رسیدن به جواب بهینه این کار ادامه می‌یافت. همین روش به صورت‌های اصلاح شده برای حل همین مسأله در سال‌های ۱۹۹۵، ۱۹۹۸ و ۲۰۰۲ به کار رفت و در بهترین حالت مسأله‌هایی با ابعاد ۱۰۰ گره را حل کرد.

در سال ۱۹۹۹، Holt و Hjorring روش L-shape را برای حل مسأله CVRPSD (با یک وسیله‌ی نقلیه) با کار بردن مجموعه‌ی جدیدی از برش‌های بهینه‌یاب پیشنهاد داده و این روش را بهبود دادند. همچنین آنها کران‌های پایین بهتری برای حل مسأله پیشنهاد کردند. [۶۷] [۲۰۰۲] Shapiro, Shabir Ahmad برای حل مسائل مسیریابی احتمالی با استفاده از روش L-shape از تکنیک SAA[2] برای حل مؤثرتر این مسأله استفاده کردند. همچنین مسأله PTSP که حالت تک وسیله‌ی نقلیه PVRP است، بیشتر مورد مطالعه قرار گرفته است. این مسأله اولین بار توسط [۳۱] [۱۹۸۸] Jaillet معرفی شد. وی تابع ارزش مرحله دوم را برای هر تور بطور مشخص[3] معرفی کرده و بدین ترتیب امید ریاضی تور توالی غیر قطعی را محاسبه نمود. نتایج محاسباتی برای حل این مسأله (با استفاده از الگوریتم‌های ابتکاری مسأله TSP که توسط Jaillet توسعه یافته بودند) بوسیله‌ی [۶۴] [۱۹۸۷] Rossi et al و همچنین توسط [۵] [۱۹۹۳] Bertsimas و Howell ارائه شدند.

[۴۱] [۱۹۹۳] Laporte et al این مسأله را به صورت مسأله‌ی عدد صحیح خطی احتمالی مدل‌سازی کرده و با استفاده از الگوریتم عدد صحیح L-shape بطور دقیق حل کرد.

---

1- branching
2- Sample Average Approximation
3- explicit

نزدیک‌ترین مطالعه صورت گرفته مربوط به مسأله‌ی مسیریابی توسعه یافته احتمالی، توسط [۶۸] [۲۰۰۷] Tang و Hooks برای مدل‌سازی و حل مسأله PGTSP صورت گرفت. این مطالعه این مسأله را به صورت عدد صحیح احتمالی در مرحله‌ای مدل‌سازی کرده و با استفاده از الگوریتم L-shapc عدد صحیح الگوریتم حل دقیق برای این مسأله ارائه کرده است. همچنین الگوریتم‌های - ابتکاری مبتنی بر مفهوم ارزش قرار دادن را برای این مسأله‌ی توسعه داده و نتایج محاسباتی گسترده‌ای را ارائه کرده است.

## ۳-۲- تعریف مسأله

### ۳-۲-۱- تشریح مسأله‌ی مسیریابی توسعه یافته احتمالی

براساس آنچه در تعریف مسأله‌ی کلاسیک مسیریابی تشریح شد، این مسأله و انواع گونه‌های آن به دنبال تعیین توالی بهینه‌ای از مسیرها برای تشکیل یک مسیر کامل است. به طور معمول در ساده‌ترین حالت، مسأله فروشنده دوره‌گرد این مسیر بهینه را در میان گره‌های مشخصی با فواصل معین می‌جوید، بطوریکه از گره مشخصی شروع شده و به همان گره باز گردد.

همچنین با در نظر گرفتن نامساوی مثلث در این مسیر هر گره تنها یک بار در مسیر قرار خواهد گرفت. مسأله مسیریابی وسیله نقلیه به جای یک مسیر به دنبال مجموعه‌ای از مسیرهاست که در مجموع بایستی همه‌ی گره‌ها را پوشش دهد و همه مسیرها از گره مشخصی شروع و به همان گره بازگردند. لیکن در مدل کردن شرایط واقعی گاهی لزومی به ملاقات همه‌ی گره‌ها نیست، بلکه تنها عبور از یکی از گره‌ها، ما را از ملاقات مجموعه‌ی مشخصی از گره‌ها بی‌نیاز می‌کند. به عبارت دیگر گره‌ها به مجموعه‌ای از خوشه‌ها[1] تقسیم‌بندی می‌شوند و در عمل، هدف، عبور از هر خوشه در هر مسیر (یا تور) ممکن خواهد بود. در واقع اگر مسیری حداقل از یکی از گره‌های هر خوشه عبور نماید کافی است. این حالت توسعه‌ی مدل کلاسیک مسیریابی به حالت خوشه‌ای است. این توسعه هم برای حالتی که تعیین یک مسیر (یک وسیله‌ی نقلیه) مدنظر است و چه برای حالتی که تعیین چند مسیر (چند وسیله‌ی نقلیه) مدنظر است قابل اعمال است. از این رو، مسأله مسیریابی کلاسیک حالت خاصی از مسأله مسیریابی توسعه یافته است. این حالت خاص در صورتی رخ می‌دهد که هر خوشه تنها شامل یک گره باشد.

---

1- Clusters

٤٦

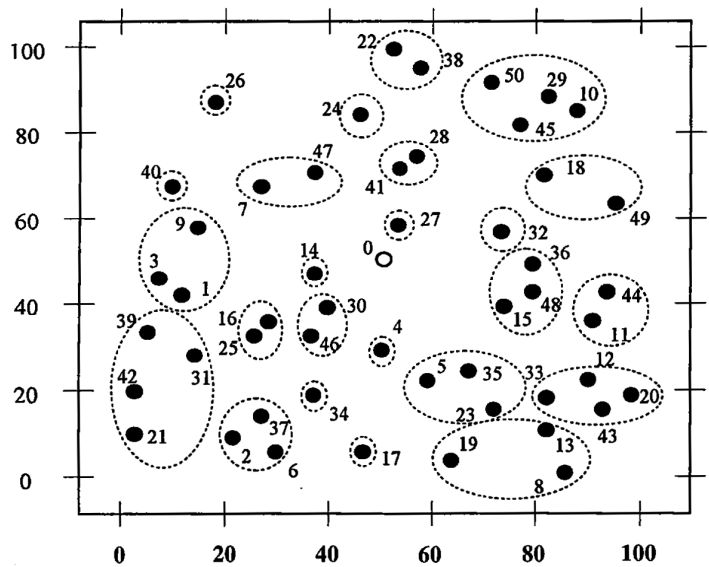

شکل ۳-۱- مجموعه گره‌ها به صورت خوشه‌بندی شده

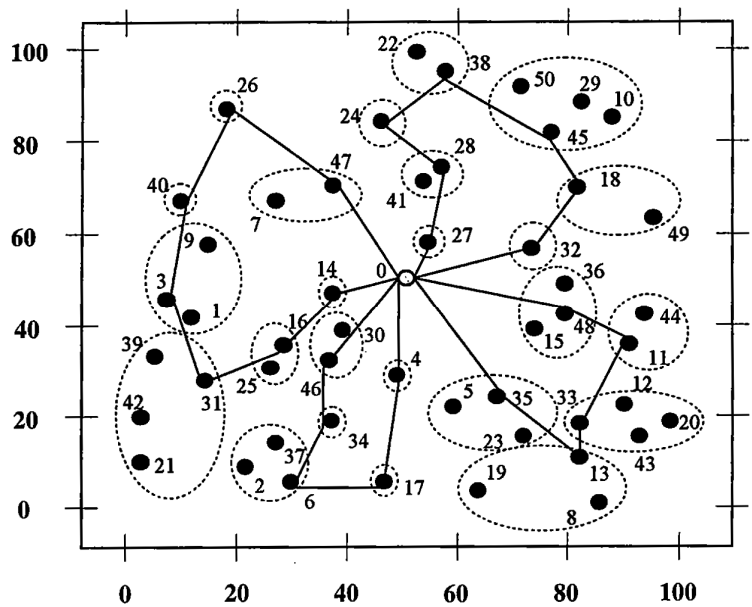

شکل ۳-۲- یک جواب شدنی برای مسأله مسیریابی توسعه‌یافته

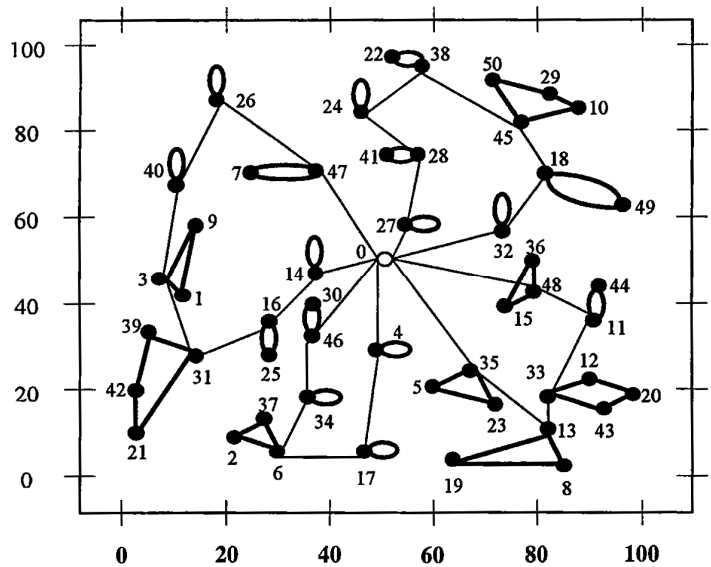

شکل ۳-۳- مسیرهای شدنی درون خوشه‌ای همراه با مسیرهای خارجی

مسائل کاربردی، غالباً برعکس مدل‌های کلاسیک با پارامترهای غیرقطعی مواجه هستند و در مواردی در نظر نگرفتن این مسأله، مدل را از شرایط واقعی به شدت دور می‌کند. در مسأله‌ی مسیریابی نیز موارد زیادی از پارامترها حالت احتمالی دارند. از عمومی‌ترین آنها، هزینه‌ی حمل و نقل بین هر دو گره است که ماهیتاً وابسته به هزینه‌ی سوخت بوده و متغیر است. همچنین تقاضاهای هر گره برای ملاقات شدن در هر تور (مسیر) غالباً متغیری تصادفی است. در این مسأله تقاضای هر گره را بصورت تصادفی مدل کرده و بدین ترتیب فرض می‌کنیم هر گره احتمال حضور $p_i$ و عدم حضور $(1-p_i)$ دارد. مسأله‌ی حاصل، توسعه مسأله مسیریابی به حالت احتمالی است و در نتیجه حالت کلاسیک مسیریابی مثالی خاص از آن خواهد بود که در آن احتمال حضور هر گره ۱ و عدم حضور ۰ باشد.

مسأله‌ای که در این رساله مدل سازی و حل شده است ترکیبی از هر دو توسعه در مسأله مسیریابی وسیله‌ی نقلیه است که در آن گره‌ها به مجموعه‌ای از خوشه‌ها دسته‌بندی شده‌اند و حضور هر خوشه به صورت احتمالی مدل می‌شود.

به این صورت، اگر حضور همه‌ی خوشه‌ها با احتمال ۱ در نظر گرفته شود با مسأله مسیریابی وسیله نقلیه توسعه یافته مواجه خواهیم بود و اگر در هر خوشه تنها یک گره در نظر گرفته شود، با مسأله‌ی مسیریابی وسیله نقلیه احتمالی. این مسأله به دنبال یافتن مسیرهای بهینه با توجه به تعداد وسائل نقلیه در دسترس برای ملاقات هر خوشه حداقل برای یک بار خواهد بود. چنان که گفته شد، ملاقات هر خوشه، با ملاقات حداقل یک گره از آن صورت می‌پذیرد.

## 3-2-2- مدلسازی مسأله مسیریابی توسعه یافته احتمالی

مسأله مسیریابی توسعه یافته احتمالی (GPVRP) بر روی یک گراف کامل به صورت $G=(V,A)$ قابل تعریف است، بطوریکه $V=\{0,1,...,n\}$ مجموعه گره‌ها و $A=\{(i,j)|i,j\in V\}$ مجموعه‌ی مسیرها هستند. همچنین هر مجموعه گره $C_k$ و $1\leq k\leq m$ را یک خوشه[1] می‌نامیم و هر گره عضو $V$ بطور منحصر بفردی عضو تنها یکی از این خوشه‌ها خواهد بود. به ازای هر مسیر $(i,j)\in A$، $d_{ij}$ طول مسیر فرض خواهد شد. GPVRP به دنبال مجموعه‌ای از تورهاست که امید ریاضی مجموع طول آنها حداقل بوده و هر خوشه حداقل یک بار توسط یک تور ملاقات شده باشد. بدیهی است برای ملاقات یک خوشه تنها عبور از یکی از اعضای آن خوشه کافی است. برای مدل سازی مسأله، علاوه بر آنچه گفته شد، تعریف متغیرها و پارامترهای زیر لازم است:

$x_{ij}$: متغیر تصمیم (صفر و یک) حضور یا عدم حضور مسیر $(i,j)\in A$ در یکی از تورهاست. اگر برابر با ۱ باشد مسیر $(i,j)$ در تور حضور خواهد داشت و اگر ۰ باشد، حضور نخواهد داشت.

$y_i$: متغیر تصمیم (صفر و یک) حضور یا عدم حضور گره $i\in V$ در یکی از تورهاست. اگر برابر با ۱ باشد گره $i$ در تور حضور خواهد داشت و اگر ۰ باشد، حضور نخواهد داشت.

$K$: تعداد وسایل نقلیه موجود

$S$: یک مجموعه از گره‌ها به طوری که $S\subset V/\{0\}$

$\delta(s)$: مجموعه مسیرهایی که از مجموعه $S$ خارج می‌شوند یا به آن وارد می‌شوند.

$\delta(s)=\{(i,j)|(i,j)\in A, i\in S, j\notin S\}$

$\mu(s)$: کاردینالیتی مجموعه‌ی $\{k|C_k\subseteq S,\forall k\in\{1,2,...,m\}\}$

$\xi$: بردار متغیرهای تصادفی مرتبط با حضور احتمالی خوشه‌ها

$$(p)\ \text{Min}\sum_i\sum_{j>i}d_{ij}x_{ij}-E_\xi[Q(x,\xi)] \qquad (3-1)$$

S.T.

$$\sum_{j\in C_k}y_j\geq 1 \qquad 1\leq k\leq m \qquad (3-2)$$

$$\sum_{i<t}x_{it}+\sum_{j>t}x_{tj}\geq 2y_t \qquad 0\leq t\leq n \qquad (3-3)$$

---

1- Cluster

$$\sum_{i \in V/\{\circ\}} x_{i\circ} = k \qquad (3\text{-}4)$$

$$\sum_{j \in V/\{\circ\}} x_{\circ j} = k \qquad (3\text{-}5)$$

$$\sum_{\substack{(i,j) \in \delta(s) \\ i<j}} x_{ij} + \sum_{\substack{(i,j) \in \delta(s) \\ j<i}} x_{ji} \geq 2 y_t \qquad \mu(V/S) \neq \circ \text{ و } \forall \mu(s) \neq \circ \qquad (3\text{-}6)$$

$$y_\circ = k \qquad (3\text{-}7)$$

$$x_{ij} = \{\circ, 1\} \qquad \circ \leq i < j \leq n \qquad (3\text{-}8)$$

$$y_j = \{\circ, 1\} \qquad \circ \leq j \leq n \qquad (3\text{-}9)$$

در این مدل محدودیت‌های (3-2) تضمین کننده این است که هر خوشه حداقل یک بار ملاقات خواهد شد. محدودیت‌های (3-3) تضمین کننده ارتباط بین تورها خواهد بود. البته این ارتباط تنها همراه با محدودیت‌های حذف زیر دور[1] (3-6) معنی‌دار خواهد بود. محدودیت (3-7) برای اعمال حضور قطعی گره مبدأ (دپو) لازم است. محدودیت‌های (3-4) و (3-5) نیز برای اعمال محدودیت استفاده از k عدد وسیله‌ی نقلیه و در نتیجه ساخت k تور لازم هستند. البته محدودیت‌های (3-2) را می‌توانستیم به صورت $\sum_{j \in C_k} y_i = 1, 1 \leq k \leq m$ نیز بنویسیم. چرا که با وجود شرط نامساوی مثلث، هر خوشه تنها یک بار ملاقات خواهد شد.

در این مدل از تابع هدف مشابهت زیادی به تابع هدف دومرحله‌ای در مقاله [37] [40] [1994] Laporte et al برای مسأله PTSP دارد. برای تشریح بیشتر این مدل و طراحی الگوریتم‌های حل دقیق و ابتکاری نیاز به تعاریف زیر داریم:

**تور توالی غیرقطعی[2]:** در واقع جواب بهینه مدل بالا، به ما مجموعه‌ای از تورهای توالی غیرقطعی می‌دهد. این نوع تور با تورهای عادی تفاوت دارد. تورهای توالی غیرقطعی، تنها بیانگر توالی گره‌هایی هستند که در خوشه بایستی ملاقات شوند. ملاقات قطعی هر گره (یا هر خوشه) موکول به معلوم شدن این است که آیا خوشه محتوی گره مورد نظر در زمان اجرای تور حضور داشته باشد یا نه. اگر این پیشامد (که احتمال حضور آن $p_i$ بود) رخ دهد این گره ملاقات خواهد شد و در غیر اینصورت مسیر واقعی طی شده به طور مستقیم از گره قبل از i به گره بعد از i در تور توالی غیرقطعی تغییر خواهد کرد.

---

1- Subtour eliminator
2- a priori tour

۵۰

عملاً تور توالی غیرقطعی تنها بیانگر اولویت ما در توالی گره‌هاست و تور واقعی که رخ می‌دهد بستگی مستقیم به این توالی و همچنین پیشامد رخ داده دارد.

از این رو طول مسیر واقعی با مسیر تور توالی غیرقطعی با هم متفاوت هستند و تنها می‌توان از امید ریاضی تور توالی غیرقطعی سخن گفت. لذا در این رساله "L" بیانگر یک تور توالی غیرقطعی خواهد بود. "$\ell$" طول این تور در حالت قطعی است و $E(L)$ امید ریاضی تور توالی غیر قطعی.

براین اساس و همچنین بنابر تعریف تابع مقدار مرحله دوم در بهینه‌سازی احتمالی دو مرحله‌ای، می‌توان تابع هدف را تشریح کرد. $\sum_i \sum_{j>i} d_{ij} x_{ij}$ مجموع طول مسیرها در حالتی است که تمام خوشه‌ها حضور داشته باشند. همچنین $Q(x, \xi)$ مقدار کاهش مجموع طول مسیرهاست، در صورتی که $\xi$ رخ دهد. لذا $E(Q(x, \xi))$، بیانگر امید ریاضی این کاهش خواهد بود.

با توجه به اینکه هدف ما حداقل کردن امید ریاضی طول تورهای توالی غیرقطعی یا همان $E(L_{PGVRP})$ است، لذا $\sum_i \sum_{j>i} d_{ij} x_{ij} - E_\xi(Q(x, \xi))$ تابع هدف را تشکیل خواهد داد. به عبارت دیگر امید ریاضی کاهش صورت گرفته یعنی $E_\xi(Q(x, \xi))$ برابر با $E_\xi(L) - \sum_i \sum_{j>i} d_{ij} x_{ij}$ می‌باشد. البته باید توجه داشت که $E_\xi(Q(x, \xi))$ در صورتی که نامساوی مثلث برقرار باشد، همواره غیرمنفی خواهد بود.

قضیه ۱- فرض کنید تورهای توالی غیرقطعی L به صورت $(0, 1, 2, ..., h, 0)$ باشد. در آن صورت
$$E(L) = \sum_{i=0}^{h} \sum_{j=0}^{h} \alpha_{ij} d_{ij}$$
اگر $\alpha_{ij} = p_i p_j \prod_{t=i+1}^{j-1}(1-p_t)$. ( $p_j$ احتمال حضور گره $j$ ام است. این احتمال برابر با احتمال خوشه‌ای است که محتوی گره $j$ ام است)

توجه داشته باشید تورهای توالی غیرقطعی مجموعه k تور را تشکیل می‌دهد که هر یک از آنها را یک وسیله‌ی نقلیه طی خواهد کرد. به همین ترتیب، نمایش L، این مجموعه‌ی تورها را بدون اولویت بر هم پشت سر هم قرار می‌دهد. همچنین گره مبدأ یا دپو احتمال حضور برابر ۱ داشته و عضو هیچ خوشه‌ای نیست.

۵۱

# فصل چهارم

# طراحی الگوریتم‌های مناسب (دقیق و ابتکاری) برای حل مسأله مسیریابی توسعه یافته احتمالی

**مقدمه**

بنابر آنچه در فصول پیش آمد، مسأله مسیریابی توسعه یافته احتمالی مدل‌سازی شده و تکنیک مناسب برای استفاده در طراحی الگوریتم حل دقیق (روش L-Shape) معرفی گردید. در این فصل با استفاده از روش L-Shape الگوریتم حل دقیق پیشنهاد شده است. همچنین برای حل مسائل بزرگ، الگوریتم حل ابتکاری مناسب نیز ارائه گردیده است.

## ۴-۱- طراحی الگوریتم دقیق

PGVRP مدل سازی شده در بخش قبل در واقع توسعه‌ای است از GVRP و PVRP. مسأله GVRP مسأله‌ای کاملاً قطعی است که در آن احتمال حضور هر خوشه برابر ۱ است و مسأله PVRP مسأله‌ای است احتمالی که معادل PGVRP با خوشه‌های تک عضوی به حساب می‌آید. همچنین PGVRP توسعه‌ای از PGTSP نیز به حساب می‌آید، چرا که اگر در PGVRP تعداد وسائل نقلیه در دسترس تنها یک عدد باشد، عملاً با PGTSP مواجه خواهیم بود.

الگوریتم پیشنهادی برای حل این مسأله بطور دقیق، در واقع توسعه‌ای است از روش L-Shape عدد صحیح پیشنهاد شده توسط [۴۱] [۱۹۹۳] Laporte و Louveaux. همچنین ابزارهایی که برای بهبود این

روش توسط [۶۸] Tang [۲۰۰۵] و Hooks پیشنهاد شده نیز برای حل مسأله‌ی PGVRP توسعه داده شده است. از این رو در روش شاخه و کران از کران‌های بالا[1] بهبود داده شده برای کاهش مسافت مورد انتظار (امید ریاضی) در این مسأله و برای هر گره درخت استفاده خواهد شد تا بدین وسیله امکان حذف شاخه‌های بیشتری وجود داشته باشد. بدین وسیله تعداد برش‌های بهینه‌یاب نیز در الگوریتم شاخه و برش[2] کاهش خواهد یافت.

### ۴-۱-۱- کران پایین برای PGVRP

[۳۲] Jaillet [۱۹۸۷] برای مسأله فروشنده دوره‌گرد احتمالی، کران پایینی به صورت زیر طراحی کرده است. L تور توالی غیرقطعی بهینه برای PTSP بوده و $\ell$ تور بهینه برای TSP روی همان گراف می‌باشد. همچنین احتمال حضور همه گره‌ها مساوی و برابر با P فرض شده است. $E(L_{PTSP}) \geq p.\ell_{TSP}$

بر همین اساس به سادگی کران پایین مزبور برای مسأله PGVRP با احتمال حضور خوشه‌های مساوی و برابر با p قابل توسعه است. L نشان‌دهنده‌ی تورهای توالی غیرقطعی بهینه برای PGVRP بوده و $\ell$ تور بهینه برای GVRP روی همان گراف است.

$E(L_{PGVRP}) \geq p.\ell_{GVRP}$

این کران برای حالتی که احتمال حضور خوشه‌ها مساوی نباشد به صورت زیر قابل توسعه است.

$E(L_{PGVRP}) \geq \underset{1 \leq i \leq m}{\text{Min}} \{p_i\} \ell_{GPVRP}$

با وجود سادگی تعریف بالا، بدلیل NP-Hard بودن مسأله GVRP و در نتیجه سخت بودن محاسبه کران پایین استفاده از آن در همه‌ی مسائل توصیه نمی‌شود.

همچنین با توجه به اینکه این کران پایین مستقل از سایر مراحل الگوریتم بوده و جواب‌های حاصل در هر شاخه از الگوریتم شاخه و کران توان بهبود آن را ندارند، لذا تنها در مواردی که ابعاد مسأله PGVRP بسیار بزرگ بوده و حالت قطعی آن بعد به مراتب پایین‌تری داشته باشد، استفاده از این کران پایین توصیه می‌شود.

برای طراحی کران پایین $E(L)$ می‌توان به شکلی دیگر عمل کرد. با توجه به اینکه

$E(L) = \sum_i \sum_{j>i} d_{ij} x_{ij} - E_\xi [Q(x, \xi)]$

---

1- upper bound
2- Branch and cut

٥٤

در صورتی که بتوانیم برای $E_\xi(Q(x,\xi))$ کران بالایی پیدا کنیم. برای $E(L)$ کران پایین در دسترس خواهد بود. لذا فرض کنید $\theta = E_\xi[Q(x,\xi)]$

همچنین U را به صورت زیر تعریف می‌کنیم:

$$U = \sum_{\substack{i,j \in \bar{C}(t) \\ i \in C(j) \\ j \in C(i)}} \max \{(1-p_t)(d_{it} + d_{tj} - d_{ij})\} \tag{4-1}$$

اگر $\bar{C}(t)$ را به عنوان مکمل خوشه‌ی $C(t)$ در نظر بگیریم، مقدار U یک کران بالای مناسب برای $\theta$ ایجاد می‌کند بطوریکه

$$0 \le \theta \le U \tag{4-2}$$

U با جمع حداکثر کاهش‌های متصور در فاصله‌ی طی شده با حذف هر خوشه از تور توالی غیرقطعی ساخته می‌شود. برای هر خوشه کاهش در فاصله طی شده برای تمام گره‌های ممکن ورودی به و خروجی از خوشه محاسبه شده و حداکثر آنها مشخص می‌شود. پس از ضرب در احتمال عدم حضور خوشه‌ی مورد نظر، مجموع برای همه‌ی خوشه‌ها (گره‌های تشکیل دهنده‌ی تور توالی غیرقطعی) محاسبه می‌شود. مشخصاً این مقدار از هر مقدار ممکن دیگر برای کاهش در فاصله‌ی طی شده بیشتر است. لیکن با وجود سادگی محاسباتی، کران بالای مؤثرتری نیز قابل محاسبه است. این کران بالا به صورت زیر قابل محاسبه است:

$$U = \sum_{k=1}^{m} u_k \tag{4-3}$$

در جایی که

$$u_k = \max_{t \in C_k} [\max_{\substack{i,j \in \bar{C}(t) \\ i \notin C(j) \\ j \notin C(i)}} [(1-p_t)(d_{it} + d_{tj} - d_{ij})]] \tag{4-4}$$

این کران بالا با وجود برآورد دقیق‌تری که از کران بالا محاسبه می‌کند، پیچیده‌تر از معادله (4-2) است. معادله (4-4) در هر خوشه ابتدا به ازای هر گره عضو خوشه حداکثر کاهش متصور را محاسبه می‌کند و سپس مقادیر حاصل از خوشه‌ها را ($u_k$) با هم جمع می‌نماید.

کران بالای دیگری که توسط [37] [40] [1994] Laporte et al پیشنهاد شده برای این مسأله بصورت زیر توسعه می‌یابد. با در نظر گرفتن مسافت طی شده‌ی مورد انتظار برای هر مسیر $(i,j)$ که ممکن است در تور توالی غیرقطعی قرار گیرد به طور خاص، برای هر $(i,j)$، $b_{ij}$ به عنوان کران پایین بر مقدار فاصله‌ی مورد انتظار که واقعاً در مرحله‌ی دوّم در راستای مسیرهایی که حداقل یک گره مشترک با $(i,j)$

۵۵

داشته باشند در نظر گرفته می‌شود. این مقدار برای هر (i , j) به صورت زیر قابل محاسبه است. (البته برای $i \neq 0$ و $j \neq 0$)

$$b_{ij} = p_i p_j d_{ij} + \frac{1}{2} p_i (1-p_j) \min_{\substack{t \notin C(j) \\ t \in C(i)}} \{d_{it}\} + \frac{1}{2} (1-p_i) p_j \min_{\substack{t \notin C(i) \\ t \in C(j)}} \{d_{tj}\} \qquad (4-5)$$

با توجه به احتمال حضور ۱ برای مبدأ، برای $i=0$ یا $j=0$، $b_{ij} = p_i p_j d_{ij}$

بدین ترتیب با تشکیل ماتریس B با عناصر $b_{ij}$ و D با عناصر $d_{ij}$ و X با عناصر $x_{ij}$، نامساوی

$$\theta \leq (D-B) X \qquad (4-6)$$

یک کران بالا برای مسأله ایجاد می‌کند.

باید توجه کرد که جنس این کران بالا با ۲ مورد قبلی متفاوت است. چرا که از متغیرهای تصمیم در آن استفاده شده و لذا براساس تور توالی غیرقطعی معین نمی‌شود.

### ۴-۲-۱-۴- الگوریتم دقیق

مراحل الگوریتم دقیق در ادامه ارائه خواهد شد. همچنین شرح هر مرحله با جزئیات بیشتر پس از آن آمده است.

**الف- گام‌های الگوریتم PGVRP**

**گام ۱-** ایندکس هر تکرار (q) را برابر با صفر قرار دهید. $q:=0$.

لیست زیر مسأله‌ها را تشکیل داده و تنها مسأله اولیه را در لیست قرار دهید. یک جواب شدنی را با استفاده از الگوریتم ابتکاری ایجاد کرده و جواب آن را در $\overline{Z}$ به عنوان بهترین جواب در دسترس قرار دهید.

**گام ۲-** اگر لیست زیر مسأله‌ها خالی است، توقف کنید، در غیراینصورت آخرین زیر مسأله‌ی لیست را انتخاب نمایید.

**گام ۳-** مقدار q را یک واحد افزایش دهید. ($q:=q+1$). q امین زیر مسأله را حل کرده و جواب آن را در چند تایی مرتب ($X^q, Y^q, \theta^q$) قرار دهید.

**گام ۴-** اگر $\bar{Z} \geq \theta^q - DX^q$، مسأله را حذف کرده و به گام ۲ باز گردید. در غیر اینصورت نامساوی‌هایی که برقرار نیستند[1] را شناسایی کرده و حداقل یکی از آنها را به مسأله جاری اضافه نمایید. به گام ۳ باز گردید.

**گام ۵-** بررسی کنید که آیا جواب $Y^q, X^q$ عدد صحیح هستند یا نه. اگر حداقل یک متغیر عدد صحیح نبود، با شاخه زدن[2] بر روی متغیر مزبور ۲ زیر مسأله‌ی جدید ایجاد کنید. برای هر مسأله، کران بالای (۴-۲) را مورد بررسی قرار داده و اگر ممکن بود $\theta$ را بر آن اساس بهبود دهید. زیر مسأله‌های جدید را به لیست زیر مسأله‌ها اضافه کنید. به گام ۲ باز گردید.

**گام ۶-** کاهش فاصله‌ی مورد انتظار ($Q(X^q)$) را با استفاده از جواب عدد صحیح $X^q$ محاسبه کنید و مقدار $Z^q$ را به صورت زیر محاسبه کنید.

$$Z^q := DX^q - Q(X^q) \qquad (۴-۷)$$

اگر $Z^q < \bar{Z}$، آنگاه مقدار $Z^q$ را در $\bar{Z}$ قرار دهید ($\bar{Z} := Z^q$)

**گام ۷-** اگر $\theta^q \leq Q(X^q)$، زیر مسأله‌ی جاری را حذف کنید و به گام ۲ باز گردید. در غیر اینصورت یک برش بهینه یاب به مسأله اضافه کرده و به گام ۳ باز گردید.

الگوریتم دقیق ۲ بخش اصلی دارد. بخش اول پروسه‌ی شاخه و برش است (گام ۱ تا ۵) و بخش دوم اضافه کردن برش بهینه یاب (گام‌های ۶ و ۷). در ابتدا لیست مسائل که تنها شامل مسأله‌ی اولیه است تشکیل می‌شود. شاخه و برش مجموعه‌ای از مسائل خطی ساده‌سازی شده[3] را ایجاد کرده، به لیست زیر مسأله‌ها اضافه کرده و آنها را حل می‌نماید تا به جوابی عدد صحیح برسد سپس مقدار واقعی تابع مقدار[4] محاسبه شده و برآن اساس برش بهینه یاب شناسایی شده و در صورت نیاز به مسأله‌ی ساده‌سازی شده اضافه می‌شود.

الگوریتم زمانی توقف می‌کند که لیست زیر مسأله‌ها خالی شده باشد.

---

1- Violated inequalities
2- branching
3- Relaxed
4- recourse function

**ب- مسأله‌ی اولیه**

الگوریتم دقیق با حل مسأله اولیه که در واقع ساده‌سازی شده مسأله‌ی p است آغاز می‌شود. مسأله به صورت زیر خواهد بود.

$$\text{Min} \sum_i \sum_{j>i} d_{ij} x_{ij} - \theta \qquad (4-8)$$

S.T.

$$\sum_{j \in C_k} y_j \geq 1 \qquad 1 \leq k \leq m \qquad (4-9)$$

$$\sum_{i<t} x_{it} + \sum_{j>t} x_{tj} \geq 2 y_t \qquad 0 \leq t \leq n \qquad (4-10)$$

$$\sum_{i \in V/\{0\}} x_{i0} = k \qquad (4-11)$$

$$\sum_{j \in V/\{0\}} x_{0j} = k \qquad (4-12)$$

$$y_0 = k \qquad (3-7)$$

$$\theta \leq (D-B)X \qquad (4-6)$$

$$0 \leq \theta \leq U \qquad (4-2)$$

$$0 \leq x_{ij} \leq 1 \qquad 0 \leq i < j \leq n \qquad (4-13)$$

$$0 \leq y_i \leq 1 \qquad 0 \leq i \leq n \qquad (4-14)$$

در این مسأله محدودیت‌های حذف کننده‌ی زیردورها و عدد صحیح بودن متغیرها، ساده‌سازی شده است. همچنین عبارت کاهش فاصله‌ی مورد انتظار از تابع هدف مسأله (4-2) با $\theta$ (که البته در محدودیت (4-2) با U محدود شده) جایگزین شده است.

محدودیت (4-6) کران بالایی به صورت تابع برای $\theta$ ایجاد می‌کند که اساساً برای ایجاد مقدار $\theta$ در (4-8) مورد استفاده قرار می‌گیرد. در طول الگوریتم این مقدار با ایجاد برش‌های بهینه‌یاب بهبود می‌یابد تا جایی که $\theta$ بیانگر مقدار دقیق کاهش فاصله‌ی مورد انتظار گردد.

۵۸

### ج- بررسی نامساوی‌هایی که برقرار نیستند

نامساوی‌های مجازی که می‌توان برای GVRP به کار برد، برای GPVRP نیز مجاز خواهند بود. با وجود اینکه انواع نامساوی‌های مختلفی می‌توان برای GVRP معرفی کرد، در این رساله تنها از نامساوی‌های حذف زیر دورها در این مرحله استفاده شده است. برای ایجاد نامساوی‌های حذف زیر دور از الگوریتم GSEC-H2 که توسط Fichetti et al [۱۹۹۷] [۱۵] معرفی شده است استفاده شده است.

### د- برش بهینه‌یاب

برش بهینه‌یاب پیشنهادی برای استفاده در گام ۷ به صورت زیر است:

$$\theta \leq U + (Q(X^q) - U)\left(\sum_{(i,j)\in E^q} x_{ij} - |E^q| + 1\right) \qquad (4\text{-}15)$$

$$E^q = \{(i,j) \in A \mid x_{ij}^q = 1\}$$

برای نشان دادن اینکه این برش بهینه یاب یک برش مجاز است فرض کنید که در گام ۷، $X = X^q$ باشد. در آن صورت $\sum_{(i,j)\in E^q} x_{ij} = |E^q|$ و در نتیجه این نامساوی به شکل $\theta \leq Q(X^q)$ در می‌آید که مجاز است. در صورتی که $X \neq X^q$ نباشد، $\sum_{(i,j)\in E^q} x_{ij} \leq |E^q| - 2$ خواهد بود و

$$U + (Q(X^q) - U)\left(\sum_{(i,j)\in E^q} x_{ij} - |E^q| + 1\right) \geq 2U - Q(X^q) \geq \theta$$

که این هم مجاز است.

## ۴-۲- الگوریتم‌های ابتکاری

از آنجا که مسأله TSP یک مسأله‌ی NP-Hard است، مسأله PGVRP که توسعه یافته این مسأله و در نتیجه پیچیده‌تر از آن است نیز NP-Hard خواهد بود. لذا حل این مسأله در ابعاد بالا با استفاده از الگوریتم حل دقیق در زمان منطقی امکان‌پذیر نخواهد بود. لذا در کنار طراحی الگوریتم دقیق نیاز به طراحی الگوریتم ابتکاری مناسب برای مسائل با ابعاد بالا وجود دارد. الگوریتم‌های ابتکاری پیشنهادشده برای تمام حالت‌های مسائل مسیریابی، توسعه یافته الگوریتم‌های ابتکاری برای حل مسأله‌ی TSP هستند. همچنان که در فصل اول اشاره شد، الگوریتم ابتکاری برای حل مسأله‌ی TSP مبتنی بر مفهوم «ارزش

قرار دادن»[1] طراحی شده‌اند. این الگوریتم‌ها در هر مرحله با محاسبه‌ی ارزش قرار دادن یک گره ملاقات نشده در تور ناقص و مقایسه‌ی این مقادیر با هم، تصمیم می‌گیرد که کدام گره را به کدام نقطه‌ی تور ناقص اضافه کرده و تور جدید را ایجاد کند. این کار تا کامل شدن تور (شامل شدن همه‌ی گره‌ها در تور) ادامه می‌یابد. مسأله اصلی در طراحی الگوریتم‌ها در وهله‌ی اول تعیین معیاری مناسب برای محاسبه‌ی «ارزش قرار دادن» خواهد بود و در وهله‌ی دوم چگونگی استفاده از این معیار در اضافه کردن گره جدید به تور ناقص است.

در مسأله TSP ارزش قرار دادن گره $t$ در میان گره‌های $i$ و $j$ به سادگی به صورت $d_{it} + d_{tj} - d_{ij}$ محاسبه می‌شود. این مقدار در حالتی که مسأله PGVRP باشد این قدر ساده نیست. در واقع ما در این مسأله با امید ریاضی ارزش قرار دادن برای مسأله‌ی PGVRP مواجه هستیم. برای محاسبه‌ی امید ریاضی ارزش قرار دادن $j$ در تور توالی غیر قطعی ($0, 1, 2, ..., i-1, i, ..., m-1, m, 0$) در حالتی که تنها یک وسیله‌ی نقلیه و در نتیجه یک تور مدنظر است می‌توان از فرمول زیر استفاده کرد.

$$\text{Extra}[j] = \left\{ \sum_{r=0}^{i-1} \left\{ p_r \prod_{t=r+1}^{i-1} (1-p_t) \left[ \sum_{z=i}^{m} p_z (d_{rj} + d_{jz} - d_{rz}) \prod_{q=i}^{z-1} (1-p_q) \right] \right\} \right\} p_j \qquad (4\text{-}16)$$

در جایی که $z = i$ باشد، $\prod_{q=i}^{z-1}(1-p_q) = 1$ خواهد بود. همچنین $p_j$ برای $j \in V$ برابر با احتمال حضور خوشه‌ای است که محتوی گره $j$ می‌باشد.

با استفاده از ارزش مورد انتظار قرار دادن، می‌توان الگوریتم‌های ابتکاری پیشنهادی برای مسأله فروشنده دوره گرد را به مسأله PGVRP توسعه داد. الگوریتم ابتدایی که پیشنهاد می‌شود در هر مرحله از الگوریتم قرار دادن، ارزش قرار دادن همه‌ی گره‌های ملاقات نشده در تمام نقاط ممکن تور ناقص را محاسبه کرده و کمترین را انتخاب نموده و گره منتخب را در مکان منتخب قرار می‌دهد و به مرحله بعد می‌رود. این کار تا تکمیل تور ادامه می‌یابد. برای توسعه این الگوریتم به مسأله‌ی PGVRP گام‌های زیر پیشنهاد می‌شود.

---

1- insertion value

## ۴-۲-۱- الگوریتم mmI برای حل PGVRP (با وسائل نقلیه ظرفیت محدود)

**گام ۱**- ایندکس وسیله‌ی نقلیه را برابر با ۱ قرار دهید ($m:=1$)

**گام ۲**- مسیر به طول ۲ گره تشکیل داده و آخرین و اولین گره مسیر را گره مبدأ (دپو) قرار می‌دهیم.

**گام ۳**- تمامی گره‌هایی را که در خوشه‌های ملاقات نشده قرار دارند، در تمام نقاط ممکن در تور ناقص قرار داده و مقدار مورد انتظار ارزش قرار دادن متناظر را محاسبه کنید.

براساس کمترین مقدار مورد انتظار ارزش قرار دادن، گره مناسب و محل قرار دادن مناسب آن در تور ناقص را تعیین کرده و تور را بازسازی می‌کنیم. تعداد گره در تور را یک واحد افزایش می‌دهیم ($mc:=mc+1$) در صورتی که همه‌ی خوشه‌ها ملاقات شده باشند تور ساخته شده را به عنوان تور وسیله‌ی نقلیه m ام معرفی کرده و توقف می‌کنیم.

اگر ظرفیت خودرو m تکمیل شده باشد ایندکس وسیله‌ی نقلیه را یک واحد افزایش می‌دهیم ($mc:=mc+1$)، تور ساخته شده را به عنوان تور وسیله‌ی نقلیه m ام معرفی کرده و به گام ۲ باز می‌گردیم.

براساس مطالعات Stewart [۲۵] [۱۹۸۵] و Golden، مشخص شده است که درمسأله فروشنده دوره گرد، در صورتی که معیار تصمیم‌گیری بزرگترین ارزش قرار دادن باشد، نتایج حاصل بهتر خواهند بود. با استفاده از نتایج این مطالعه و تعمیم آن به مسأله PGVRP، الگوریتم‌های ابتکاری متنوعی را می‌توان توسعه داد. اصلی‌ترین الگوریتم قابل ارائه که در این رساله پیشنهاد شده است، مبتنی بر الگوریتم (Max-min-Insertion) MmI می‌باشد.

## ۴-۲-۲- الگوریتم MmI برای حل PGVRP (با وسائل نقلیه ظرفیت محدود)

**گام ۱**- ایندکس وسیله‌ی نقلیه را برابر با ۱ قرار دهید ($m:=1$)

**گام ۲**- مسیر به طول ۲ گره تشکیل داده و آخرین و اولین گره مسیر را مبدأ (دپو) قرار دهید.

**گام ۳**- خوشه‌های ملاقات نشده را در نظر بگیرید. در هر خوشه، گره‌های موجود را در تمام نقاط ممکن تور ناقص قرار داده و کمترین مقدار مورد انتظار ارزش قرار دادن را در آن خوشه محاسبه کنید. مقادیر حاصل برای همه‌ی خوشه‌های ملاقات نشده را در نظر گرفته و بیشترین آنها را انتخاب کنید. گره و نقطه‌ی قرار دادن در تور ناقص مربوط به مقدار به دست آمده را در نظر گرفته و مطابق آن گره را در نقطه‌ی مورد نظر قرار دهید و تور را بازسازی کنید.

تعداد گره در تور را یک واحد افزایش می‌دهیم ($mc := mc+1$) در صورتی که همه‌ی خوشه‌ها ملاقات شده باشند تور ساخته شده را به عنوان تور وسیله‌ی نقلیه $m$ ام معرفی کرده و توقف می‌کنیم.

اگر ظرفیت خودرو $m$ تکمیل شده باشد ایندکس وسیله‌ی نقلیه را یک واحد افزایش می‌دهیم ($mc := mc+1$)، تور ساخته شده را به عنوان تور وسیله‌ی نقلیه $m$ ام معرفی کرده و به گام ۲ باز می‌گردیم.

الگوریتم‌های پیشنهادی در مرحله‌ی قبل برای حل مسأله مسیریابی توسعه یافته احتمالی با وسایل نقلیه با ظرفیت محدود و برابر ممکن است. برای حل این مسأله در حالتی که مسأله بطور دقیق حل شده است الگوریتم mmI به صورت زیر بازنویسی شده و پیشنهاد می‌گردد.

الگوریتم‌های پیشنهادی در مرحله قبل مسأله مسیریابی توسعه‌یافته احتمالی را در حالتی که وسائل نقلیه ظرفیت محدود دارند حل می‌کند. با توجه به اینکه الگوریتم حل دقیق برای حل مدل مسیریابی توسعه یافته احتمالی با وسائل نقلیه نامحدود ظرفیت طراحی شده، الگوریتم ابتکاری جدیدی، مبتنی بر «ارزش قرار دادن» پیشنهاد می‌شود. بدین وسیله در فصل بعد نتایج حاصل از بکار بردن این دو الگوریتم قابل مقایسه و ارزیابی خواهد بود. گام‌های الگوریتم مزبور به شرح زیر پیشنهاد می‌گردد.

### ۴-۳-۲-۴- الگوریتم ابتکاری برای حل PGVRP (با وسائل نقلیه ظرفیت نامتناهی)

گام ۱- ایندکس وسیله نقلیه را برابر با ۱ قرار دهید. ($m := 1$)

گام ۲- مسیر به طول ۲ گره تشکیل داده و آخرین و اولین گره مسیر را مبدأ (دپو) قرار دهید.

گام ۳- تمامی گره‌هایی را که در خوشه‌های ملاقات نشده قرار دارند، در تمام نقاط ممکن در تور ناقص قرار داده و مقدار مورد انتظار ارزش قرار دادن متناظر را محاسبه کنید. براساس کمترین مقدار مورد انتظار ارزش قرار دادن، گره مناسب و محل قرار دادن مناسب آن در تور ناقص را تعیین کرده و تور را بازسازی می‌کنیم. اگر ایندکس وسیله‌ی نقلیه برابر $k$ (تعداد وسائل نقلیه) نیست، یک واحد به ایندکس وسیله‌ی نقلیه اضافه کرده و به گام ۲ بروید، در غیراینصورت به گام ۴ بروید.

گام ۴- تمام گره‌هایی را که در خوشه‌های ملاقات نشده قرار دارند، در تمام نقاط ممکن در تور ناقص همه‌ی وسائل نقلیه قرار داده و مقدار مورد انتظار ارزش قرار دادن متناظر را محاسبه کنید. براساس کمترین مقدار مورد نظر ارزش قرار دادن، وسیله‌ی نقلیه، گره و محل قرار دادن در تور ناقص را مشخص کرده و تور ناقص وسیله‌ی نقلیه مزبور را بر این اساس به روز کنید. در صورتی که همه‌ی خوشه‌ها ملاقات شده‌اند تورهای وسائل نقلیه را به عنوان تورهای کامل معرفی کرده و توقف کنید، در غیراینصورت گام ۴ را تکرار کنید.

# فصل پنجم

# نتایج محاسباتی به کارگیری الگوریتم‌های پیشنهادی برای حل مسأله‌ی مسیریابی توسعه یافته احتمالی

**مقدمه**

در این فصل با به کارگیری الگوریتم حل دقیق پیشنهادی و همچنین الگوریتم ابتکاری پیشنهادی، مدل ارائه شده در فصل سوّم در ابعاد مختلف حل شده است. شیوه تولید مسائل نمونه به صورتی بوده است که گستره مسائل کوچک تا بزرگترین مسائل حل شده تاکنون و برای موارد مشابه را دربرگیرد. همچنین سعی شده است قابلیت‌های الگوریتم‌های پیشنهادی با حل این مسائل به آزمون گذارده شود همچنان که در بخش بعد، نتایج محاسباتی بطور خلاصه در جدولی ارائه می‌شود، تحلیل نتایج بدست آمده در نمودارها و جداول دیگر ارائه خواهد شد.

## ۵-۱- روش‌شناسی اعتبارسنجی الگوریتم‌های حل (دقیق و ابتکاری)

برای اعتبارسنجی الگوریتم‌های حل ۲ معیار مورد نظر خواهد بود:

۱- کیفیت جواب

۲- مؤثر بودن الگوریتم

در این رساله ۲ الگوریتم برای حل مسأله پیشنهاد شده است. الگوریتم اوّل، الگوریتم حل دقیق مبتنی بر روش شاخه و کران و استفاده از تکنیک L-Shape است. با توجه به اثبات همگرایی این الگوریتم بهینه بودن جواب الگوریتم نیاز به اعتبارسنجی مجدد ندارد. برای بررسی مؤثر بودن الگوریتم، زمان اجرای الگوریتم در مقایسه با زمان اجرای الگوریتم ابتکاری اعتبارسنجی خواهد شد. لازم به ذکر است، الگوریتم ابتکاری با استفاده

از زبان ++C برنامه‌نویسی و اجرا شده است. همچنین الگوریتم دقیق با استفاده از نرم‌افزار Cplex 10.1 برنامه‌نویسی شده و توسط Visual Studio 2005 اجرا شده است. با توجه به محدودیت این نرم‌افزار برای کاربردهای آموزشی و در دسترس نبودن نسخه تجاری آن، امکانات مختلف این نرم‌افزار برای بهینه‌سازی عملکرد الگوریتم شاخه و کران و همچنین بهینه‌سازی موازی[1] مورد استفاده قرار نگرفته است.

از این رو عملکرد الگوریتم دقیق با استفاده از نسخه آموزش این نرم‌افزار تفاوت جدی با نسخه‌ی تجاری آن دارد. برای سنجش عملکرد الگوریتم‌ها بایستی مسائل نمونه از ابعاد کوچک تا ابعاد بزرگ مسأله، ایجاد می‌شد. در مسأله مسیریابی توسعه یافته احتمالی پیچیدگی مسأله وابسته به ۲ پارامتر است.

پارامتر اوّل تعداد گره‌ها و پارامتر دوّم تعداد خوشه‌هاست. از این رو در طراحی مسائل نمونه ۱۶ مسأله تولید شد که در آنها ترکیب این ۲ فاکتور از ابعاد کوچک تا ابعاد بزرگ در نظر گرفته شود. تعداد گره‌ها حداقل ۱۰ گره بوده و به تدریج تا ۳۰۰ گره افزایش می‌یابد. همچنین همزمان با افزایش تعداد گره‌ها، تعداد خوشه‌ها (کلاسترها) نیز از ۲ خوشه تا ۱۵۰ خوشه افزایش می‌یابد. سایر فاکتورهای دخیل در تولید مسائل نیز که تأثیر کمتری در پیچیدگی مسأله دارند (تعداد خودروها) بطور منطقی تعیین شده است. همچنین برای همگونی مسائل، ظرفیت خوشه‌ها در هر مسأله بطور مساوی تقسیم شده است. همچنان که در جدول نتایج محاسباتی مشخص است، با افزایش تعداد گره به ۳۰۰ گره و ۱۰۰ خوشه، الگوریتم دقیق در زمان کمتر از ۱۰ ساعت قابل حل نبوده و در نتیجه به جای جواب بهینه در جدول نتایج محاسباتی مقدار کران پایین محاسبه و قرار داده شده است. الگوریتم ابتکاری پیشنهادی نیز برای حل مسائل مزبور به کار رفته و جواب‌های حاصل از به کارگیری این الگوریتم نیز در جدول آمده است. در رابطه با کیفیت جواب‌های حاصل از الگوریتم ابتکاری نیز در گستره‌ی مسائل نمونه، جواب‌های حاصل از الگوریتم با جواب بهینه یا کران پایین مقایسه خواهد شد. همچنین مؤثر بودن الگوریتم نیز با توجه به زمان اجرای الگوریتم برای حل مسائل قابل ارزیابی است.

## ۵-۲- نتایج محاسباتی

همچنان که در جدول مشخص است و در نمودارهایی که در پی آمده است قابل ارزیابی است، الگوریتم ابتکاری در اکثر موارد از دقت بسیار بالایی (خصوصاً در ابعاد پایین) در برآورد جواب بهینه برخوردار است. همچنین زمان اجرای الگوریتم دقیق نیز با وجود عدم استفاده از ابزارهای بهینه‌سازی نرم‌افزار Cplex (خصوصاً بهینه‌سازی موازی) در مواردی حتی از الگوریتم ابتکاری نیز بهتر است.

در ادامه جداول ارزیابی محاسباتی الگوریتم‌های پیشنهادی و مقایسه آنها از جنبه‌های مختلف آمده است.

---

1- Parallel Optimization

۶۵

جدول ۵-۱- نتایج محاسباتی به کارگیری الگوریتم‌ها برای حل مسائل نمونه

| شماره مسأله | تعداد گره | تعداد کلاستر | ظرفیت کلاسترها | تعداد وسیله نقلیه | جواب الگوریتم MMI (ابتکاری) | زمان انتظار برای اجرای MMI (ابتکاری) | جواب الگوریتم L-Shape (دقیق) | زمان انتظار برای اجرای الگوریتم L-Shape (دقیق) | انحراف |
|---|---|---|---|---|---|---|---|---|---|
| ۱ | ۱۰ | ۲ | ۵ | ۱ | ۲۰ | ۰ | ۲۰ | ۰/۰۴۶۸۷۵ | ۰ |
| ۲ | ۱۰ | ۵ | ۲ | ۲ | ۸۲ | ۰ | ۷۸ | ۰/۰۷۸۱۲۵ | ۰/۰۵۱۲۸ |
| ۳ | ۳۰ | ۵ | ۶ | ۲ | ۳۰ | ۰ | ۲۸ | ۰/۰۶۲۴۹۷ | ۰/۰۷۱۴۲ |
| ۴ | ۳۰ | ۱۰ | ۳ | ۵ | ۱۳۶ | ۰ | ۱۲۶ | ۰/۰۴۶۸۷ | ۰/۰۷۹۳۶ |
| ۵ | ۵۰ | ۱۰ | ۵ | ۴ | ۵۰ | ۰ | ۴۸ | ۰/۱۴۰۶۲۵ | ۰/۰۴۱۶۶ |
| ۶ | ۵۰ | ۲۵ | ۲ | ۷ | ۱۷۴ | ۰ | ۱۵۸ | ۳/۴۶۸۷۵ | ۰/۱۰۱۲۶ |
| ۷ | ۸۰ | ۲۰ | ۴ | ۶ | ۷۸ | ۱ | ۷۰ | ۰/۵۶۲۵ | ۰/۱۱۴۲۸ |
| ۸ | ۸۰ | ۴۰ | ۲ | ۱۲ | ۲۴۸ | ۱ | ۲۴۲ | ۰/۶۵۶۲۵ | ۰/۰۲۴۷۹ |
| ۹ | ۱۰۰ | ۲۵ | ۴ | ۸ | ۹۶ | ۰ | ۹۲ | ۱/۰۶۲۵ | ۰/۰۴۳۴۷ |
| ۱۰ | ۱۰۰ | ۵۰ | ۲ | ۱۴ | ۲۶۲ | ۲ | ۲۵۰ | ۱/۰۱۵۶۳ | ۰/۰۴۸ |
| ۱۱ | ۲۰۰ | ۵۰ | ۴ | ۱۵ | ۱۴۰ | ۲ | ۱۲۶ | ۴/۵۳۱۲۵ | ۰/۱۱۱۱ |
| ۱۲ | ۲۰۰ | ۱۰۰ | ۲ | ۱۸ | ۲۲۸ | ۲۵ | ۱۹۶ | ۱۳۱/۵ | ۰/۱۶۳۲۶ |
| ۱۳ | ۲۵۰ | ۵۰ | ۵ | ۱۸ | ۱۶۰ | ۳ | ۱۱۸ | ۲/۷۳۴ | ۰/۳۵۵۹۳ |
| ۱۴ | ۲۵۰ | ۱۲۵ | ۲ | ۲۵ | ۲۸۸ | ۹۴ | ۲۸۶ | ۵۰۶/۹۶۹ | ۰/۰۰۰۶۹ |
| ۱۵ | ۳۰۰ | ۱۰۰ | ۶ | ۲۰ | ۱۸۴ | ۴۶ | ۱۳۴ | L | ۰/۳۷۳۱۳ |
| ۱۶ | ۳۰۰ | ۱۵۰ | ۲ | ۳۰ | ۲۳۲ | ۳۲۱ | ۱۴۰ | L | ۰/۶۸۵۷۱ |

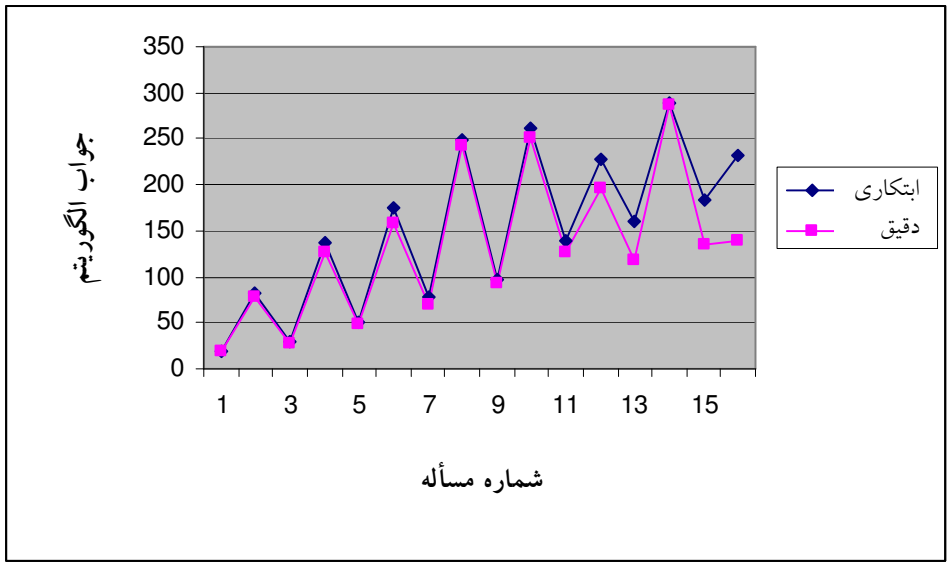

شکل ۵-۱- نمودار مقایسه‌ای پاسخ الگوریتم‌های دقیق و ابتکاری

برای حل مسائل نمونه

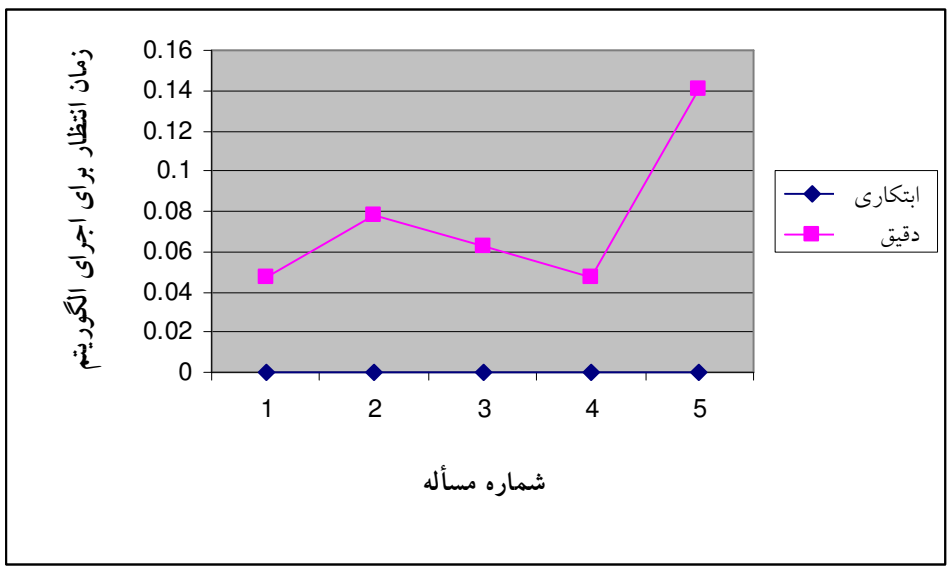

شکل ۵-۲- نمودار مقایسه‌ای زمان انتظار الگوریتم‌های دقیق و ابتکاری

برای حل مسائل نمونه (مسائل ۱ تا ۵)

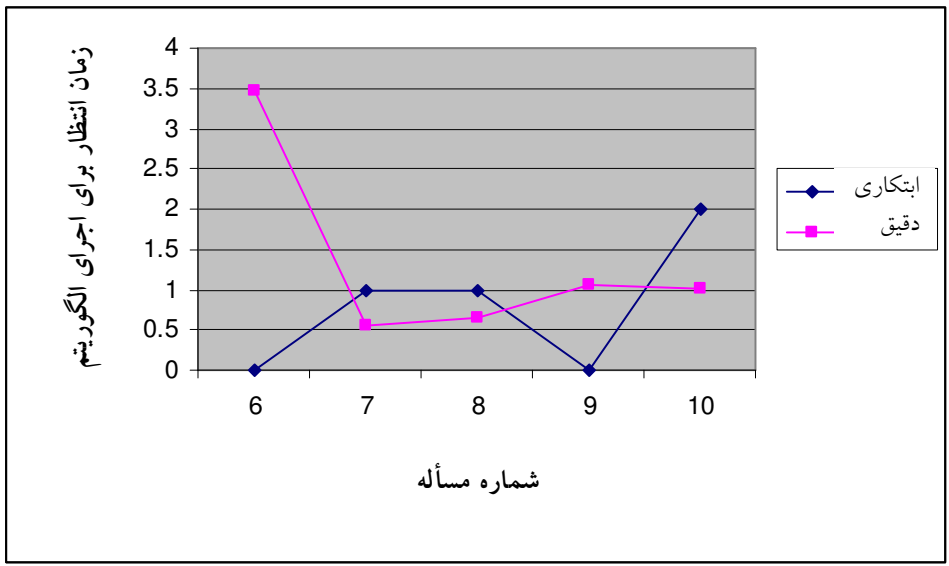

شکل ۵-۳- نمودار مقایسه‌ای زمان انتظار الگوریتم‌های دقیق و ابتکاری

برای حل مسائل نمونه (مسائل ۶ تا ۱۰)

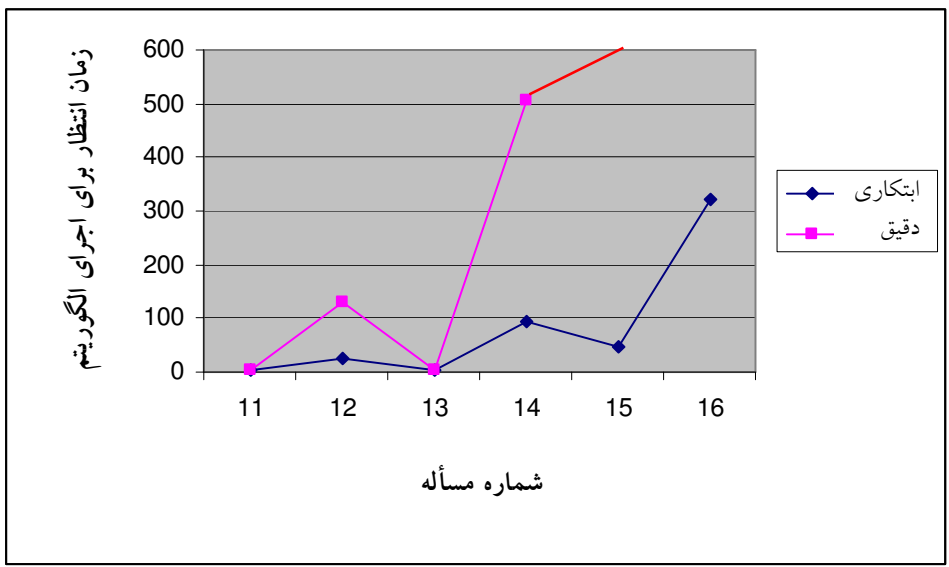

شکل ۵-۴- نمودار مقایسه‌ای زمان انتظار الگوریتم‌های دقیق و ابتکاری

برای حل مسائل نمونه (مسائل ۱۱ تا ۱۶)

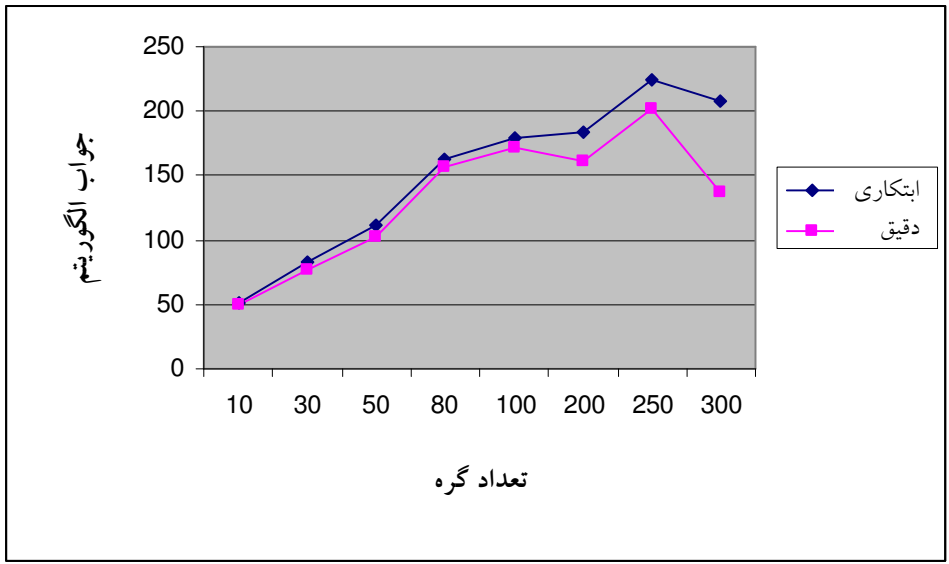

شکل ۵-۵- نمودار مقایسه‌ای پاسخ الگوریتم‌های دقیق و ابتکاری

با توجه به تعداد گره مسائل نمونه

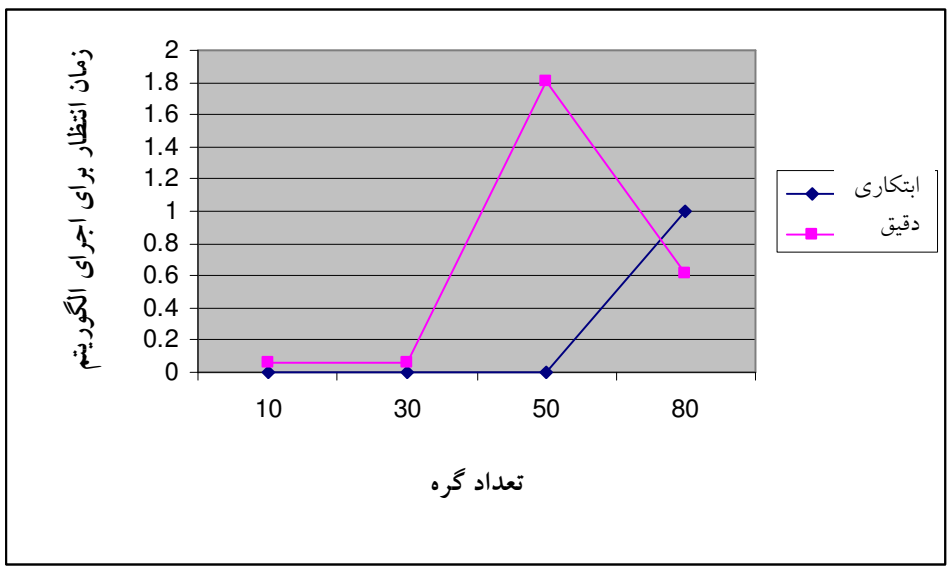

شکل ۵-۶- زمان انتظار برای اجرای الگوریتم‌های دقیق و ابتکاری

با توجه به تعداد گره مسائل نمونه (۱ تا ۸۰ گره)

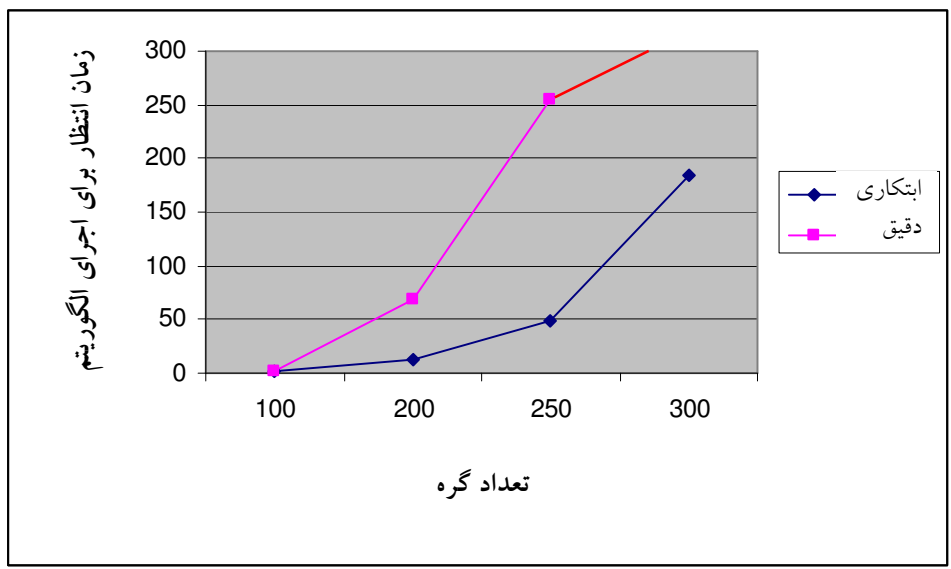

شکل ۵-۷- زمان انتظار برای اجرای الگوریتم‌های دقیق و ابتکاری

با توجه به تعداد گره مسائل نمونه (۱۰۰ تا ۳۰۰ گره)

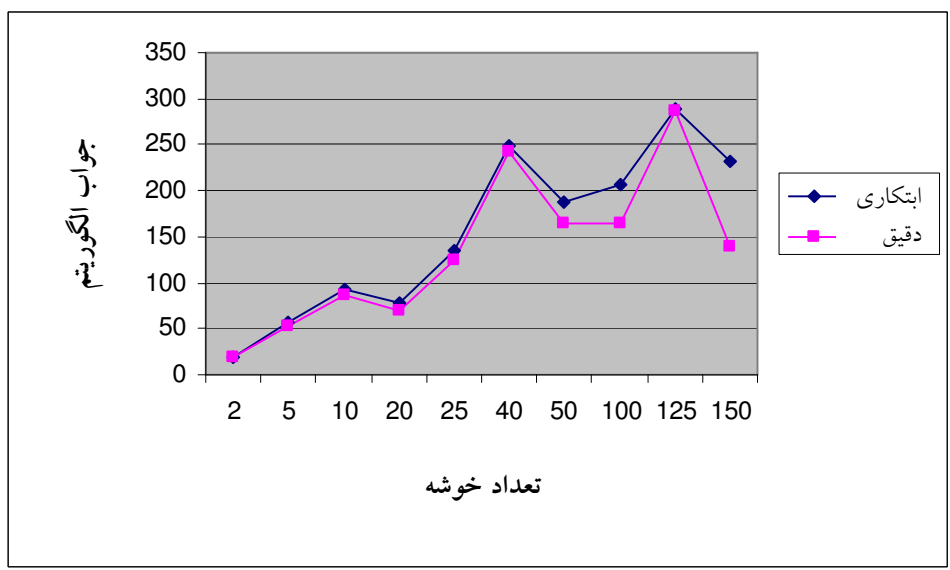

شکل ۵-۸- نمودار مقایسه‌ای پاسخ الگوریتم‌های دقیق و ابتکاری

با توجه به تعداد خوشه مسائل نمونه

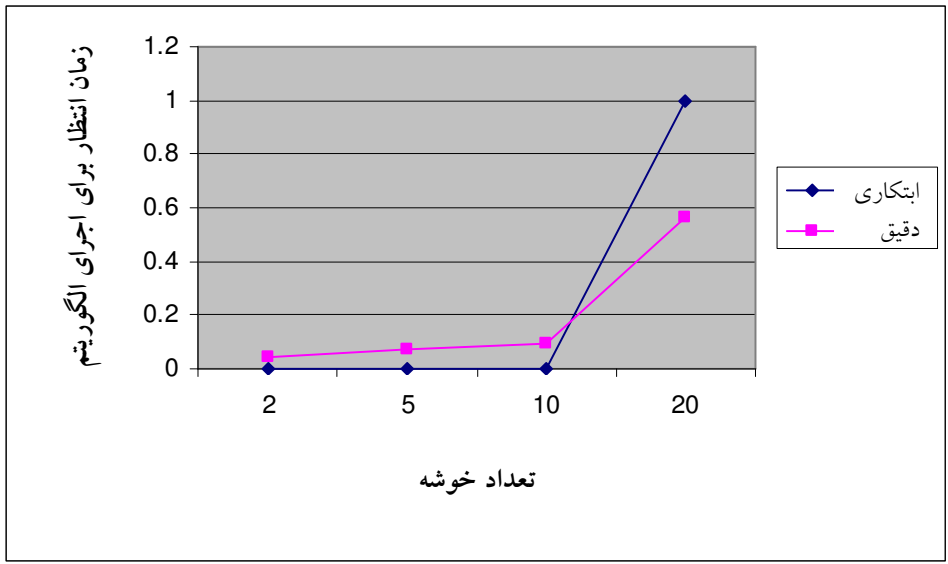

**شکل ۵-۹- زمان انتظار برای اجرای الگوریتم‌های دقیق و ابتکاری با توجه به تعداد خوشه مسائل نمونه (۲ تا ۲۰ خوشه)**

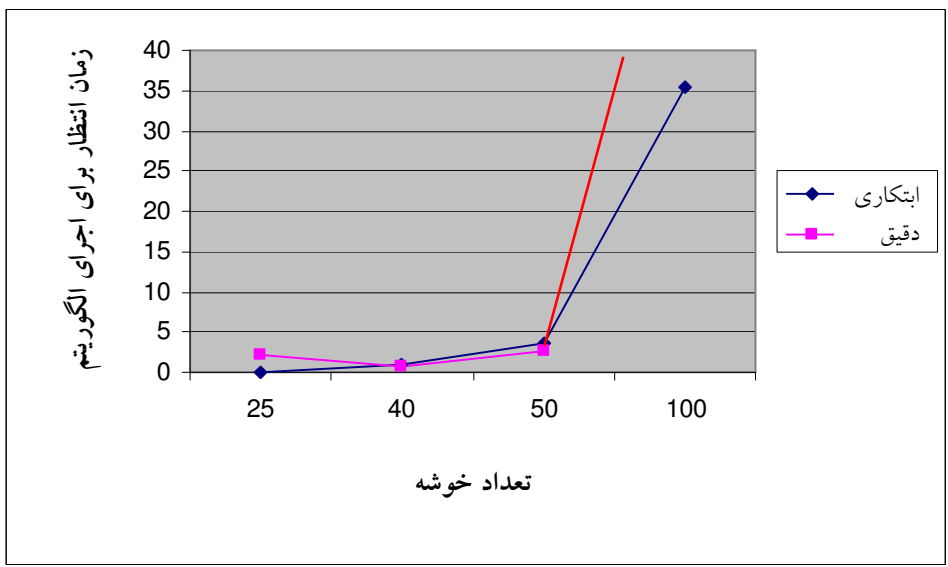

**شکل ۵-۱۰- زمان انتظار برای اجرای الگوریتم‌های دقیق و ابتکاری با توجه به تعداد خوشه مسائل نمونه (۲۵ تا ۱۰۰ خوشه)**

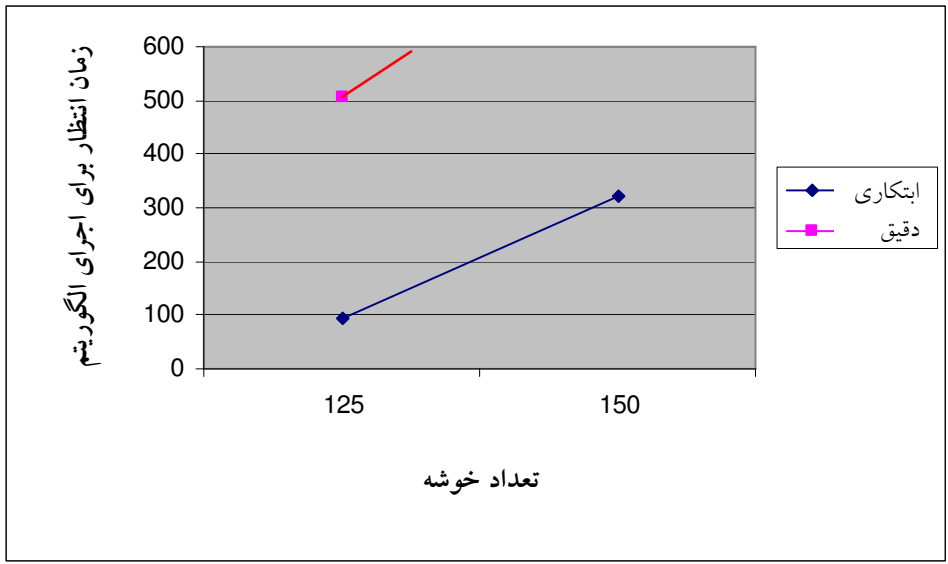

شکل ۵-۱۱- زمان انتظار برای اجرای الگوریتم‌های دقیق و ابتکاری

با توجه به تعداد خوشه مسائل نمونه (۱۲۵ تا ۱۵۰ خوشه)

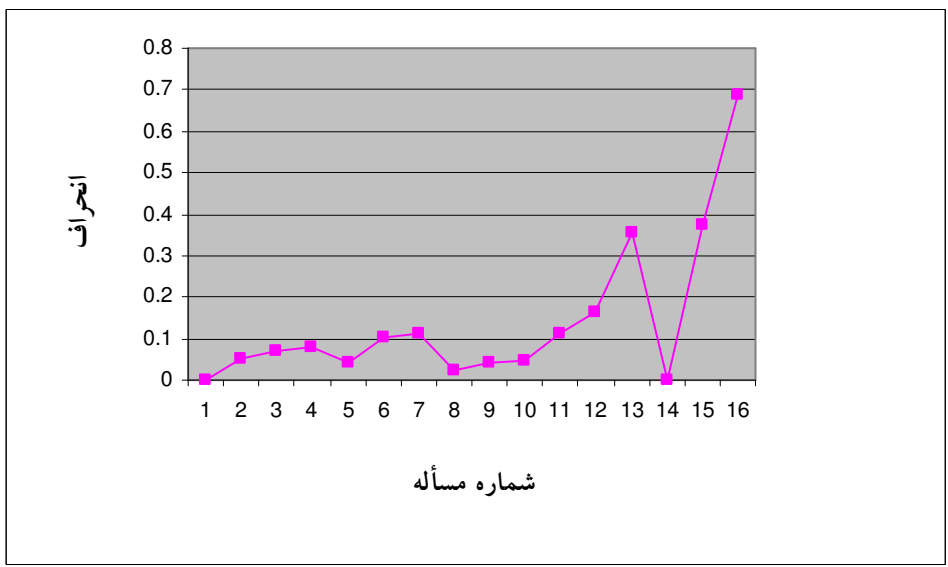

شکل ۵-۱۲- درصد انحراف جواب حاصل از الگوریتم ابتکاری

از بهینه جواب

# فصل ششم

# نتیجه‌گیری و پیشنهاد تحقیقات آتی

## ۶-۱- جمع بندی و نتیجه‌گیری

در این رساله سعی شد مسأله‌ی مسیریابی، به عنوان مسأله‌ای محوری در ادبیات بهینه‌سازی ترکیبیاتی[1]، با پیچیده‌ترین مسأله‌ی ترکیبیاتی (مسأله مسیریابی- مکان‌یابی) [Salhi, [in print] و Nagy ترکیب شده و در عین حال با در نظر گرفتن پارامترهای احتمالی جنبه‌ای واقعی‌تر به آن داده شود. با توجه به اینکه این مسأله تاکنون مدلسازی و حل نشده بود، لذا تمام مراحل مدلسازی آن صورت گرفته و در رساله تشریح شده است. همچنین با استفاده از الگوریتم‌های حل برای مسائل مشابه و به کارگیری روش L-Shape برای بهینه‌سازی احتمالی، الگوریتم حل دقیق برای این مسأله توسعه داده شده است. با توجه به بزرگ بودن ابعاد مسائل واقعی روش حل دقیق برای حل مسائل کاربردی تا ابعادی محدود (۲۵۰ گره و ۱۵۰ خوشه) قابل استفاده است و لذا لزوم طراحی الگوریتم ابتکاری مناسب برای دسترسی به حل مناسب احساس می‌شد. از این رو با استفاده از الگوریتم‌های ابتکاری کلاسیک برای حل مسأله فروشنده دوره گرد، الگوریتم ابتکاری برای حل این مسأله پیشنهاد شد که براساس نتایج محاسباتی امکان حل بزرگترین مسائل (در مقایسه با مسائل مشابه حل شده که حالت خاص این مسأله به حساب می‌آیند) را در کمتر از ۱۰ دقیقه دارد. همچنین دقت این الگوریتم در مقایسه با کران‌های بالای محاسبه شده بسیار خوب ارزیابی

---

1- Combinatorial Optimization

می‌شود. البته به دلیل محدودیت موجود در به کارگیری نرم‌افزار Cplex ۱۰.۱ برای کاربردهای آموزشی، امکان تنظیم برنامه به طور بهینه و در نتیجه افزایش راندمان اجرای آن برای حل مسائل بزرگتر ممکن نشد. با این وجود مقایسه نتایج محاسباتی بیانگر کارکرد مناسب الگوریتم دقیق برای حل مسائل واقعی است و حتی استفاده از تکنیک‌های بهینه‌سازی در ابعاد بزرگ امکان حل مسائل با ابعاد متوسط را در زمانی کوتاه‌تر از الگوریتم ابتکاری نشان می‌دهد.

در مواجهه با مسائل واقعی به طور روز افزون با پیچیدگی‌هایی چون تصادفی بودن و بروز محدودیت‌های مختلف مواجه هستیم. همچنین رقابت میان صنایع تولیدی و خدماتی مختلف ارزش بهینه‌سازی مسائل را دوچندان می‌کند. از این رو سعی در مدل‌سازی واقعی‌تر، امکان دسترسی به پاسخ‌های بهتر را ایجاد کرده و همچنین بهره بردن از تکنیک‌های بهینه‌سازی در ابعاد بزرگ، دسترسی به پاسخ‌های قابل اتکا را تسهیل می‌نماید. روش‌های بهینه‌سازی ابتکاری نیز به عنوان مکملی مؤثر در تسهیل شرایط تولید و خدمت‌رسانی مفید خواهد بود. مسأله مسیریابی نیز به عنوان یکی از پرکاربردترین مسائل حوزه بهینه‌سازی که در حوزه‌های مختلف بطور مستقیم و غیرمستقیم مورد استفاده قرار می‌گیرد از این قاعده مستثنی نیست. مدل‌سازی واقعی‌تر و توسعه الگوریتم‌های مؤثر برای حل مسائل در ابعاد بزرگ تأثیر جدی در بهبود کاربردهای آن خواهد داشت و لذا بایستی جدی گرفته شود.

## ۶-۲- پیشنهاد برای تحقیقات آتی

با توجه به تعریف مسأله مورد نظر و قابلیت تعمیم آن در شبکه‌های توزیع متفاوت، می‌توان زمینه‌های تحقیقاتی آتی را برای آن متصور شد:

### ۱- تغییر در مدل‌سازی مسأله

- اضافه کردن محدودیت‌های دیگر از جمله محدودیت پنجره‌ی زمانی[1]، محدودیت تحویل دادن و تحویل گرفتن[2] و ...

و ساختن مدل‌های متنوع مسیریابی احتمالی توسعه‌یافته، مشابه حالت کلاسیک

- توسعه مدل براساس مدل‌سازی مبتنی بر کمان
- توسعه مدل بر روی گراف جهت‌دار[3]

---

1- Time Windows
2- Pick up and Delivery
3- asymmetric graph

۲- تغییر در مدل‌سازی احتمالی مسأله

● هزینه‌ی حمل و نقل، تعداد وسائل نقلیه در دسترس، ظرفیت وسائل نقلیه و ... را می‌توان به طور احتمالی مدل کرد.

● می‌توان با در نظر گرفتن چند پارامتر تصادفی، مسأله را بطور چند مرحله‌ای مدل‌سازی کرد.

۳- بهبود روش‌های حل

۱- استفاده از تکنیک‌های شبیه‌سازی برای بهبود کارایی الگوریتم در دسترس به جواب بهینه در زمان کمتر از جمله استفاده از تکنیک SAA[1]

۲- به کارگیری الگوریتم‌های کلاسیک دیگر، به کار رفته در حل مسائل مسیریابی برای بهینه‌سازی احتمالی این مدل

---

1- Sample Avarage Approximation

# فهرست مراجع

# پیوست

پیوست الف

# اصل تجزیه

در مسائل عملی، بسیاری از مسائل برنامه‌ریزی خطی از نظر ابعاد آن قدر بزرگ می‌باشند، که نمی‌توان حتّی در کامپیوترهای پیشرفته امروزی حل نمود. مسائلی را می‌توان به صورت برنامه‌ریزی خطّی فرموله نمود که دارای هزاران سطر و ده‌ها هزار ستون باشند. در چنین مسائلی، باید روشهائی به کار گرفت که مسائل بزرگ را به تعدادی از مسائل کوچک تبدیل نماید، به طوری که آن مسائل کوچک با کامپیوترهـای امروزی قابل حل باشند و از طریـق حـل ایـن مسـائل کوچـک بـه حـل مسـأله‌ی بـزرگ اصـلی رسـید. خوشبختانه، تکنیکی موجود است، که به اصل تجزیه معروف می‌باشد و دقیقاً همین کار را انجام می‌دهد.

حتی اگر ابعاد مسأله‌ی برنامه‌ریزی خطی قابل دستکاری باشد، بعضی از قیدهای مسأله ممکـن اسـت دارای ساختمان خاصّی باشند که بتوان مسأله را با راندمان بیشتری دستکاری نمود. در چنین حالت‌هائی می‌توان مسأله را تجزیه نمود به یک مسأله در حالت کلی و یک مسأله با ساختمان خـاص، کـه در ایـن صورت روش مؤثرتری برای حل مسأله می‌توان به کار گرفت. در این جا نیز اصل تجزیه را می‌توان برای حل چنین مسأله‌ی برنامه‌ریزی خطّی با راندمان بیشتری به کار برد.

اصل تجزیه یک روش سیستماتیک برای حل مسائل بزرگ برنامه‌ریزی خطّی، یا مسائل برنامـه‌ریـزی خطّی با ساختمان خاص قیدها می‌باشد. قیدها به دو مجموعه افراز می‌گردند، قیدهای عمومی (یا قیدهای پیچیده) و قیدهائی با ساختمان خاص. در تئوری که ذیلاً ارائه می‌گردد، نشان داده می‌شود کـه ضـرورتی ندارد بعضی از قیدها ساختمان خاصّی داشته باشند تا بتوان از ایـن روش اسـتفاده کـرد؛ اگرچـه داشـتن ساختمان خاصی برای قیود مسأله، راندمان کار را در حل مسأله با به کار بردن این روش افزایش می‌دهد.

استراتژی روش تجزیه آن است که روی دو مسأله‌ی برنامه‌ریزی خطّی جدا عمل می‌نماید، یکی از مسائل روی قیدهای عمومی و دیگری روی قیدهائی با ساختمان خاص می‌باشد. اطلاعات بین ایـن دو مسـأله رد و بدل می‌گردد، تا به نقطه‌ای می‌رسیم که جواب بهینه مسأله‌ی اصلی می‌باشد. مسأله برنامـه‌ریـزی خطّـی روی قیدهای عمومی به مسأله مادر و مسأله‌ی برنامه‌ریزی خطّی روی قیودی با ساختمان خاص را، مسأله‌ی فرعـی می‌نامیم. مسأله‌ی مادر، مجموعه‌ی جدیدی از ضرایب قیمت‌ها را به مسأله‌ی فرعی رد مـی‌کنـد و در عـوض ستون جدیدی که مبتنی بر ضرایب قیمت‌ها می‌باشد دریافت می‌نماید. در شروع مطلب فرض بر آن است که مجموعه‌ی تعریف شده به وسیله‌ی قیدهائی با ساختمان خاص، محدود می‌باشد. پس از اینکـه اصـل تجزیـه مورد مطالعه قرار گرفت این فرض رها خواهد شد و حالت کلّی مورد بحث قرار می‌گیرد و مطالعه به حـالتی که چندین مسأله‌ی فرعی موجود خواهد بود بسط و گسترش داده می‌شود.

## پ-۱- الگوریتم تجزیه

مسأله‌ی برنامه‌ریزی خطّی زیر را در نظر می‌گیریم.

Min z = cx
S.T.
  Ax = b
  Dx = d
  x ≥ ○

اگر فرض کنیم $x = \{x | Dx = d \ \& \ x \geq ○\}$

مسأله به صورت زیر درمی‌آید.

Min cx
S.T.
  Ax = b
  x ∈ x

برای سهولت فرض بر آن است که x مجموعه‌ی محدود می‌باشد. چون x یک مجموعه‌ی محدّب محدود می‌باشد (چند وجهی محدّب محدود)، آنگاه هر x ∈ x را می‌توان به صورت ترکیب محدّب از رئوس x نوشت. اگر رئوس x را به $x_1, x_2, ..., x_t$ نشان دهیم، برای هر x ∈ x داریم:

$$x = \sum_{j=1}^{t} \lambda_j x_j$$

$$\sum_{j=1}^{t} \lambda_j = 1$$

$$\lambda_j \geq ○ \quad j = 1, 2, ..., t$$

با قرار دادن x در مسأله‌ی بهینه‌سازی قبلی، مسأله‌ی فوق الذکر به صورت زیر درمی‌آید؛ که در این مسأله، متغیرها، $\lambda_1, ..., \lambda_t$ می‌باشند.

$$\text{Min} \sum_{j=1}^{t} (cx_j) \lambda_j$$

S.T.

$$\sum_{j=1}^{t} (Ax_j) \lambda_j = b \tag{پ-۱}$$

$$\sum_{j=1}^{t} \lambda_j = 1 \tag{پ-۲}$$

$$\lambda_j \geq ○ \quad j = 1, 2, ..., t \tag{پ-۳}$$

چون t ؛ تعداد نقاط رأسی مجموعه‌ی x، معمولاً بزرگ می‌باشد، کوشش برای پیدا نمودن مختصات تماس رئوس x یعنی مختصات نقاط $x_1, x_2, ..., x_t$ و حل مسأله بدین طریق، کار بس مشکلی است. به جای این کار سعی می‌شود جواب بهینه‌ای برای این مسأله (در نتیجه جواب بهینه‌ای برای مسأله‌ی اصلی) بدون محاسبه‌ی کلیه مختصات رئوس x پیدا نمود.

### کاربرد روش سیمپلکس اصلاح شده

حل مسأله‌ی فوق الذکر را به وسیله الگوریتم سیمپلکس اصلاح شده در نظر می‌گیریم فرض می‌کنیم یک جواب اساسی شدنی به صورت $\lambda = (\lambda_B, \lambda_N)$ باشد. بعلاوه فرض می‌کنیم ماتریس $B^{-1}$ با ابعاد $(m+1)\times(m+1)$ در دست می‌باشد. (بحث اینکه به چه صورتی جواب اساسی شدنی شروع مسأله را به دست می‌آوریم، در قسمت ۳ بیان خواهد شد). متغیرهای دوآل نظیر معادلات (پ-۱) و (پ-۲) را به w و $\alpha$ نشان می‌دهیم. با توجه به آن چه گذشت داریم.

$$(w, \alpha) = \hat{c}_B B^{-1}$$

که در آن $\hat{c}_B$ ضریب متغیرهای اساسی در تابع مقصود است و $\hat{c}_j = cx_j$ برای هر یک از متغیرهای اساسی $\lambda_j$. عکس ماتریس پایه، متغیرهای دوآل، مقادیر متغیرهای اساسی و مقدار تابع مقصود در جدول زیر نشان داده شده است که در آن $\overline{b} = B^{-1}\begin{pmatrix} b \\ 1 \end{pmatrix}$.

| پایه | RHS |
|---|---|
| $(w, \alpha)$ | $\hat{c}_B \overline{b}$ |
| $B^{-1}$ | $\overline{b}$ |

با استفاده از اطلاعات این جدول می‌توان گفت که آیا جواب حاصل، بهینه است یا نه. اگر بهینه نیست کدام یک از متغیرهای غیراساسی را باید افزایش دهیم. این کار با محاسبه‌ی زیر صورت می‌گیرد.

$$z_k - \hat{c}_k = \underset{1 \leq j \leq t}{\text{Max}} (z_j - \hat{c}_j) = \underset{1 \leq j \leq t}{\text{Max}} (w, \alpha) \begin{bmatrix} Ax_j \\ 1 \end{bmatrix}$$

$$9 - cx_j = \underset{1 \leq j \leq t}{\text{Max}} \, wAx_j + \alpha - cx_j$$

$$= \underset{1 \leq j \leq t}{\text{Max}} \left[(wA - c) x_j + \alpha\right] \tag{پ-۴}$$

چون $z_j - \hat{c}_j = 0$ برای متغیرهای غیراساسی می‌باشد، پس مقدار ماکزیمم قبلی بزرگتر یا مساوی صفر می‌باشد. بنابراین، اگر $z_k - \hat{c}_k = 0$، آنگاه $z_j - \hat{c}_j \leq 0$ برای تمام متغیرهای غیراساسی و جواب بهینه به دست آمده است. از طرف دیگر اگر $z_k - \hat{c}_k > 0$، آنگاه متغیر غیراساسی $\lambda_k$ افزایش داده می‌شود.

مشخص نمودن ایندکس k، ایندکس متغیر وارد شونده، با به کار بردن معادله (پ-۴) در این مرحله از کار شدنی نیست، زیرا t عدد بسیار بزرگی است و نقاط رأسی $x_j$ها متناظر متغیرهـای اساسـی $\lambda_j$ قبـلاً محاسبه نشده است. چون X یک چندوجهی محدود می‌باشد، بنابراین ماکزیمم مقدار هر تابع خطـی بـه عنوان تابع مقصود بر روی X در یکی از نقاط رأسی آن اتفاق می‌افتد. بنابراین:

$$\underset{1 \leq j \leq t}{\text{Max}} (wA - c)x_j + \alpha = \underset{x \in X}{\text{Max}} (wA - c)x + \alpha$$

برای یک جواب اساسی داده شده $(\lambda_B, \lambda_N)$ با متغیرهای دوآل $(w, \alpha)$، مسأله فرعـی برنامـه‌ریـزی خطّی زیر را که مسأله ساده‌ای می‌باشد حل می‌نماییم (ساده از این نظر که X دارای ساختمان خاص است).

$$\text{Max } (wA - c)x + \alpha$$
S.T.
$$x \in X$$

لازم به ذکر است که تابع مقصود، مقدار ثابتی را همراه دارد. پس در شروع به جای اینکه مقـدار تـابع مقصود را که معمولاً صفر قرار می‌دادیم، $\alpha$ قرار خواهیم داد. فرض می‌کنیم $x_k$ جواب بهینه مسأله قبلی با مقدار تابع مقصود $z_k - \hat{c}_k$ باشد. اگر $z_k - \hat{c}_k = 0$، آنگاه جواب اساسی شدنی $(\lambda_B, \lambda_N)$ جواب بهینه می‌باشند. در غیر این صورت اگر $z_k - \hat{c}_k > 0$، آنگاه متغیر $\lambda_k$ وارد پایه می‌گـردد. سـتون محـوری بـا ضرب $B^{-1}$ در $\begin{pmatrix} Ax_k \\ 1 \end{pmatrix}$ به صورت زیر محاسبه می‌گردد.

$$\overline{A}_{.k} = B^{-1} A_{.k} = B^{-1} \begin{pmatrix} Ax_k \\ 1 \end{pmatrix}$$

چون X محدود می‌باشد پس هیچگاه $\overline{A}_{.k} \leq 0$ نخواهد بود. محـدود بـودن X، باعـث بـه وجـود آمـدن مسأله‌ی محدود خواهد بود. ستون به روز شده $\begin{pmatrix} z_k - \hat{c}_k \\ \overline{A}_{.k} \end{pmatrix}$ به جدول فوق الذکر الحـاق مـی‌شـود. متغیـر اساسی $\lambda_{B_r}$ به وسیله قاعده‌ی مینیموم مشخص می‌گردد تا پایه را ترک نماید. عکس ماتریس پایه جدیـد، مقادیر متغیرهای دوآل و بردار سمت راست به روز می‌گردد (این کار به وسیله‌ی عمل محـوری صـورت می‌گیرد که $\overline{a}_{rk}$ عضو محوری می‌باشد). پس از انجام عمل محوری، دوباره تکرار می‌گردد.

مرحله مادر یک جواب شدنی بهتری که باعث کاهش تابع مقصود می‌باشد به ما می‌دهد و مرحله فرعی (مسأله فرعی) بررسی می‌نماید که آیا $z_j - \hat{c}_j \leq 0$ برای هر $\lambda_j$ می‌باشد یا نه و در صورت منفی بودن جواب مثبت‌ترین $z_k - \hat{c}_k$ را مشخص می‌نماید.

## خلاصه الگوریتم تجزیه

### قدم اولیه

یک جواب اساسی شدنی برای معادلات (پ-۱)، (پ-۲) و (پ-۳) پیدا می‌کنیم. فرض می‌کنیم $B$ ماتریس پایه و $B^{-1}$ عکس آن باشد. جدول زیر را که در آن $(w, \alpha) = \hat{c}_B B^{-1}$ و $\bar{b} = \bar{B}^{-1}\begin{pmatrix} b \\ 1 \end{pmatrix}$ می‌باشد تشکیل می‌دهیم (یادآور می‌شویم که $\hat{c}_j = cx_j$)

| RHS | عکس پایه |
|---|---|
| $\hat{c}_B \bar{b}$ | $(w, \alpha)$ |
| $\bar{b}$ | $B^{-1}$ |

### قدم اصلی

۱- مسأله‌ی فرعی زیر را حل می‌نماییم.

$$\text{Max } (wA - c)x + \alpha$$
$$\text{S.T.}$$
$$x \in X$$

فرض می‌کنیم $x_k$ جواب اساسی شدنی بهینه برای مسأله فوق با مقدار تابع مقصود $z_k - \hat{c}_k$ باشد. اگر $z_k - \hat{c}_k = 0$ توقف می‌کنیم. جواب اساسی شدنی مرحله‌ی مادر جواب بهینه مسأله می‌باشد. در غیر این صورت به مرحله‌ی ۲ می‌رویم.

۲- فرض می‌کنیم $\overline{A}_{.k} = B^{-1} \begin{pmatrix} Ax_k \\ 1 \end{pmatrix}$ و ستون به روز شده الحاق $\begin{pmatrix} z_k - \hat{c}_k \\ \overline{A}_{.k} \end{pmatrix}$ برای مسأله مادر باشـد. $\overline{a}_{rk}$ را عضو محوری قرار می‌دهیم. که در آن ایندکس $r$ از رابطه زیر مشخص شده است.

$$\frac{\overline{b}_r}{\overline{a}_{rk}} = \min_{1 \leq i \leq m+1} \left\{ \frac{\overline{b}_i}{\overline{a}_{ik}} : \overline{a}_{ik} > 0 \right\}$$

با انجام عمل محوری جدول را به روز می‌نماییم. این عمل عکس مـاتریس پایـه، متغیرهـای دوآل و بردار سمت راست را به دست می‌دهد.

پس از عمل محوری ستون $\lambda_k$ را حذف نموده و مرحله‌ی ۱ را تکرار می‌نماییم.

### تذکر بعضی نکات در رابطه با الگوریتم تجزیه

۱- الگوریتم فوق الذکر کاربرد مستقیم سیمپلکس اصلاح شده می‌باشد، جز محاسبه‌ی $z_k - \hat{c}_k$ کـه بـا حل مسأله‌ی فرعی صورت می‌گیرد. بنابراین الگوریتم در تعداد متناهی از مراحل متقـارب مـی‌باشـد، بـه شرط آن که قاعده‌ی جلوگیری از به دور افتادن در مرحله‌ی مـادر و مـسأله‌ی فرعـی در صـورت وجـود تبهگنی به کار برده شود.

۲- در هر تکرار، مرحله‌ی مادر یک جواب اساسی شدنی جدید بهتر برای معادلات (پ-۱) و (پ-۲) و (پ-۳) به وسیله‌ی معرفی متغیر غیراساسی $\lambda_k$ به پایه که توسط مسأله فرعی تولید شده؛ به ما می‌دهد. در هر تکرار مسأله‌ی فرعی یک نقطه‌ی رأسی $x_k$ تولید می‌کند؛ که متناظر ستون بـروز شـده $\begin{pmatrix} z_k - \hat{c}_k \\ \overline{A}_{.k} \end{pmatrix}$ می‌باشد، از این رو، بعضی مواقع این روش را روش تولید ستون نیز می‌نامند.

۳- در هر تکرار بردار دوآل متفاوتی از مرحله‌ی مادر به مسأله‌ی فرعـی رد مـی‌شـود. بـه جـای حـل مسأله‌ی فرعی از نو، پایه بهینه‌ی آخرین تکرار را می‌توان با تعدیل سطر قیمت، به عنوان پایـه شـروع بـه کار گرفت.

۴- اگر قیدهای مادر از نوع نامساوی باشند، آنگاه مقادیر $z_j - \hat{c}_j$ برای متغیرهـای زایـدی کـه جهـت برگرداندن این قیود به قیود تساوی به کـار گرفتـه شـده؛ بایـستی محاسـبه گردیـده و جهـت وارد شـدن احتمالی به پایه مورد توجه قرار گیرد. برای قید مادر $I$ از نوع $\leq$ که در آن $S_i$ متغیر زایدی می‌باشـد کـه جهت برگرداندن این قید به قید تساوی به کار رفته داریم:

$$z_{s_i} - c_{s_i} = (w, \alpha) \begin{pmatrix} e_i \\ 0 \end{pmatrix} - 0 = w_i$$

۸۸

بنابراین برای مسأله مینیمم نمودن، متغیر زاید متناظر قید از نوع ≤، زمانی می‌تواند وارد پایه گردد کـه $w_i > 0$ باشد (اگر قید از نوع ≥ باشد قید زاید متناظر می‌تواند وارد پایه شود که $w_i < 0$).

۵- در هر تکرار، مسأله‌ی فرعی ضرورتی ندارد که کاملاً بهینه گردد. کافی است $z_k - \hat{c}_k$ هـایی یافت شود که مثبت باشد. در این حالت $\lambda_k$ متغیر غیراساسی می‌باشد که می‌تواند وارد پایه شـود (در مـسأله‌ی مادر).

**محاسبه و کاربرد بندهای پائین**

الگوریتم تجزیه متوقف می‌گردد اگر $\text{Max } z_j - \hat{c}_j = 0$. با توجه به اینکه تعداد متغیرهـای $\lambda_1, ..., \lambda_t$ زیاد می‌باشد، ادامه محاسبه تا زمانی که شرط فوق برقرار گردد، ممکـن اسـت وقت‌گیـر و بـرای مـسائل بزرگ بسیار گران باشد. از این رو باند پائین برای تابع مقصود را برای هر جواب شدنی به دسـت آورد که در تصمیم‌گیری نهائی کمک فراوانی می‌نماید.

چون الگوریتم تجزیه نقاط شدنی که باعث بهتر شدن مقدار تابع مقصود مـی‌گـردد، ممکـن اسـت در مرحله‌ی متوقف کند که اختلاف باند پـائین در آن مرحلـه و مقـدار تـابع مقصـود در آن مرحلـه از عـدد مشخّصی کمتر گردد (در بازه‌ی قابل قبولی باشد). در این صورت جواب حاصل، جـواب بهینـه نخواهـد بود، امّا تضمین آن را دارد که در حد قابل قبول جواب شدنی خوبی باشد و با دقّت معینـی قابـل قبـول است. مسأله فرعی زیر را در نظر می‌گیریم؛

$\text{Max } (wA - c)x + \alpha$
S.T.
$\quad x \in X$

که در آن $w$ بردار دوآلی است که از مرحله‌ی مادر دریافت شده است. فرض می‌کنیم مقدار بهینه تـابع مقصود مسأله‌ی فوق $z_k - \hat{c}_k$ باشد. حال فرض می‌کنیم $x$ یک جواب شدنی مسأله‌ی اصلی باشد، یعنی $x$ در قیدهای $Ax = b$ صدق می‌کند و $x \in X$. طبق تعریف $z_k - \hat{c}_k$ و چون $x \in X$ داریم

$$(wA - c)x + \alpha \leq (z_k - \hat{c}_k)$$

(زیرا $z_k - \hat{c}_k$ ماکزیمم مقدار $(wA - c)x + \alpha$ به ازاء هر $x \in X$ می‌باشد) چون $Ax = b$ است و در نتیجه نامساوی فوق را می‌توان بصورت زیر نوشت:

$$cx \geq wAx - (z_k - \hat{c}_k) + \alpha = wb + \alpha - (z_k - \hat{c}_k) = \hat{c}_B \bar{b} - (z_k - \hat{c}_k)$$

۸۹

چون رابطه‌ی فوق برای هر $x \in X$ که در رابطه $Ax = b$ صدق می‌کند برقرار است پس داریم

$$\underset{\substack{ax=b \\ x \in x}}{Max} cx \geq \hat{c}_B \bar{b} - (z_k - \hat{c}_k)$$

به عبارت دیگر $\hat{c}_B \bar{b} - (z_k - \hat{c}_k)$ یک باند پائین برای مجموعه مقادیر تابع مقصود مسأله‌ی اصلی می‌باشد.

## پ-۲- روش به دست آوردن جواب شروع

در این قسمت روش به دست آوردن یک جواب اساسی شدنی برای شروع مسأله را برای مسأله مادر، با به کار بردن متغیرهای مصنوعی را توضیح می‌دهیم. حذف متغیرهای مصنوعی با به کار بردن فاز اول روش دو فازی یا با به کارگیری روش M- بزرگ صورت خواهد گرفت. اگر در نهایت، متغیر مصنوعی با مقدار مثبت موجود باشد، آنگاه مسأله‌ی اصلی جواب شدنی نخواهد داشت.

### قیدهای نامساوی

مسأله‌ی زیر را در نظر می‌گیریم:

$$Min \sum_{j=1}^{t}(cx_j)\lambda_j$$

S.T.

$$\sum_{j=1}^{t}(Ax_j)\lambda_j \leq b$$

$$\sum_{j=1}^{t}\lambda_j = 1$$

$$\lambda_j \geq 0 \quad j = 1, 2, ..., t$$

اگر یک $x_1 \in X$ طوری یافت شود که $Ax_1 \leq b$، آنگاه، پایه زیر در دست می‌باشد، که در آن ماتریس واحد متناظر بردار متغیرهای زاید $s \geq 0$ می‌باشد.

$$B = \begin{bmatrix} I & Ax_1 \\ 0 & 1 \end{bmatrix} \quad \text{و} \quad B^{-1} = \begin{bmatrix} I & -Ax_1 \\ 0 & 1 \end{bmatrix}$$

جدول شروع چنین است.

| | عکس پایه | | RHS |
|---|---|---|---|
| Z | ۰ | $cx_1$ | $cx_1$ |
| S | I | $-Ax_1$ | $b - Ax_1$ |
| $\lambda_1$ | ۰ | ۱ | ۱ |

حال فرض می‌کنیم که هیچ $x \in X$ که در رابطه $Ax \leq b$ صدق کند وجود نداشته باشد. در این حالت، پس از برگرداندن مسأله‌ی مادر به صورت تساوی با افزودن متغیرهای زاید، قیدها را به صورتی دست کاری می‌نمائیم که بردار سمت راست نامنفی باشد یعنی $b \geq 0$. آنگاه، با افزودن متغیرهای تصنعی یک ماتریس واحد در قیدها ایجاد می‌نمائیم. این ماتریس واحد، پایه شروع را تشکیل می‌دهد. (با به کار بردن روش M- بزرگ یا دو فازی می‌توان متغیرهای مصنوعی را حذف نمود.)

**قیدهائی به صورت تساوی**

در این حالت $m+1$ متغیر مصنوعی را به قیدها می‌افزاییم که پایه اولیه شدنی را تشکیل دهند. در این جا نیز با به کار بردن روش M- بزرگ یا دو فازی می‌توان متغیرهای مصنوعی را حذف کرد.

**پ-۴- حالتی که ناحیه‌ی X نامحدود می‌باشد**

برای ناحیه‌ی نامحدود X، الگوریتم تجزیه بایستی جزئی تعدیل شود. در این حالت تمام نقاط X را نمی‌توان به صورت ترکیب محدّب از نقاط رأسی X نوشت، در صورتی که هر $x \in X$ را می‌توان به صورت ترکیب محدّب نقاط رأسی بعلاوه ترکیب خطّی نامنفی نقاط رأسی همگن X نوشت یعنی اگر $x_1, ..., x_t$ نقاط رأسی X و $d_1, ..., d_l$ نقاط رأسی همگن X باشند آنگاه هر $x \in X$ را می‌توان به صورت زیر نوشت

$$x = \sum_{j=1}^{t} \lambda_j x_j + \sum_{j=1}^{l} \mu_j d_j$$

$$\sum_{j=1}^{t} \lambda_j = 1$$

$$\lambda_j \geq 0 \qquad j = 1, 2, ..., t$$

$$\mu_j \geq 0 \qquad j = 1, 2, ..., l$$

۹۱

مسأله را می‌توان به مسأله‌ای که در آن $\mu_1, \mu_2, ..., \mu_l$ و $\lambda_1, \lambda_2, ..., \lambda_t$ متغیرهای مسأله می‌باشند تبدیل نمود، به صورت زیر:

$$\text{Min} \sum_{j=1}^{t}(cx_j)\lambda_j + \sum_{j=1}^{l}(cd_j)\mu_j$$

S.T.

$$\sum_{j=1}^{t}(Ax_j)\lambda_j + \sum_{j=1}^{l}(Ad_j)\mu_j = b \qquad (پ-5)$$

$$\sum_{j=1}^{t}\lambda_j = 1 \qquad (پ-6)$$

$$\lambda_j \geq 0 \qquad j=1,2,...,t$$
$$\mu_j \geq 0 \qquad j=1,2,...,l$$

چون t و l معمولاً بزرگ هستند، سعی خواهد شد که این مسأله با روش سیمپلکس اصلاح شده حل شود. فرض می‌کنیم که یک پایه شدنی برای مسأله‌ی بالائی جهت شروع در دست باشد و فرض می‌کنیم w و $\alpha$ متغیرهای دوآل مسأله باشند که متناظر قیدهای (پ-5) و (پ-6) هستند. بعلاوه فرض می‌کنیم $(w, \alpha) = \hat{c}_B B^{-1}$، $B^{-1}$ ($\hat{c}_B$ ضریب متغیرهای اساسی در مسأله می‌باشد) و $\bar{b} = B^{-1}\begin{pmatrix} b \\ 1 \end{pmatrix}$؛ معلوم و در جدول زیر نشان داده شده‌اند.

| RHS | عکس پایه |
|---|---|
| $\hat{c}_B \bar{b}$ | $(w, \alpha)$ |
| $\bar{b}$ | $B^{-1}$ |

یادآور می‌شویم که در جدول فوق الذکر، جواب حاصل از آن بهینه می‌باشد اگر، برای تمام متغیرها داشته باشیم $z_j - \hat{c}_j \leq 0$ یعنی برای جواب بهینه باید داشته باشیم.

برای $\lambda_j$های غیراساسی

$$0 \geq z_j - \hat{c}_j = (w, \alpha)\begin{bmatrix} Ax_j \\ 1 \end{bmatrix} - cx_j = wAx_j + \alpha - cx_j \qquad (پ-7)$$

و برای $\mu_j$های غیراساسی

$$0 \geq z_j - \hat{c}_j = (w, \alpha)\begin{bmatrix} Ad_j \\ 0 \end{bmatrix} - cd_j = wAd_j - cd_j \qquad (پ-8)$$

چون تعداد متغیرهای غیراساسی معمولاً بسیار بزرگ می‌باشند، چک نمودن شرایط (پ-۷) و (پ-۸)، به وسیله‌ی تولید نمودن نقاط رئوس متناظر آنها و جوابهای همگن رأسی مورد نظر، از نظر محاسباتی شدنی نمی‌باشد. به هر حال، می‌توان به وسیله‌ی حل مسأله‌ی فرعی زیر نشان داد که این شرایط برقرار است یا نه؟ مهمتر از همه اینکه وقتی مسأله‌ی فرعی حل می‌شود، اگر شرایط (پ-۷) و (پ-۸) برقرار نباشد، و $z_k - \hat{c}_k > 0$ آنگاه متغیر غیراساسی متناظر وارد پایه می‌گردد.

$$\text{Max } (wA - c) x + \alpha$$
$$\text{S.T.}$$
$$x \in x$$

ابتدا فرض می‌کنیم، که جواب بهینه مسأله فرعی فوق الذکر نامحدود می‌باشد. یادآور می‌شویم که زمانی این مطلب صادق است که برای یک جواب همگن رأسی مانند $d_k$، داشته باشیم $z_k - \hat{c}_k > 0$. در این حالت $\begin{pmatrix} Ad_k \\ 0 \end{pmatrix}$ محاسبه و بروز می‌گردد و با ضرب نمودن $B^{-1}$ از سمت چپ این ستون، ستون $\begin{pmatrix} z_k - \hat{c}_k \\ \bar{A}_{.k} \end{pmatrix}$ محاسبه می‌گردد و به جدول فوق الذکر الحاق می‌گردد، و روش سیمپلکس اصلاح شده ادامه پیدا می‌کند.

حال فرض می‌کنیم که جواب بهینه محدود باشد (جواب بهینه‌ی مسأله فرعی) شرط لازم و کافی برای اینکه مسأله محدود باشد، آن است که، برای هر نقطه‌ی همگن رأسی داشته باشیم $(wA - c) d_j \leq 0$ و بررسی می‌کنیم که آیا شرایط (پ-۷) برقرار می‌باشد یا نه؟ فرض می‌کنیم $x_k$ جواب بهینه مسأله باشد و مقدار مقصود $z_j - \hat{c}_j$ را در نظر می‌گیریم، اگر $z_j - \hat{c}_j \leq 0$، با توجه به شرط بهینگی $x_k$ برای هر نقطه رأسی $x_j$ داریم؛

$$(wA - c) x_j + \alpha \leq (wA - c) x_k + \alpha = z_k - \hat{c}_k \leq 0$$

از این رو شرایط (پ-۷) برقرار می‌باشد و متوقف می‌شویم و یک جواب برای مسأله‌ی اصلی به دست آمده است. اگر $z_k - \hat{c}_k > 0$ آنگاه $\lambda_k$ وارد پایه می‌گردد. این کار بدین طریق صورت می‌گیرد که ستون $\bar{A}_{.k} = B^{-1} \begin{pmatrix} Ax_k \\ 1 \end{pmatrix}$ به جدول الحاق می‌شود و عمل محوری صورت می‌گیرد. اگر مسأله‌ی مادر (مانند حالتی که مسأله محدود بوده) شامل قیدهائی از نوع نامساوی باشد، آنگاه $z_j - \hat{c}_j$ بایستی برای متغیرهای زاید متناظر آنها محاسبه گردد.

در نتیجه حل مسأله قبل یا به یک جواب بهینه می‌رسد، یا یک متغیر اساسی پیدا کرده و جهت وارد شدن به پایه معرفی می‌شود.

## پ-۳- ساختمان قطری بلوکی

در این قسمت حالت مهمی که در آن x دارای ساختمان خاص قطری بلوک می‌باشد مورد بحث قرار می‌گیرد. در این حالت خود x را می‌توان به مجموعه‌های $x_1, x_2, ..., x_T$ تجزیه نمود که هر کدام از $x_i$ ها، شامل مجموعه‌ای از متغیرها می‌باشند، که در دیگری ظاهر نمی‌گردند. اگر بردار x برطبق همین مجموعه‌ها تجزیه شود، بردارهائی نظیر $x_1, x_2, ..., x_T$ و بردار c به بردارهائی نظیر بردارهای $c_1, c_2, ..., c_T$ و ماتریس A در قید $Ax = b$ به ماتریس‌هائی مانند $A_1, ..., A_T$ تجزیه می‌شود و مسأله به صورت زیر درمی‌آید:

$$\text{Min } c_1x_1 + c_2x_2 + ... + c_Tx_T$$
$$\text{S.T.}$$
$$A_1x_1 + A_2x_2 + ... + A_Tx_T = b$$
$$B_1x_1 \leq b_1$$
$$B_2x_2 \leq b_2$$
$$\ddots \quad \vdots$$
$$B_Tx_T \leq b_T$$
$$x_1, x_2, ..., x_T \geq 0$$

که در آن

$$x_i = \{x_i : B_ix_i \leq b_i , x_i \geq 0\}$$
$$i = 1, 2, ..., T$$

از ساختمان قطری بلوکی x می‌توان جهت بالا بردن راندمان محاسبات بهره‌برداری فراوانی نمود.
برای مسأله‌ی i ام $x_i \in x_i$ اگر و فقط اگر

$$x_i = \sum_{j=1}^{t_i} \lambda_{ij} x_{ij} + \sum_{j=1}^{l_i} \mu_{ij} d_{ij}$$
$$\sum_{j=1}^{t_i} \lambda_{ij} = 1$$
$$\lambda_{ij} \geq 0 \quad j = 1, 2, ..., t_i$$
$$\mu_{ij} \geq 0 \quad j = 1, 2, ..., l_i$$

که در آن $x_{ij}$ نقاط رأسی $x_i$ و $d_{ij}$ جوابهای همگن رأسی مربوطه می‌باشند (اگر $x_i$ جوابهای همگن رأسی داشته باشد).

با قرار دادن $x_i$ در مسأله‌ی فوق الذکر، مسأله‌ی اصلی به صورت زیر درمی‌آید.

$$\text{Min} \sum_{i=1}^{T} \sum_{j=1}^{t_i} (c_i x_{ij}) \lambda_{ij} + \sum_{i=1}^{T} \sum_{j=1}^{l_i} (c_i d_{ij}) \mu_{ij}$$

S.T.

$$\sum_{i=1}^{T} \sum_{j=1}^{t_i} (A_i x_{ij}) \lambda_{ij} + \sum_{i=1}^{T} \sum_{j=1}^{l_i} (A_i d_{ij}) \mu_{ij} = b \qquad (پ-9)$$

$$\sum_{j=1}^{t_i} \lambda_{ij} = 1 \qquad i = 1, 2, ..., T \qquad (پ-10)$$

$$\lambda_{ij} \geq 0 \qquad j = 1, 2, ..., t_i \qquad i = 1, 2, ..., T$$

$$\mu_{ij} \geq 0 \qquad j = 1, 2, ..., l_i \qquad i = 1, 2, ..., T$$

در این جا T ترکیب محدّب می‌باشد یعنی برای هر مسأله‌ی فرعی یک ترکیب محدّب وجـود دارد. این مطلب انعطاف‌پذیری مسأله را بالا می‌برد ولی در عین حال به تعداد قیدهای مسأله مـی‌افزایـد. یعنـی تعداد قیدها از $m+1$ به $m+T$ افزایش پیدا می‌کند.

فرض می‌کنیم یک پایه شدنی برای مسأله فوق الذکر با ابعاد $(m+T) \times (m+T)$ در دست باشد. هـر یک از بردارهای متغیرهای اساسی بایستی حداقل یکی از $\lambda_{ij}$ را به عنوان متغیر اساسی داشته باشد (بـرای هر بلوک حداقل یک $\lambda_{ij}$). بعلاوه فرض می‌کنیم

$$\bar{b} = B^{-1} \binom{b}{1}, B^{-1}$$

$$(w, \alpha) = (w_1, ..., w_m, \alpha_1, ..., \alpha_T) = \hat{c}_B B^{-1}$$

که در آن $\hat{c}_B$ بردار ضریب متغیرهای اساسی می‌باشد.

(برای $\lambda_{ij} ... : \hat{c}_{ij} = c_i x_{ij}$ و برای $\mu_{ij} \quad \hat{c}_{ij} = c_i d_{ij}$)

اینها در جدول زیر نشان داده شده است.

| RHS | عکس پایه |
|---|---|
| $\hat{c}_B \bar{b}$ | $(w, \alpha)$ |
| $\bar{b}$ | $B^{-1}$ |

این جواب بهینه است اگر $z_{ij} - \hat{c}_{ij} \leq 0$ برای کلیه‌ی متغیرها برقرار باشد (واضح است که $z_{ij} - \hat{c}_{ij} = 0$ برای متغیرهای اساسی) به خصوص در جواب بهینه بایستی شرایط زیر برقرار باشد.

(پ-۱۱)     برای $\lambda_{ij}$ غیر اساسی     $z_{ij} - \hat{c}_{ij} = wA_i x_{ij} + \alpha_i - c_i x_{ij} \geq 0$

(پ-۱۲)     برای $\mu_{ij}$ غیراساسی     $z_{ij} - \hat{c}_{ij} = wA_i d_{ij} - c_i d_{ij} \geq 0$

که در آن شرایط (پ-۱۱) و (پ-۱۲) برقرار است یا نه را می‌توان به وسیله‌ی حل مسأله‌ی فرعی زیر بررسی نمود.

Max $(wA_i - c_i) x_i + \alpha_i$
   S.T.
     $x_i \in X_i$

اگر جواب بهینه نامحدود باشد، آنگاه جواب همگن رأسی $d_{ik}$ یافت می‌شود، به طوری که $(wA_i - c_i) d_{ik} > 0$، یعنی شرط (پ-۱۲) برقرار نمی‌باشد. در این حالت متغیر غیراساسی $\mu_{ik}$ وارد پایه می‌شود که باعث بهتر شدن تابع مقصود می‌گردد، زیرا $z_{ik} - \hat{c}_{ik} > 0$ اگر جواب بهینه محدود باشد، آنگاه برای تمام مسائل فرعی شرایط (پ-۱۲) برقرار خواهد بود. فرض کنیم $x_{ik}$، نقطه رأسی بهینه باشد. اگر مقدار بهینه $z_{ik} - \hat{c}_{ik} = wA_i x_{ik} + \alpha_i - c_i x_{ik} \leq 0$ آنگاه شرایط (پ-۱۱) برای مسأله‌ی فرعی $i$ برقرار می‌باشد. در غیر این صورت $\lambda_{ik}$ وارد پایه می‌گردد. اگر برای تمام مسائل فرعی $z_{ik} - \hat{c}_{ik} \leq 0$، آنگاه جواب بهینه برای مسأله‌ی اصلی به دست آمده است. اگر قیدهای مادر از نوع نامساوی باشند آنگاه $z_j - \hat{c}_j$ برای متغیرهای زاید غیراساسی بایستی در هر مرحله محاسبه گردد.

به طور خلاصه، هر مسأله‌ی فرعی $i$ به نوبه خود حل می‌شود. اگر مسأله‌ی فرعی $i$ جواب نامحدود به دست بدهد، آنگاه یک جواب همگن رأسی $d_{ik}$ یافت می‌شود، آنگاه $\mu_{ik}$ وارد پایه می‌گردد. اگر مسأله‌ی فرعی $i$ جواب بهینه محدود داشته باشد و $wA_i x_{ik} + \alpha_i - c_{ik} x_{ik} > 0$، آنگاه، متغیر $\lambda_{ik}$ کاندید وارد شدن به پایه می‌باشد. اگر هیچ یک از شرایط فوق الذکر برقرار نباشد. در این صورت هیچ یک از متغیرهای مسأله‌ی فرعی $i$ وارد پایه در مسأله‌ی اصلی نمی‌گردد. اگر هیچ یک از مسائل فرعی، متغیری جهت کاندید شدن برای ورود به پایه در مسأله‌ی اصلی نداشته باشند آنگاه، برای مسأله‌ی اصلی جواب بهینه به دست آمده است. در غیر این صورت مثبت‌ترین $z_{ik} - \hat{c}_{ik}$ را در نظر می‌گیریم (یا اولین $z_{ik} - \hat{c}_{ik}$ ممکن است در نظر گرفته شود) متغیر متناظر وارد پایه می‌شود و مراحل تکرار می‌گردند.

**محاسبه‌ی کرانهای پائین برای مسائل فرعی که محدود هستند**

فرض می‌کنیم $x_1, x_2, ..., x_T$ جوابهای شدنی مسأله‌ی اصلی باشند، به طوری که $x_i \in X_i$ و $\sum_{i=1}^{T} A_i x_i = b$ با توجه به تعریف $z_{ik} - \hat{c}_{ik}$ داریم.

$$(wA_i - c_i) x_i + \alpha_i \leq (z_{ik} - \hat{c}_{ik})$$

یا

$$c_i x_i \geq wA_i x_i + \alpha_i - (z_{ik} - \hat{c}_{ik})$$

اگر عبارت فوق را برای هر $i$ محاسبه نموده و جمع نماییم، داریم:

$$\sum_{i=1}^{T} c_i x_i \geq w \sum_{i=1}^{T} A_i x_i + \sum_{i=1}^{T} \alpha_i - \sum_{i=1}^{T} (z_{ik} - \overline{c}_{ik})$$

امّا:

$$\sum_{i=1}^{T} A_i x_i = b \quad \text{و} \quad \sum_{i=1}^{T} c_i x_i = cx$$

بنابراین:

$$cx \geq wb + \alpha_1 - \sum_{i=1}^{T} (z_{ik} - \hat{c}_{ik})$$

یا

$$cx \geq \hat{c}_B \overline{b} - \sum_{i=1}^{T} (z_{ik} - \hat{c}_{ik})$$

که تعمیم یافته کران پائین مسأله‌ی فرعی در قسمت گذشته می‌باشد.

$$(c_i - wA_i) x_{ik} - \alpha_i \geq 0$$

**پیوست ب**

**مدل دوگان**

برای معنی داشتن مفاهیم جبر مجرد، چون تعامد و Polarity بایستی در چارچوبی عمل کنیم که ضرب داخلی[1] معنی داشته باشد.

- یک حلقه‌ی جابجایی بنیادین[2] تعریف می‌کنیم: R؛ که در تمام کاربردهای مورد نظر در این مقاله مجموعه اعداد حقیقی با دو عمل معمولی جمع و ضرب خواهد بود.

همه‌ی مجموعه‌ها مورد نظر فرض می‌شود که زیرمجموعه‌های $\mathbb{R}^n$، یعنی R- مدل آزاد از درجه‌ی n[3] خواهند بود.

تعریف ۱-۱- فرض کنیم R یک حلقه‌ی جابجایی است که $X \subseteq R ; \{0,1\} \subseteq X$؛ $D \subseteq X$؛ برای هر زیرمجموعه‌ی $S \subseteq X^n$؛

«$(X,D)$ دوگان یا $(X,D)$ Dual» مجموعه‌ی S را به صورت زیر تعریف می‌کنیم:

$S^* = \{x \in X^n : Sx \in D\}$

در جایی که $Sx \in D$ به معنی این است که $\forall s \in S ; sx \in D$

$(X,D) \text{ Dual } S = S^* = \{x \in X^n ; \forall s \in S ; s.x \in D\}$

هنگامی‌که:

$R = \mathbb{R} , X = \mathbb{R} , D = \{0\} = 1$

$S^* = (X,D) \text{ Dual of subset } S \subseteq \mathbb{R}^n = \{x \in \mathbb{R}^n ; \forall s \in S ; sx = 0\}$

$\Rightarrow S^* = S^\perp = S$ زیر فضای متعامد بر =

Orthogonal Complement of a Vector Subspace of $\mathbb{R}^n$

$\underbrace{\text{مکمل متعامد یک زیرفضای برداری } \mathbb{R}^n}_{S \subseteq \mathbb{R}^n}$

---

1- Inner Product
2- Underlying Commutative Ring R
3- R-module of rank n

هنگامی‌که

$R = \mathbb{R}, X = \mathbb{R}, D = \mathbb{R}_+ \Rightarrow$

$S^* = (X,D)$ دوگان زیرمجموعه $\subseteq \mathbb{R}^n = \{x \in \mathbb{R}^n, \forall s \in S; Sx \in \mathbb{R}_+\}$

$\Rightarrow S^* = S^\perp =$ Cone Polar to S

از این به بعد فرض می‌کنیم X، D و R ثابت هستند و به سادگی تنها به $S^*$ به صورت «دوگان S» ارجاع خواهیم داد.

تعریف ۲-۱- برای یک زیرمجموعه‌ی $S \subseteq X^n$، اگر:

S محدود شده است $\Rightarrow S = T^*$ برای بعضی $T \subseteq X^n$

علاوه بر آن اگر:

S به طور متناهی محدود شده است $\Rightarrow S = T^*; |T|$ متناهی باشد

نتیجه:

مجموعه‌های متناهیاً محدود شده برای $A_{m \times n} \in X^{m \times n}$ به شکل $\{x \in X^n : Ax \in D^m\}$ هستند.

تعریف ۳-۱- زیرمجموعه $S \subseteq X^n$ تحت دوگان $(X,D)$ بسته است اگر $\underbrace{S = S^{**}}_{S=(S^*)^*}$

نتایج حاصل از تعریف ۳-۱:

1-4- Proposition
(a) $S \subseteq T \Rightarrow T^* \subseteq S^*$
(b) $(S \cup T)^* = S^* \cap T^*$
(c) $S \subseteq S^{**}$
(d) $S^* = S^{***}$
(e) S is closed $\Leftrightarrow$ S is Constrained
(f) $S^{**}$ is the smallest closed set containings

برای تعریف یک مجموعه (مثلاً $S^*$) تاکنون از روشی استفاده کردیم که از بیرون مجموعه را تعریف می‌کرد. (یعنی از مجموعه‌ی S که ممکن است خارج از مجموعه‌ی $S^*$ باشد برای تعریف آن استفاده می‌شود.) همانطور می‌توان برای تعریف هر زیرفضا یا مخروط محدبی به طور معادل از تعریف مناسب دیگری از درون نیز استفاده کرد. مثلاً این بار از ترکیب خطی مجموعه‌ای از مولدها برای ایجاد یک زیرفضا یا ترکیب غیرمنفی برای ایجاد مخروط.

× نمایش داخلی[1] یا خارجی[2] یک مجموعه مشابه در واقع همان «Duality» است. در نتیجه ما به تعریف مناسبی از «Generated Sets» نیازمندیم که با استفاده از انواع مشخصی از ترکیبات خطی بیان شده باشد. توجه کنید که در مورد مخروط‌ها و زیرفضاها، انواعی از ترکیبات خطی در نظر گرفته می‌شوند که تحت آنها مجموعه‌های محدود بسته باشند.

(are those under which the constrained sets are closed)

$S \subseteq X^n$ is closed with respect to T-linear combinations, where $T \subseteq X^m$

در حالی که

$(y_1 a_1 + ... + y_m a_m) \in S \,;\, a_1, ..., a_m \in S \,;\, y \in T$

به عنوان مثال:

تحت ترکیب خطی closed است. $\quad , \quad S \subseteq \mathbb{R}^n \quad , \quad \{x \in \mathbb{R}^n : Sx = 0\}$

تحت ترکیب خطی غیرمنفی closed است. $\quad , \quad S \subseteq \mathbb{R}^n \quad , \quad \{x \in \mathbb{R}^n : Sx \geq 0\}$

$(D^m)^* = (X, D) - \text{dual of } D^m$

قضیه ۱-۵

برای هر $m \geq 1$، $(D^m)^*$ یک مجموعه از بردارهای $m$ مؤلفه‌ای است که ضرایب ترکیب خطی را به نحوی که همه‌ی مجموعه‌های محدود تحت آنها، closed باشند تأمین می‌کند.

$(D^m)^* = \{y \in X^m : \forall n \geq 1, S \subseteq X^n, a_1, ..., a_m \in S^* \Rightarrow (y_1 a_1 + ... + y_m a_m) \in S^*\}$

اثبات:

فرض کنید $y \in (D^m)^*$، $S \subseteq X^n$، $a_1, ..., a_m \in S^*$ در حالی که $m \geq 1$: آنگاه برای هر $s \in S$ خواهیم داشت:

$(y_1 a_1 + ... + y_m a_m) s = y_1(a_1 s) + ... + y_m(a_m s) = yz$

با توجه به اینکه $z_i = a_i s$ برای $1 \leq i \leq m$ و از آنجایی که $a_i \in S^*$ و $s \in S \Rightarrow \underbrace{1 \leq i \leq m, z_i \in D}_{z \in D^m}$

$\Rightarrow$ با توجه به اینکه $\Rightarrow yz \in D \Rightarrow (y_1 a_1 + ... + y_m a_m) s \in D, \forall s \in S$

$y \in (D^m)^* \qquad\qquad (y_1 a_1 + ... + y_m a_m) \in S^*$

---

1- Internal
2- External

۱۰۰

علاوه بر این بنابر تعریف، $0, 1 \in X$ ⇐ $e_i \in X^m$ ; $e_i = (0, 0, ..., 1, 0, 0)$ ; $1 \leq i \leq m$

- همینطور بدیهی است که $e_i \in (D^m)^*$؛

در نتیجه برای $y \in X^m$، اگر شروط قضیه را با استفاده از $S = D^m, m = n$ و $e_i = a_i$ ; $1 \leq i \leq m$ تأمین نماییم به $y \in (D^m)^*$ می‌رسیم و بدین ترتیب اثبات پایان می‌پذیرد. □

این بدان معنی است که با توجه به $D \subseteq X \subseteq R$، هر گاه از عناصر $(D^m)^*$ به عنوان ضرایب ترکیب خطی برای عناصر $S^*$ (نماینده‌ی دوگان هر $S$ ای در هر فضای $X^n$ با هر $n \geq 1$ ای) استفاده کنیم، حاصل ترکیب خطی مجدداً عضوی از خود $S^*$ خواهد بود. همینطور برعکس اگر مجموعه‌ای را از $y \in X^m$ در نظر بگیریم که به ازای هر $n$ و هر $S \subseteq X^n$ و هر $m$ عضوی از $S^*$ ترکیب خطی آن عضو $S^*$ باشد، اعضای مجموعه‌ی مزبور عضو $(D^m)^*$ خواهند بود.

تعریف ۶-۱

برای $S \subseteq X^n$، مجموعه ساخته شده توسط $S$ را با $\sigma(s)$ نمایش می‌دهیم، بطوریکه:

$\sigma(s) = \{x \in X^n : x = y_1 a_1 + ... + y_m a_m ,\text{ where } m \geq 1, a_1, ..., a_m \in s \text{ and } y \in (D^m)^*\}$

خود $S$ هم می‌تواند یک مجموعه Generated باشد، مشروط بر اینکه برای $T \subseteq X^n$، $S = \sigma(T)$، اگر finitc.lTl باشد، آنگاه

$S$ به طور متناهی ساخته شده است.

× از آنجایی که $(D^m)^*$ شامل بردارهای واحد است ⇐ $S \subseteq \sigma(s)$
× روشن است که: $S \subseteq T \Rightarrow \sigma(s) \subseteq \sigma(T)$

می‌توان نتیجه گرفت که مجموعه‌هایی به طور متناهی ساخته شده را می‌شود بصورت زیر نمایش داد:

$S = \{yA : y \in (D^m)^*\}$ ; where $A \in X^{m \times n}$

قضیه ۷-۱- فرض کنیم $S \subseteq X^n$، گزاره‌های زیر معادل همدیگرند.

(a)  $\sigma(s) \subseteq S$
(b)  $\sigma(s) = S$
(c)  $S$ یک مجموعه ساخته شده است.

نتیجه‌ی ۸-۱: اگر $S \subseteq X^n$ باشد، آنگاه $\sigma(s)$ کوچکترین مجموعه ساخته شده شامل $S$ خواهد بود.
اثبات: فرض کنیم $T$ کوچکترین مجموعه ساخته شده شامل $S$ باشد:

۱۰۱

$$S \subseteq T \implies \sigma(s) \subseteq \underbrace{\sigma(T)}_{\sigma(T)=T}$$
$$\implies \sigma(s) \subseteq T$$

که با توجه به کوچکترین بودن T تنها حالت مساوی برقرار خواهد بود.

$$\implies \sigma(s) = T$$

نتیجه‌ی ۹-۱: هر مجموعه‌ی محدود شده یک مجموعه‌ی تولید شده است.

این نتیجه بسیار شبیه تئوری Minkowski برای مخروط‌های متناهی است: هر مخروط چند وجهی[1] به صورت متناهی تولید شده[2] است. البته نتیجه‌ی ما مشروط به متناهی بودن نیست.

۱۰-۱: عکس نتیجه‌ی بالا در حالت کلی صادق نیست.

قضیه ۱۱-۱- فرض کنیم که $S \subseteq X^n$، در نتیجه: $\sigma(S^*) = S^* = (\sigma(s))^*$

اثبات: می‌دانیم که $S^* = \sigma(S^*)$ ، $S \subseteq \sigma(s)$ ، $S \subseteq S^{**}$

چون مجموعه محدود است $\Leftarrow$ مجموعه ساخته شده هم هست

بنابر نتیجه‌ی ۸-۱ $\sigma(s)$ کوچکترین مجموعه ساخته شده شامل S است $\Leftarrow$ $S \subseteq \sigma(s) \subseteq S^{**}$

بنابر ۴-۱ بند a: $\quad S^* = (S^{***} \subseteq (\sigma(s))^* \subseteq S^*)$

$$\implies S^* = (\sigma(s))^*$$
$$\implies S^* = (\sigma(s))^* = \sigma(s^*)$$

نتیجه ۱۲-۱: اگر $\sigma(T) = S$ $\Leftarrow$ $S^* = T^*$

اثبات: $\sigma(T^*) = S^* \implies (\sigma(T))^* = S^*$

$$\implies T^* = S^*$$

نتیجه‌ی ۱۳-۱: فرض کنیم $A \in X^{m \times n}$ و تعریف کنیم $T = \{yA : y \in (D^m)^*\}$، آنگاه $T^* = \{x \in X^n : Ax \in D^m\}$ در نتیجه هر مجموعه ایجاد شده متناهی یک دوگان به طور متناهی محدود شده دارد.

---

1- finitely constrained
2- finitely generated